\documentclass[3p,times]{elsarticle}

\usepackage{graphicx}
\usepackage{psfrag}
\usepackage[dvipsnames]{xcolor}

\usepackage{amsthm}
\usepackage{dsfont}
\usepackage{rotating}
\usepackage{multirow}
\usepackage{epstopdf}
\usepackage{newtxmath}

\DeclareMathAlphabet\mathbfcal{OMS}{cmsy}{b}{n}  

\newcommand{\R}{\mathds{R}}
\newcommand{\DD}{\mathbf{D}}
\newcommand{\QQ}{\mathbf{Q}}
\renewcommand{\vv}{\mathbf{v}}
\newcommand{\xx}{\mathbf{x}}
\newcommand{\nn}{\mathbf{n}}
\newcommand{\q}{\mathbf{q}}

\newcommand{\xxi}{\boldsymbol{\xi}}
\newcommand{\dt}{{\Delta t}}

\newcommand{\Q}{\mathbf{Q}}

\newcommand{\U}{\mathbf{U}}

\newcommand{\Fr}{\textnormal{Fr}}
\newcommand{\Sr}{\textnormal{Sr}}
\newcommand{\CFL}{\textnormal{CFL}}

\newcommand{\bzero}{\mathbf{0}}
\newcommand{\w}{\mathbf{w}}

\newcommand{\Topp}{\vmathbb{T}} 
\newcommand{\Vop}{\vmathbb{V}} 
\newcommand{\Sop}{\vmathbb{S}} 
\newcommand{\Rop}{\vmathbb{R}} 

\newtheorem{theorem}{Theorem}

\DeclareMathOperator{\argmin}{argmin}

\begin{document}
	
	\begin{frontmatter}
		
		\journal{Applied Numerical Mathematics}
		
		

		\title{An all Froude high order IMEX scheme for the shallow water equations on unstructured Voronoi meshes}
		
		\author[UNIFE]{Walter Boscheri$^*$}
		\ead{walter.boscheri@unife.it}
		\cortext[cor1]{Corresponding author}
		\author[UNIBZ]{Maurizio Tavelli}
		\ead{maurizio.tavelli@unibz.it}
		\author[UNITAR]{Cristóbal E. Castro}
		\ead{ccastro@academicos.uta.cl}
		\address[UNIFE]{Department of Mathematics and Computer Science, University of Ferrara, Ferrara, Italy}
		\address[UNIBZ]{Faculty of Computer Science, Free University of Bozen, Bozen, Italy}
		\address[UNITAR]{Departamento de Ingeniería Mecánica, Facultad de Ingeniería, Universidad de Tarapacá, Arica, Chile}
%
\begin{abstract}
We propose a novel numerical method for the solution of the shallow water equations in different regimes of the Froude number making use of general polygonal meshes. The fluxes of the governing equations are split such that advection and acoustic-gravity sub-systems are derived, hence separating slow and fast phenomena. This splitting allows the nonlinear convective fluxes to be discretized explicitly in time, while retaining an implicit time marching for the acoustic-gravity terms. Consequently, the novel schemes are particularly well suited in the low Froude limit of the model, since no numerical viscosity is added in the implicit solver. Besides, stability follows from a milder CFL condition which is based only on the advection speed and not on the celerity. High order time accuracy is achieved using the family of semi-implicit IMEX Runge-Kutta schemes, while high order in space is granted relying on two discretizations: (i) a cell-centered finite volume (FV) scheme for the nonlinear convective contribution on the polygonal cells; (ii) a staggered discontinuous Galerkin (DG) scheme for the solution of the linear system associated to the implicit discretization of the pressure sub-system. Therefore, three different meshes are used, namely a polygonal Voronoi mesh, a triangular subgrid and a staggered quadrilateral subgrid. The novel schemes are proved to be Asymptotic Preserving (AP), hence a consistent discretization of the limit model is retrieved for vanishing Froude numbers, which is the given by the so-called "lake at rest" equations. Furthermore, the novel methods are well-balanced by construction, and this property is also demonstrated. Accuracy and robustness are then validated against a set of benchmark test cases with Froude numbers ranging in the interval $\Fr \approx [10^{-6};5]$, hence showing that multiple time scales can be handled by the novel methods. 
\end{abstract}
%
\end{frontmatter}

\section{Introduction} \label{sec.intro}

Shallow water equations are extensively used in modeling physical processes that affect environmental and geophysical phenomena \cite{Brocchini2008,GarciaNavarro2019,FrankGiraldo2020,ToroBookSWE,Zeitlin2018,Mandli2021}. They are designed to describe the dynamics of shallow incompressible and inviscid fluid flows. Scenarios of storm surges, tsunami wave propagation-inundation, dam breaks and river floods as well as atmospheric processes are some of the most challenging cases. All these physical applications involve the description of multiscale phenomena where advection and acoustic-gravity waves coexist. The Froude number $\Fr=u/\sqrt{g h}$, which measures the ratio between convective velocity and pressure wave speed, is typically used to represent the time scale of the flow under consideration. In order to properly solve the shallow water equations, numerical schemes are designed to accurately and efficiently compute the solution of the hyperbolic governing equations \cite{Beisiegel2021}. However, due to the multiscale nature of the physical process, this goal is still a very challenging problem \cite{Vater2018}. For example, in a tsunami case, the advection process describes the motion of a fluid parcel at small finite velocity while the acoustic-gravity wave travels proportionally to the square root of gravity times water depth which can be in the order of 100-800 kilometers an hour in open sea. 

Therefore, to deal with multiple time scales, numerical methods must be designed and constructed with the so-called Asymptotic Preserving (AP) property, meaning that the numerical schemes can capture the behavior of the governing equations in the asymptotic limits of the model, i.e. when $\Fr \to 0$ \cite{JINAP1,JinPar}. Explicit Godunov-type solvers \cite{Godunov1959, HLL1983, LW1960, Munz1994} are very popular and behave quite well for high Froude number flows. Shock-capturing schemes for the shallow water equations have been reviewed in \cite{ToroBookSWE}. In the low Froude regime, explicit schemes are inaccurate and do not satisfy the AP property because of the severe CFL-type stability condition which is based on the acoustic-gravity wave speed that becomes dominant. Indeed, in \cite{Dellacherie1} the effect of numerical viscosity on the slow waves introduced by upwind-type schemes is proven to degrade the accuracy. As a consequence, in low Froude regimes, explicit numerical methods are forced to perform a huge number of small time steps to keep tracking of acoustic waves while the fluid barely moves. On the other hand, when fully implicit time step schemes are used, larger time steps can be employed without losing significant information and preserving the quality of the numerical solution. The major drawback of fully implicit time stepping techniques is given by the need of solving a system which might become strongly nonlinear due to the presence of the convective fluxes in the governing equations. 

To overcome this problem, a class of semi-implicit discretizations has started to gain visibility in the last decades \cite{ParkMunz2005,Klein,Casulli1990,Casulli1999,BosFil2016,Boscarino22,ADERFSE,DegTan,Chalons}. In this context, advection is discretized explicitly while pressure is taken implicitly, thus the resulting stability condition is only constrained by the main flow speed, which in the low Froude limit vanishes. Thus, semi-implicit schemes are much more efficient compared to explicit methods in the low Froude limit, and they also exhibit less numerical viscosity and accurate resolution because the implicit terms do not need any numerical stabilization. Furthermore, in the semi-implicit context, the resulting system for the unknown pressure typically results to be linear or mildly nonlinear \cite{Dumbser_Casulli16,BosPar2021}, hence avoiding the solution of strongly nonlinear algebraic systems. The idea of separating the slow and fast time scale has been effectively interpreted as a splitting of the fluxes, see \cite{Toro_Vazquez12}. Following this approach, a lot of research has been carried out to devise numerical methods able to deal with multiple time scales \cite{StagDG_Dumbser2013,BDLTV2020,Busto_SWE2022,Ioriatti_SWE2019,ChertockKurganov19,Bermudez2020,SIINS22}, which are often referred to as \textit{all Mach} solvers, recalling the hydrodynamic analogy. Another strategy to deal with multiple time scales is given by the class of implicit-explicit (IMEX) methods \cite{AscRuuSpi,BP2017,BosRus,PR_IMEX} or, more in general, by the so-called partitioned schemes \cite{Hofer}. In \cite{BosFil2016}, semi-implicit and IMEX time stepping techniques have been unified in a single framework leading to efficient all Mach solvers \cite{BosPar2021,SICNS22,Boscarino22} with linearly implicit algebraic systems. 

An important aspect regarding shallow water models for variable bottom topography is the balance between fluxes and the geometrical source, where extensive research has been developed. If the numerical scheme preserves this balance it is called \textit{well-balanced} or preserving the C-property \cite{BermudezWB,WBLeVeque,Pares2021}. Without this preserving equilibrium capability a numerical scheme is not useful as for example in tsunami propagation scenarios \cite{Castro2012}. A rather general approach to design well-balanced algorithms relies in the family of path-conservative schemes \cite{CastroPares2006,CastroPares2013}, that was originally proposed for the definition of weak solutions in the context of non-conservative hyperbolic systems \cite{CastroPares2008,CastroPares2010}.

If the physical process that needs to be modeled and simulated is related to long distance and time wave propagation, a very accurate numerical method is mandatory in order to preserve the information that is propagated by the numerical scheme. Such methods need to comply with space and time accuracy in the form of high order spatial and time discretization of the numerical solution \cite{Castro2008,Shu2014,Toro2020}. High order semi-implicit discontinuous Galerkin (DG) schemes for the shallow water equations have been recently forwarded in \cite{StagDG_Dumbser2013,TavelliSWE2014}, while high order IMEX finite volume schemes for hydrodynamics can be found for instance in \cite{BosPar2021,SICNS22}.

The aim of this work is to design a high order implicit-explicit scheme for the shallow water equations that can capture the flow behavior at all Froude numbers while respecting the well-balanced property. To achieve high order of accuracy in space, a CWENO reconstruction technique is employed on general polygonal grid, differently from what has been recently presented in \cite{WENOFD_SWE22} where a WENO finite difference scheme was designed on Cartesian meshes. Moreover, in our approach a robust finite volume method is used which can deal with very general control volumes. An asymptotic preserving scheme for the shallow water system with Coriolis forces has been derived in \cite{ChertockKurganov19}, which applies to low Froude flows and achieves up to second order of accuracy. Here, we will show higher accurate time stepping techniques based on the usage of semi-implicit IMEX schemes. Instead of using a finite element paradigm for the solution of the implicit part of the governing system \cite{Busto_SWE2022}, we design a discontinuous Galerkin solver applied to a staggered triangular subgrid. In this way, high order spatial accuracy can be easily achieved on polygonal grids using a compact stencil. Data are transferred between different meshes and discretizations by means of high order $L_2$-projection operators that will be specifically designed. The new algorithms do not require any orthogonality property of the computational mesh, therefore they can be applied to any unstructured conforming mesh. 

This article is organized as follows. In Section \ref{sec.pde} the governing two-dimensional shallow water equations are presented, studying the multiscale nature of the equation by deriving the dimensionless form, followed by splitting the system into an advection and pressure sub-system and studying the asymptotic behavior in the low Froude limit. In Section \ref{sec.numscheme} the new numerical scheme is presented considering the discretization in both space and time. Section \ref{sec.numtest} is devoted to show a set of numerical results which demonstrate the accuracy and robustness of the novel scheme considering different flow regimes with Froude number ranging from $\Fr=10^{-6}$ up to $\Fr=5.73$. Finally, Section \ref{sec.concl} finalizes this article by summarizing the work and giving an outlook to future investigations.

\section{Governing equations} \label{sec.pde}
\begin{figure}[!htbp]
	\begin{center}
		\begin{tabular}{c}		\includegraphics[width=0.6\textwidth]{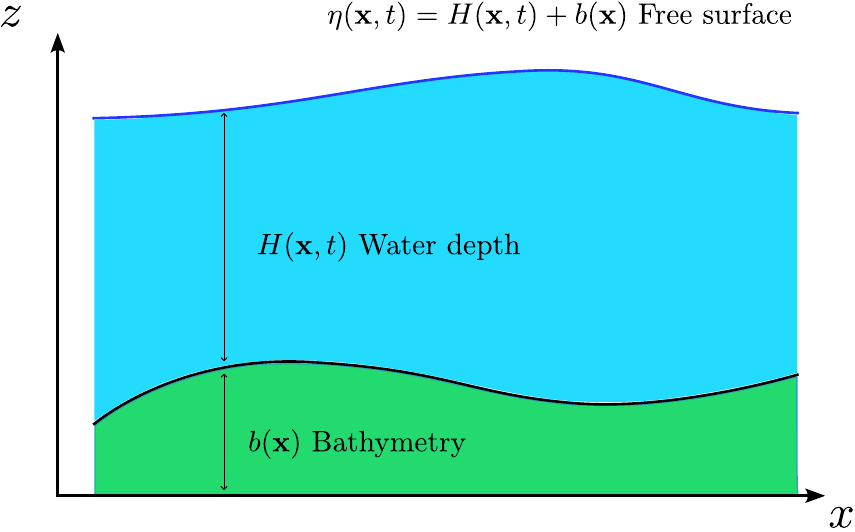} \\
		\end{tabular}
		\caption{Schematic of the computational domain and the notation used for the shallow water equations over a fixed bottom.}
		\label{fig.SWEnotation}
	\end{center}
\end{figure}

Let us consider a two-dimensional bounded domain $\Omega \in \R^2$, which is defined by the space coordinates $\xx=(x,y)$, and a time interval with the time coordinate $t \in \R_0^+$. The frictionless shallow water system is described by the following set of partial differential equations (PDE):
\begin{subequations}
	\label{eqn.SWE}
	\begin{align}
	\frac{\partial \eta}{\partial t} + \nabla \cdot \q &= 0, \label{eqn.SWE_eta} \\
	\frac{\partial \q}{\partial t} + \nabla \cdot \left( \vv \otimes \q \right) + gH \nabla \eta &= \bzero, \label{eqn.SWE_q} \\
	\frac{\partial b}{\partial t} &= 0, \label{eqn.SWE_b} 
	\end{align}
\end{subequations}
where $\eta(\xx,t) \geq 0$ is the free surface elevation, $b(\xx)$ is a prescribed bottom bathymetry, $H(\xx,t)=\eta(\xx,t)-b(\xx)\geq 0$ represents the total water depth and $g$ is gravity acceleration. The velocity of the water is described by the vector field $\vv(\xx,t)= (u,v)$, while $\q(\xx,t)=H(\xx,t) \vv(\xx,t)$ denotes the flow discharge, which is the corresponding conservative variable. A schematic of the domain and the notation of the governing PDE is shown in Figure \ref{fig.SWEnotation}.

\subsection{Scaling of the shallow water equations}  \label{ssec.scaling}
To analyze the multiscale nature of the governing equations, it is convenient to derive the corresponding dimensionless form \cite{ChertockKurganov19,BDLTV2020,BosPar2021}. Therefore, the governing PDE \eqref{eqn.SWE} can be rescaled and represented in dimensionless form by considering the rescaled variables:
\begin{equation}
	\label{eqn.scale_var}
	\tilde{\xx}=\xx/L_0, \qquad \tilde{t}=t/T_0, \qquad \tilde{\eta}=\eta/H_0, \qquad \tilde{H}=H/H_0, \qquad \tilde{\vv} = \vv / U_0, \qquad \tilde{b}=b/H_0,
\end{equation}
where $L_0$, $T_0$, $H_0$, $U_0$ are the characteristic length, time, depth and velocity, respectively. Substitution of \eqref{eqn.scale_var} into system \eqref{eqn.SWE}, dropping the tilde superscripts to ease notation, yields the rescaled shallow water equations:
\begin{subequations}
	\label{eqn.SWEr}
	\begin{align}
		\Sr \cdot \frac{\partial \eta}{\partial t} + \nabla \cdot \q &=0, \label{eqn.SWEr_eta} \\
		\Sr \cdot \frac{\partial \q}{\partial t} + \nabla \cdot \left( \vv \otimes \q \right) + \frac{H}{\Fr^2} \nabla \eta &= \bzero, \label{eqn.SWEr_q} \\
		\frac{\partial b}{\partial t} &= 0, \label{eqn.SWEr_b}
	\end{align}
\end{subequations}
with the Strouhal number $\Sr$ and the Froude number $\Fr$ defined as
\begin{equation}
	\label{eqn.Sr_Fr}
	\Sr = \frac{L_0}{T_0 \, U_0}, \qquad \Fr = \frac{U_0}{\sqrt{g \, H_0}}.
\end{equation} 
Without loss of generality, we assume that the characteristic time results to be $T_0=L_0/U_0$, hence leading to a convective time scale with $\Sr=1$. Furthermore, a reference Froude number $\Fr=\varepsilon$ is considered, with $\varepsilon$ representing the asymptotic expansion parameter which will be used to study the asymptotic limit of the governing equations. Consequently, the rescaled system \eqref{eqn.SWEr} can be rewritten as
\begin{subequations}
	\label{eqn.SWErr}
	\begin{align}
		\frac{\partial \eta}{\partial t} + \nabla \cdot \q &=0, \label{eqn.SWErr_eta} \\
		\frac{\partial \q}{\partial t} + \nabla \cdot \left( \vv \otimes \q \right) + \frac{H}{\varepsilon^2} \nabla \eta &= \bzero, \label{eqn.SWErr_q} \\
		\frac{\partial b}{\partial t} &= 0. \label{eqn.SWErr_b}
	\end{align}
\end{subequations}  
The rescaled shallow water system is hyperbolic and its eigenvalues $\boldsymbol{\lambda}=(\lambda_1,\lambda_2,\lambda_3,\lambda_4)$ in the normal  direction $\nn=(n_x,n_y)$ are given by
\begin{equation}
	\label{eqn.lambda}
	\lambda_1 = \vv \cdot \nn - a/\varepsilon, \qquad \lambda_2 = 0, \qquad \lambda_3=\vv \cdot \nn, \qquad \lambda_4 = \vv \cdot \nn + a/\varepsilon,
\end{equation}
with the rescaled celerity $a=\sqrt{H}$ (the dimensional celerity is $a=\sqrt{g H}$). 

\subsection{Flux splitting} \label{ssec.fluxsplit}
We can make some considerations about the eigenstructure \eqref{eqn.lambda} of the rescaled system. The zero eigenvalue $\lambda_2$ is related to the bottom jump, while the transport of the transverse velocity corresponds to the eigenvalue $\lambda_3$. The eigenvalues $\lambda_1$ and $\lambda_4$ carry information about the propagation of acoustic-gravity waves and they are responsible of a severe time step restriction if fully explicit numerical schemes are used to discretize the governing PDE. Indeed, the stability condition on the time step $\dt=t^{n+1}-t^n$ is given by
\begin{equation}
	\label{eqn.timestep_full} 
	\dt \leq \CFL \min \limits_{\Omega} \frac{h}{|\vv \cdot \nn \pm a/\varepsilon|},
\end{equation}  
where $h$ represents the characteristic mesh size of the computational cell and the CFL number must be chosen such that $\CFL < 1/2$ on two-dimensional unstructured grids to ensure stability. Let us notice that the time step $\dt$ is of order $\varepsilon$, thus the time step goes to zero when $\varepsilon \to 0$, which is the so-called \textit{asymptotic limit} of the system. Apart from being extremely inefficient in the asymptotic limit due to vanishing time steps, explicit schemes are also not able to correctly capture the asymptotic regime as discussed in \cite{Guillard,GH,Dellacherie1}. 

Therefore, we proceed adopting a \textit{flux splitting} technique, that has been widely used in the literature for the shallow water equations \cite{Casulli1990,Casulli1992,Toro_Vazquez12,TavelliSWE2014,Busto_SWE2022} but also for the compressible Euler and Navier-Stokes equations \cite{Dumbser_Casulli16,TavelliCNS,BDLTV2020,BosPar2021} as well as for incompressible fluids \cite{TavelliIncNS,SIINS22}. Specifically, the shallow water system is divided into a convective and a pressure sub-system, which will be discretized explicitly and implicitly in time, respectively. The two sub-systems write as follows. 
\begin{itemize}
	\item Convective sub-system:
	\begin{equation}
    	\label{eqn.SWEc}		
		\left\{ \begin{array}{rcl}
			\partial_t \eta &=& 0 \\ [2mm]
			\partial_t \q + \nabla \cdot \left( \vv \otimes \q \right) &=& \bzero \\ [2mm]
			\partial_t b &=& 0
		\end{array} \right. , \qquad \boldsymbol{\lambda}^c=\left( 0, \, 0, \, \vv \cdot \nn, \, 2 \vv \cdot \nn \right).
	\end{equation}
\item Pressure sub-system
    \begin{equation}
    	\label{eqn.SWEp}
    	\left\{ \begin{array}{rcl}
    		\partial_t \eta + \nabla \cdot \q &=& 0 \\ [2mm]
    		\partial_t \q + \frac{H}{\varepsilon^2} \nabla \eta &=& \bzero \\ [2mm]
    		\partial_t b &=& 0
    	\end{array} \right. , \qquad \boldsymbol{\lambda}^p=\left( -a/\varepsilon, \, 0, \, 0, \, a/\varepsilon \right).
    \end{equation} 
\end{itemize}

It is clear that the eigenvalues of the convective sub-system \eqref{eqn.SWEc} do not contain the celerity $a$, which is indeed present in the eigenvalues of the pressure sub-system \eqref{eqn.SWEp}. However, since the pressure sub-system will be discretized implicitly, those terms will no longer appear in the stability condition \eqref{eqn.timestep_full}, making the resulting numerical method extremely efficient. Furthermore, the absence of the free surface wave speed $a$ in the convective eigenstructure leads to numerical schemes which are particularly well suited for applications in the asymptotic regime, i.e. when $\varepsilon \to 0$, because the numerical dissipation will be drastically reduced since it is only proportional to the (very low) water speed. For example, low Froude flows take place in tidal motions, or in river flows across flat lands, or even mud and debris floods produced by landslides when bed load transport is included in the model.

\subsection{Low Froude limit of the shallow water equations} \label{ssec.lowFroude}
To investigate the asymptotic limit of the PDE system \eqref{eqn.SWErr}, let us assume the computational domain $\Omega(\xx)$ to be assigned with periodic boundary conditions on $\partial \Omega$ and let us introduce the $k$-th order 
Chapman-Enskog expansion of a generic variable $\phi(\xx,t)$ in powers of the non-dimensional stiffness parameter $\varepsilon$, that reads
\begin{equation}
	\phi(\xx,t) = \phi_{(0)}(\xx,t) + \varepsilon \phi_{(1)}(\xx,t) + \varepsilon^2 \phi_{(2)}(\xx,t) + \ldots + \mathcal{O}(\varepsilon^k).
	\label{eqn.exp}
\end{equation}
Notice that the bottom elevation $b(\xx)$ is time-independent, therefore it is not affected by the asymptotic expansion and it only contributes to zeroth order terms in the definition of the total water depth, namely
\begin{equation}
	H(\xx,t) = \eta_{(0)}(\xx,t) - b(\xx) + \varepsilon \eta_{(1)}(\xx,t)  + \varepsilon^2 \eta_{(2)}(\xx,t) + \ldots + \mathcal{O}(\varepsilon^k).
\end{equation} 
Application of the expansion \eqref{eqn.exp} to the rescaled governing PDE \eqref{eqn.SWErr} and collection of the like powers of $\varepsilon$ yields the following $k$-th leading order equations for $k\in\{0,-1,-2\}$:
\begin{itemize}
	\item $\mathcal{O}(\varepsilon^{0})$
	\begin{subequations}
		\label{eqn.SWE_eps0}		
		\begin{align}
			\partial_t \eta_{(0)} + \nabla \cdot \left( \left(\eta_0-b\right) \, \vv_0 \right)&= 0, \\ 
			\partial_t \q_{(0)} + \nabla \cdot \left( \vv_{(0)} \otimes \q_{(0)} \right) + \eta_{(2)} \nabla \eta_{(0)} + \eta_{(1)} \nabla \eta_{(1)} + \left(\eta_{(0)}-b \right) \nabla \eta_{(2)} &= \bzero,
		\end{align}
	\end{subequations}

    \item $\mathcal{O}(\varepsilon^{-1})$
    \begin{equation}
    	\label{eqn.SWE_eps1}		
    		\eta_{(1)} \nabla \eta_{(0)} + \left(\eta_{(0)}-b \right) \nabla \eta_{(1)} = \bzero,
    \end{equation}

    \item $\mathcal{O}(\varepsilon^{-2})$
    \begin{equation}
    	\label{eqn.SWE_eps2}		
    		\left(\eta_{(0)}-b \right) \nabla \eta_{(0)} = \bzero.
    \end{equation}
\end{itemize} 
From \eqref{eqn.SWE_eps2} we immediately get that
\begin{equation}
	\eta_{(0)} \equiv \eta_{(0)}(t),
\end{equation}
hence the free surface elevation is constant in space. Using this information in \eqref{eqn.SWE_eps1} and assuming no dry area in the domain, namely $H(\xx,t)>0$, allows us to conclude that 
\begin{equation}
	\eta_{(1)} \equiv \eta_{(1)}(t),
\end{equation}
because the quantity $\left(\eta_{(0)}-b \right) \neq 0$ is constant in space as well. From \eqref{eqn.SWE_eps0} it follows that
\begin{subequations}
	\label{eqn.SWE_lowFr}		
	\begin{align}
		\partial_t \eta_{(0)}  + \nabla \cdot \left( \left(\eta_0-b\right) \, \vv_0 \right) &=0, \label{eqn.SWE_lowFr_eta}\\ 
		\partial_t \q_{(0)} + \nabla \cdot \left( \vv_{(0)} \otimes \q_{(0)} \right) + \left(\eta_{(0)}-b \right) \nabla \eta_{(2)} &= \bzero. \label{eqn.SWE_lowFr_q}
	\end{align}
\end{subequations}
Integration of the mass equation \eqref{eqn.SWE_lowFr_eta} over the computational domain and application of Gauss theorem leads to
\begin{equation}
	\label{eqn.SWE_lowFr_eta2}
	\partial_t \eta_{(0)} = -\frac{1}{|\Omega|} \int \limits_{\Omega} \nabla \cdot \left( \left(\eta_0-b\right) \, \vv_0 \right) \, d\xx = -\frac{1}{|\Omega|} \int \limits_{\partial \Omega} \left(\eta_0-b\right) \, \vv_0 \cdot \nn \, dS,
\end{equation} 
where $\nn$ is the outward pointing unit normal vector defined on the domain boundary $\partial \Omega$. Since we have assumed periodic boundaries, the right hand side of \eqref{eqn.SWE_lowFr_eta2} vanishes, implying that $\eta_{(0)}$ is constant both in space and time, hence the total water depth as well ($H(\xx,t)=const$). The \textit{low Froude shallow water system} then writes
\begin{subequations}
	\label{eqn.LowFr}		
	\begin{align}
		\nabla \cdot \left( \left(\eta_0-b\right) \, \vv_0 \right) &=0, \label{eqn.LowFr_eta}\\ 
		\partial_t \q_{(0)} + \nabla \cdot \left( \vv_{(0)} \otimes \q_{(0)} \right) + \left(\eta_{(0)}-b \right) \nabla \eta_{(2)} &= \bzero. \label{eqn.LowFr_q}
	\end{align}
\end{subequations}
%
\section{Numerical scheme} \label{sec.numscheme}

\subsection{Discretization of the space-time computational domain} \label{ssec.spacetimeD}
Let us fix some notation related to the space and time computational domains.

\paragraph{Time computational domain} The time coordinate is defined in the interval $[0;t_f]$, which is approximated by a sequence of discrete points $t^n$. Thus the time computational domain is discretized such that
\begin{equation}
	t = t^n + \alpha \dt, \qquad  \alpha \in [0,1], 
	\label{eqn:time}
\end{equation} 
with the time step $\dt=t^{n+1}-t^n$.
Because of the \textit{implicit} discretization of the pressure sub-system \eqref{eqn.SWEp}, the time step is limited by a classical CFL stability condition which is only based on the maximum convective eigenvalue \eqref{eqn.SWEc}, that is
\begin{equation}
	\label{eqn.timestep} 
	\dt \leq \CFL \min \limits_{\Omega} \frac{h}{|\vv \cdot \nn|},
\end{equation}
hence yielding a milder stability restriction compared to \eqref{eqn.timestep_full}, especially in the asymptotic regime when $\varepsilon \to 0$.

\paragraph{Space computational domain} The computational domain is discretized by a set of 
non-overlapping unstructured control volumes $P_i$ with boundary $\partial P_i$, that are given by arbitrary shaped polygons. We will use \textit{Voronoi tessellations} \cite{ArepoTN}, even though the orthogonality property of the grid is not necessary in our framework, differently from \cite{Casulli1990,Casulli1992,Voronoi,Voronoi-DivFree,ADERFSE}. The total number of cells is $N_P$, thus $i=1,\ldots,N_P$, and the union of all elements is called the tessellation
$\mathcal{D}_{\Omega}$ of the domain 
\begin{equation}
	\mathcal{D}_{\Omega} = \bigcup \limits_{i=1}^{N_P}{P_i}. 
	\label{trian}
\end{equation}	
The surface of each polygon $P_i$ is denoted with $|P_i|$, while $P_j$ represents the Neumann neighbor of $P_i$ which shares the edge $\partial P_{ij}$ of length $|\partial P_{ij}|$. The outward pointing normal vector on the edge $\partial P_{ij}$ is addressed with $\nn_{ij}$, and the characteristic mesh size of each element is measured by $h_i=\sqrt{|P_i|}$. The element $P_i$ counts a total number of vertexes $N_{v_i}$, which also corresponds to the total number of edges. The center of mass $\ \xx_{i}$ is computed as
\begin{equation}
\xx_{i} = \frac{1}{|P_i|} \int_{P_i} \xx \, d\xx.
\end{equation}	
Starting from the Voronoi tessellation, a \textit{triangular subgrid} is introduced, by connecting the center of mass with all the vertexes of each cell. Each sub-triangle of $P_i$ is labeled with $T_{m(ij)}$, meaning that it covers the area $|T_{m(ij)}|$ defined by the center of mass connecting the two vertexes of the edge $\partial P_{ij}$. The mono-index $m$ counts the global number of the subcells over the entire mesh, and it can be used to ease the notation, thus simply writing $T_m$ to address a generic subcell. Therefore, the computational domain is covered by a total number of $N_T= \sum \limits_{i=1}^{N_P} N_{v_i}$ sub-triangles with $m=1,\ldots,N_T$. The center of mass of each sub-triangle is consequently defined as 
\begin{equation}
	\xx_{m} = \frac{1}{|T_{m}|} \int_{T_{m}} \xx \, d\xx.
\end{equation}	
The total number of edges (without repetition) of the sub-triangulation is $N_e$, with $\Gamma_e$ denoting the $e$-th edge and $|\Gamma_e|$ its length. Let us define for each edge $\Gamma_e$ a standard normal vector $\nn_e$ that points from the arbitrarily chosen left $\ell(e)$ and right $r(e)$ sub-triangle sharing the common edge $\Gamma_e$. Therefore, a sub-triangle $T_{m}$ which shares the edge $\Gamma_e$ could be either the right or the left neighbor with respect to that edge. This sub-triangulation will be used to numerically integrate any quantity inside the cell $P_i$ as well as for the implicit discretization of the pressure sub-system.

\begin{figure}[!htbp]
	\begin{center}
		\begin{tabular}{cc}		\includegraphics[width=0.52\textwidth]{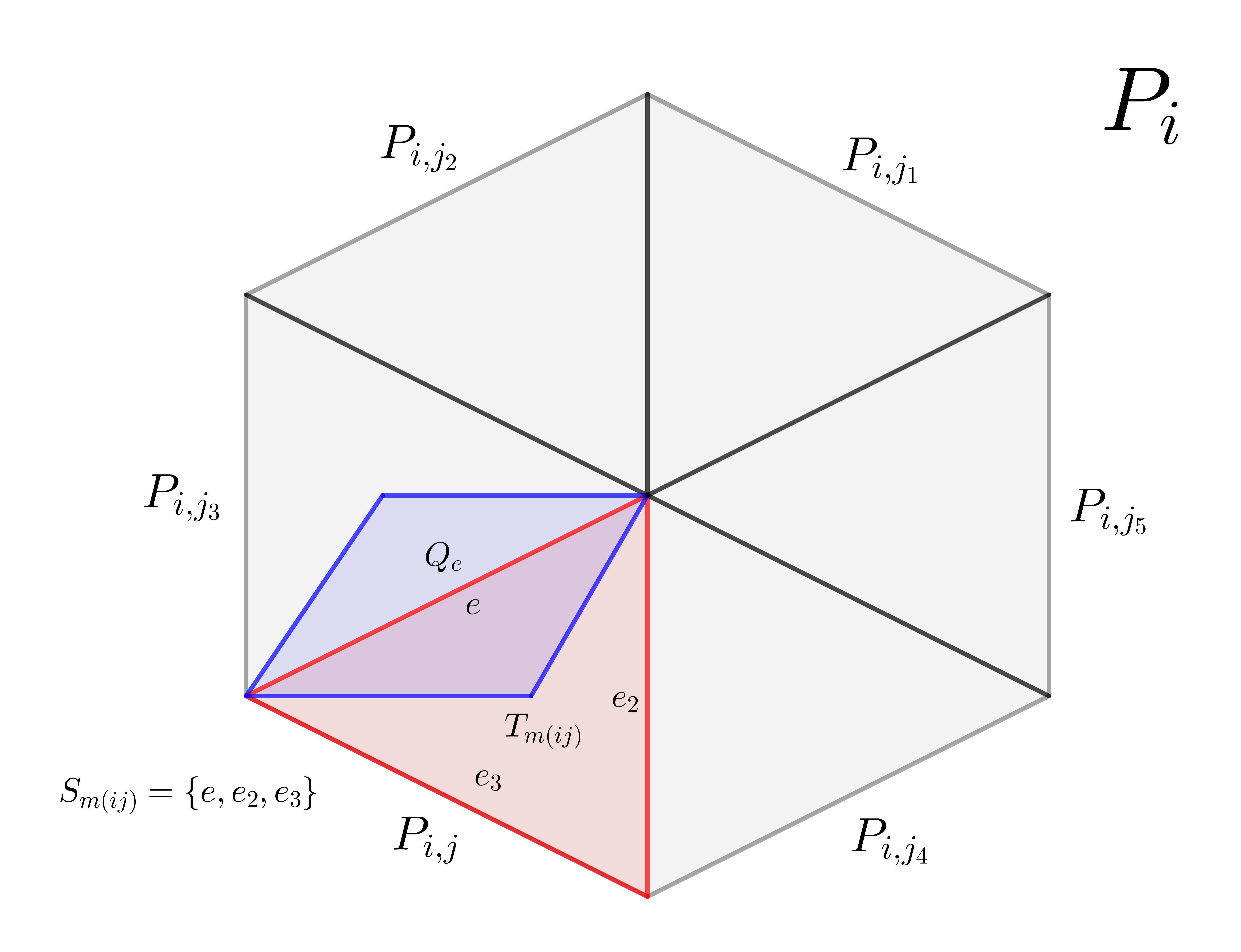}  &	\includegraphics[width=0.42\textwidth]{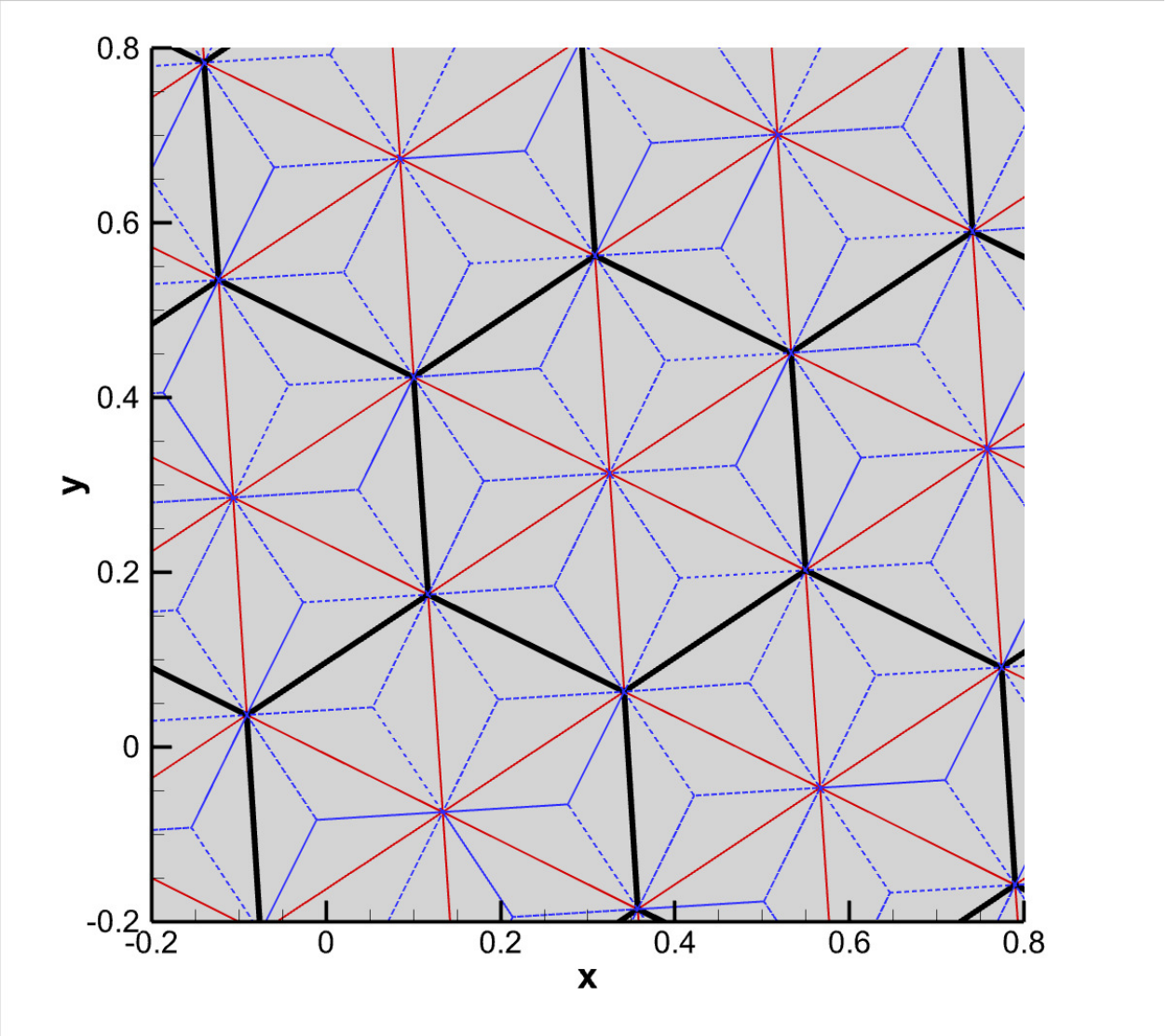}  \\
		\end{tabular}
		\caption{Left: mesh notation. Right: example of the computational meshes used for the numerical discretization of the shallow water equations: Voronoi tessellation (solid black lines), sub-triangulation (solid red lines) and staggered sub-triangulation (dashed blue lines).}
		\label{fig.mesh}
	\end{center}
\end{figure}

Finally, let us define a \textit{staggered quadrilateral subgrid} which is built upon the sub-triangulation. More precisely, the dual cell $Q_{e}$ is constructed by connecting the center of mass $\xx_{m_{\ell(e)}}$ to the two vertexes of the edge $\Gamma_e$ and then to the center of mass $\xx_{m_{r(e)}}$ of the neighbor sub-triangle. Obviously, the total number of staggered sub-elements is $N_e$. For every sub-triangle $m \in N_T$ we denote the set of edges of $T_m$ as $S_m$. A sketch of the used notation as well as the complete mesh is reported in Figure \ref{fig.mesh}.

\subsection{Basis functions and projection operators} \label{ssec.basefunc}
The solution of the shallow water equations will be numerically approximated as an expansion using a polynomial basis. We consider a polynomial space $\mathcal{V}_h$ up to degree $M$ with a total number of degrees of freedom $\mathcal{M}=\frac{1}{2}(M+1)(M+2)$. Because of the different nature of the three grids which are employed, namely the Voronoi tessellation, the triangular subgrid and the staggered quadrilateral subgrid, we need to define the associated basis functions as well as a set of projection operators to transfer data from one mesh to another.

On the Voronoi mesh, a set of \textit{conservative} Taylor functions are employed as \textit{modal} basis functions $\beta_l$, which are given by a truncated Taylor series of degree $M$ around the center of mass $\xx_{i}$ of the physical element $P_i$:
\begin{equation}
\label{eqn.Voronoi_modal}
\beta_l^{(i)}(\xx):= \beta_l(\xx)|_{P_i} = \frac{(x - x_{i})^{r_l}}{h_i^{r_l}} \, \frac{(y - y_{i})^{q_l}}{h_i^{q_l}} - \frac{1}{|P_i|}\int_{P_i}\frac{(x - x_{i})^{r_l}}{h_i^{r_l}} \, \frac{(y - y_{i})^{q_l}}{h_i^{q_l}} \, dx, \qquad 0 \leq {r_l}+{q_l}\leq M,
\end{equation}
where $l = 1, \dots, \mathcal{M}$ represents a mono-index that counts the total degrees of freedom of the expansion. The basis functions are also normalized by the characteristic length $h_i$ to avoid ill-conditioned approximations induced by low quality polygonal cells that might occur on unstructured meshes. We remark that the conservation property of the basis functions \eqref{eqn.Voronoi_modal} means that
\begin{equation}
	\label{eqn.modal_cons}
	\frac{1}{|P_i|}\int_{P_i} \sum_{l=1}^{\mathcal{M}}\beta_l{^{(i)}}(\xx) \, d\xx = 1,
\end{equation}
thus the first degree of freedom of each element $P_i$ (i.e. the one identified by $l=1$) represents the cell average value, in the finite volume sense. To make notation easier, the subscript referring to $P_i$ will be dropped, bearing in mind that the modal basis \eqref{eqn.Voronoi_modal} are defined in the physical space and therefore they are element-dependent, hence we compactly write $\beta_l^{(i)}$. The total number of the modal basis functions is then $N_{\beta}=\mathcal{M}$, thus $\{\beta_l\}_{l\in [1,N_{\beta}]}$.

The sub-triangulation allows the definition of \textit{nodal} basis functions. Each sub-triangle can be easily mapped to the reference triangular element $T_{std}$ in the reference coordinate system $\xxi=(\xi_1,\xi_2)$ defined as $T_{std} = \{ \xxi \in \R^2 : 0 \leq \xi_1 \leq 1, \, 0 \leq \xi_2 \leq 1-\xi_1  \}$. The transformation between physical $\xx$ and reference $\xxi$ coordinates of the subcell $T_{m}$ is given by the following linear mapping:
\begin{equation} 
	\label{eqn.xietaTransf} 
	\mathbf{x} := \xx(T_{m},\xxi) = \xx_{1,m} + (\xx_{2,m}-\xx_{1,m}) \, \xi_1 + (\xx_{3,m}-\xx_{1,m}) \, \xi_2,
\end{equation} 
with $\xx_{p,m}$ ($p=\{1,2,3\}$) being the vector of physical spatial coordinates of the $p$-th vertex of the sub-triangle with counterclockwise orientation. We denote with $\xxi(T_m,\xx):T_m \rightarrow T_{std}$ the inverse mapping. The coordinates of the nodes associated with the basis functions are defined on $T_{std}$ as
\begin{equation}
	\label{eqn.NCnodes}
	\xxi_k = \left( \frac{l_1}{M}, \frac{l_2}{M} \right), \qquad 0 \leq l_1 \leq M, \quad 0 \leq l_2 \leq M-l_1,
\end{equation}
with the multi-index $l=(l_1,l_2)=1,\ldots,\mathcal{M}$ already used in the definition of the modal basis \eqref{eqn.Voronoi_modal}. The nodal basis is then constructed by means of the Lagrange interpolation polynomials, hence imposing the interpolation condition
\begin{equation}
	\phi_k(\xxi_l) = \delta_{kl},
\end{equation}
where $\delta_{kl}$ denotes the Kronecker symbol, thus obtaining $N_{\phi}=\mathcal{M}$ basis functions $\{\phi_k\}_{k\in[1,N_{\phi}]}$. 

Likewise, another set of \textit{nodal} basis functions is defined on the staggered quadrilateral subgrid, where the reference element $Q_{std}$ is now given by the unit square, that is $Q_{std} = \{ \xxi \in \R^2 : 0 \leq \xi_1 \leq 1, \, 0 \leq \xi_2 \leq 1 \}$. The following transformation can be used to map the physical element $Q_e$ from the reference square $Q_{std}$:
\begin{equation} 
	\label{eqn.xietaTransf2} 
	\mathbf{x} := \xx(Q_{e},\xxi) = (1-\xi_1)(1-\xi_2) \, \xx_{1,m} + \xi_1(1-\xi_2) \, \xx_{2,m} + \xi_1\xi_2 \, \xx_{3,m} + (1-\xi_1)\xi_2 \, \xx_{4,m},
\end{equation} 
with $\xx_{p,m}$ ($p=\{1,2,3,4\}$) being the vector of physical spatial coordinates of the $p$-th vertex of the quadrilateral subcell with counterclockwise orientation. {In analogy, $\xxi(Q_e,\xx):Q_e \rightarrow Q_{std}$ will indicate the inverse mapping.} The coordinates of the degrees of freedom are again computed using \eqref{eqn.NCnodes}, but with the index ranges given by $0 \leq l_1 \leq M, \quad 0 \leq l_2 \leq M$. Therefore, a set of $N_{\psi}=(M+1)^2$ nodal basis functions $\{\psi_k\}_{k\in[1,N_{\psi}]}$ are obtained. Notice that the nodal points defined by \eqref{eqn.NCnodes} correspond to the one-dimensional Newton-Cotes quadrature points (see \cite{stroud}), thus the nodal basis on $Q_{std}$ is constructed by a tensor product of the one-dimensional basis made of $(M+1)$ nodes.
{From the basis functions $\{\phi_k\}_k$ and $\{\psi_k\}_k$ defined on the reference space, it is easy to obtain the basis functions on the physical space using the element-based transformations defined above:
\begin{equation}
	\phi_k^{(m)}(\xx)=\phi_k(\xxi(T_m,\xx)), \qquad \psi_k^{(e)}(\xx)=\psi_k(\xxi(Q_e,\xx)).
\end{equation}

To transfer numerical data from one basis to another we make use of $L_2$-projection operators. Let $\boldsymbol{\gamma}_i$ be a generic quantity which is numerically represented by the Taylor modal basis functions \eqref{eqn.Voronoi_modal} on the Voronoi cell $P_i$:
\begin{equation}
	\label{eqn.poly_Voronoi}
	\boldsymbol{\gamma}_i = \sum \limits_{l=1}^{\mathcal{M}} \beta^{(i)}_l(\xx) \, \hat{\gamma}_l := \beta^{(i)}_l \, \hat{\gamma}_{l,i},
\end{equation}
where $\hat{\gamma}_{l,i}$ are the expansion coefficients, i.e. the degrees of freedom, and Einstein summation convention is assumed over repeated indexes. The evaluation of the corresponding degrees of freedom on the sub-triangulation, that is the projection of the quantity $\boldsymbol{\gamma}_i$ from the Voronoi element $P_i$ to the subcell $T_m$, relies on the following projection operator $\Topp(\boldsymbol{\gamma}_i)$:
\begin{equation}
	\label{eqn.Topp}
	\hat{\boldsymbol{\gamma}}_{m} = \left[ \left( \int \limits_{T_m} \phi_k^{(m)} \, \phi_l^{(m)} \, d\xx \right)^{-1} \, \int \limits_{T_m} \phi_k^{(m)} \, \beta_l^{(i)} \, d\xx \right] \, \hat{\gamma}_{l,i}:=\Topp(\boldsymbol{\gamma}_i), \qquad \forall i \qquad \forall m=m(ij) \in P_i,
\end{equation}
with $\hat{\boldsymbol{\gamma}}_{m}=(\hat{\gamma}_{1,m}, \ldots \hat{\gamma}_{\mathcal{M},m})$ being the sought degrees of freedom of the nodal basis defined on the subcell $T_m$. The quantity $\boldsymbol{\gamma}_m$ can be expressed in the nodal basis on the sub-triangle as
\begin{equation}
	\label{eqn.poly_subtri}
	\boldsymbol{\gamma}_m = \phi_l^{(m)} \, \hat{\gamma}_{l,m}.
\end{equation}

Similarly, we now want to detail the projection operator $\Vop(\{ \boldsymbol{\gamma}_{m}\}_m)$ from a sub-triangle $T_m$ to a Voronoi cell $P_i$, that must retrieve the starting expansion coefficients $\hat{\gamma}_{l,i}$ used in \eqref{eqn.poly_Voronoi}. To that aim, all the subcells belonging to $P_i$ must be considered, that is $\{T_m: m(ij) \in P_i\}$, hence involving a total number of $N_{v_i}$ sub-triangles for the associated cell $P_i$. This $L_2$-operator writes
\begin{equation}
	\label{eqn.Vop}
	\hat{\boldsymbol{\gamma}}_{i} = \left( \sum \limits_{m \in P_i} \, \,  \int \limits_{T_m} \beta_k^{(i)} \, \beta_l^{(i)} \, d\xx \right)^{-1} \, \sum \limits_{m \in P_i} \, \, \int \limits_{T_m} \beta_k^{(i)} \, \phi_l^{(m)} \hat{\gamma}_{l,m}\, d\xx :=\Vop(\{ \boldsymbol{\gamma}_{m}\}_m) .
\end{equation} 
The summation over all subcells $T_{m}$ with $m \in P_i$ implies that each subcell must be associated to its expansion coefficients $\hat{\gamma}_{l,m}$, thus they can not be collected outside the rightmost integral in \eqref{eqn.Vop}. The first term on the right hand side of \eqref{eqn.Vop} is nothing but the modal mass matrix of element $P_i$, which can be computed only once and saved for all elements $P_{i \in [1,N_P]}$ at the price of some memory consumption. {Let us remark that the operator $\Vop$ can be seen as a high order average of the sub-triangular elements of $P_i$, hence requiring all the values $\{ \boldsymbol{\gamma}_{m}\}_m$ for $m \in P_i$}.

Finally, the last projection operators $\bar{\Sop}(\{\boldsymbol{\gamma}_{m}\}_m)$ and ${\Sop}(\{\boldsymbol{\gamma}_{e}\}_e)$ are introduced, which are used to transfer data from the subcell $T_m$ to the staggered quadrilateral element $Q_e$ and viceversa, respectively. Also in this case, the starting sub-triangles must be the neighbor elements of the edge $\Gamma_e$, which indeed contain the staggered cell $Q_e$ (see Figure \ref{fig.mesh}). Therefore, one has $m=\{\ell(e),r(e)\}$, and the first operator is defined as
\begin{equation}
	\label{eqn.Sop}
	\hat{\boldsymbol{\gamma}}_{e} = \left( \int \limits_{Q_e} \psi_k^{(e)} \, \psi_l^{(e)} \, d\xx \right)^{-1} \, \left( \int \limits_{T_{\ell(e),e}} \psi_k^{(e)} \phi_l^{(\ell(e))} \, d\xx \, \hat{\gamma}_{l,\ell(e)} + \int \limits_{T_{r(e),e}} \psi_k^{(e)} \phi_l^{(r(e))} \, d\xx \, \hat{\gamma}_{l,r(e)} \right):=\bar{\Sop}(\{\boldsymbol{\gamma}_{m}\}_m).
\end{equation}  
where $T_{m,e}:=T_m \cap Q_e$ is the intersection triangle between $T_m$ and $T_e$ with $e \in S_m$. The inverse map, that gives the second operator, is simply defined as
\begin{equation}
	\hat{\boldsymbol{\gamma}}_{m} = \left( \int \limits_{T_m} \phi_k^{(m)} \, \phi_l^{(m)} \, d\xx \right)^{-1} \,  \sum\limits_{e \in S_m}\,\,
	\int \limits_{T_{m,e}} \phi_k^{(m)} \psi_l^{(e)} \, d\xx \, \hat{\gamma}_{l,e} :={\Sop}(\{\boldsymbol{\gamma}_{e}\}_e).
\end{equation}

\subsection{First order semi-discrete scheme in time} \label{ssec.semidiscr}
The time discretization is based on the class of semi-implicit IMEX schemes proposed in \cite{BosFil2016}, which have been recently used in all Mach solvers for compressible flows \cite{BosPar2021,SICNS22,Boscarino22}. Let us consider the following first order in time semi-discrete scheme for the shallow water system \eqref{eqn.SWE}:  
\begin{subequations}
	\label{eqn.FOtime}
	\begin{align}
		\frac{\eta^{n+1}-\eta^n}{\dt} + \nabla \cdot \q^{n+1} &= 0, \label{eqn.FOtime_eta} \\
		\frac{\q^{n+1}-\q^n}{\dt}  + \nabla \cdot \left( \vv^n \otimes \q^n \right) + g H^n \nabla \eta^{n+1} &= \bzero, \label{eqn.FOtime_q}
	\end{align}
\end{subequations}
where the equation for the bottom elevation \eqref{eqn.SWE_b} has been neglected since $b(\xx)$ is time-independent. We can easily observe that the semi-discrete scheme \eqref{eqn.FOtime} is concerned with an implicit discretization of the pressure sub-system \eqref{eqn.SWEp} and an explicit treatment of the convective sub-system \eqref{eqn.SWEc}. According to \cite{Casulli1990}, system \eqref{eqn.FOtime} is solved by substitution. Indeed, inserting the discharge equation \eqref{eqn.FOtime_q} into the mass equation \eqref{eqn.FOtime_eta} leads to the following wave equation where the only unknown is the free surface elevation $\eta^{n+1}$:
\begin{equation}
	\label{eqn.Poisson}
	\eta^{n+1} + \dt^2 \, g \, \nabla \cdot \left( H^n \, \nabla \eta^{n+1}\right) = \eta^n - \dt \nabla \cdot \q^{*}, \qquad \q^* = \q^n - \dt \nabla \cdot \left( \vv^n \otimes \q^n \right),
\end{equation} 
where the contribution of the nonlinear convective terms is compactly written with the abbreviation $\q^*$. Once the linear system \eqref{eqn.Poisson} is solved, the new free surface elevation $\eta^{n+1}$ is used to update the flow discharge from \eqref{eqn.FOtime_q}, hence 
\begin{equation}
	\label{eqn.Vupdate}
	\q^{n+1} = \q^* - \dt g H^n \nabla \eta^{n+1}.
\end{equation}

\begin{theorem}\label{th1}(Well-balance property). 
Assuming periodic boundary conditions on $\partial \Omega \in \R$ and assuming the following initial condition
\begin{equation}
	\label{eqn.QWB}
	\eta(\xx,0) = \eta_0, \qquad \vv(\xx,0)=\bzero, \qquad b(\xx) \neq 0,
\end{equation}
the semi-discrete scheme \eqref{eqn.Poisson}-\eqref{eqn.Vupdate} is well-balanced in the sense of \cite{WBLeVeque}.
\end{theorem}

\begin{proof}
From the initial condition it follows that $\eta^n=\eta_0$ and $\vv^n=0$, thus $\q^n=\bzero$. Therefore, the nonlinear convective contribution vanishes as well, i.e. $\q^*=\bzero$, since the numerical solution does not present any discontinuity. Indeed, the numerical dissipation associated to the numerical flux of the convective term $\nabla \cdot \left(\vv^n \otimes \q^n \right)$ is exactly zero for any constant state, included $\q^n=\bzero$. Consequently, the wave equation \eqref{eqn.Poisson} reduces to
\begin{equation}
	\label{eqn.Poisson2}
	\eta^{n+1} - \dt^2 \, g \, \nabla \cdot \left( H^n \, \nabla \eta^{n+1}\right) = \eta^n, \qquad \q^* = \bzero,
\end{equation}
which admits the solution $\eta^{n+1}=\eta^n=\eta_0$, implying that $\nabla \eta^{n+1}=0$. The discharge equation is then updated according to \eqref{eqn.Vupdate}, thus obtaining
\begin{equation}
	\q^{n+1} = \bzero - \dt g H^n \cdot 0 = \bzero.
\end{equation}
Therefore, the semi-discrete scheme \eqref{eqn.Poisson}-\eqref{eqn.Vupdate} can preserve stationary solutions of the shallow water system of the form given by \eqref{eqn.QWB} with arbitrary bathymetry.
\end{proof}

\begin{theorem}\label{th2}(Asymptotic Preserving property). Assuming periodic boundary conditions on $\partial \Omega \in \R$, the semi-discrete scheme \eqref{eqn.Poisson}-\eqref{eqn.Vupdate} is a consistent approximation of the low Froude shallow water system \eqref{eqn.LowFr} at the leading order asymptotic expansion in the asymptotic limit ($\varepsilon \to 0$).	
\end{theorem}

\begin{proof}
Using the rescaled variables \eqref{eqn.scale_var} and the expansions \eqref{eqn.exp}, the semi-discrete scheme \eqref{eqn.FOtime} in non-dimensional form writes	
\begin{subequations}
	\label{eqn.FOtime_scaled}
	\begin{align}
		\frac{\eta^{n+1}-\eta^n}{\dt} + \nabla \cdot \q^{n+1} &= 0, \label{eqn.FOtime_scaled_eta} \\
		\frac{\q^{n+1}-\q^n}{\dt}  + \nabla \cdot \left( \vv^n \otimes \q^n \right) + \frac{H^n}{\varepsilon^2} \nabla \eta^{n+1} &= \bzero. \label{eqn.FOtime_scaled_q}
	\end{align}
\end{subequations}
Let us assume that the following expansions hold true for the discrete variables at any generic time $t^n$:
\begin{equation}
	\label{eqn.exp2}
	\eta^n(\xx) = \eta_{(0)}^n (\xx) + \varepsilon^2 \eta_{(2)}^n(\xx), \qquad \vv^n(\xx) = \vv_{(0)}^n (\xx) + \varepsilon \vv_{(1)}^n(\xx),
\end{equation}
where $\eta_{(0)}^n(\xx)=\eta_{(0)}$ is constant in space and time because periodic boundaries are assumed (see Section \ref{ssec.lowFroude}), and $\varepsilon^2 \eta_{(2)}^n(\xx)$ is a  perturbation of the free surface level, thus the total water depth at zeroth order is given by $H_{(0)}^n(\xx)=\eta_{(0)}(\xx)-b(\xx)$, which is also constant. Inserting \eqref{eqn.exp2} into the semi-discrete scheme \eqref{eqn.FOtime_scaled} and retaining only zeroth order terms of the expansions lead to
\begin{subequations}
	\label{eqn.FOtime_scaled2}
	\begin{align}
		\nabla \cdot \left( H_{(0)} \vv_{(0)} \right)^{n+1} &= 0, \label{eqn.FOtime_scaled2_eta} \\
		\frac{\left( H_{(0)} \vv_{(0)}\right)^{n+1}-\left( H_{(0)} \vv_{(0)} \right)^n}{\dt}  + \nabla \cdot \left( \vv^n \otimes (H_{(0)} \vv_{(0)})^n \right) + H_{(0)}^n \nabla \eta_{(2)}^{n+1} &= \bzero, \label{eqn.FOtime_scaled2_q}
	\end{align}
\end{subequations} 
that is a consistent discretization at first order in time of the low Froude shallow water system \eqref{eqn.LowFr}. Formal substitution of \eqref{eqn.FOtime_scaled2_q} into \eqref{eqn.FOtime_scaled2_eta} gives the the corresponding rescaled version of \eqref{eqn.Poisson}-\eqref{eqn.Vupdate} in the asymptotic limit, namely
\begin{subequations}
	\begin{align}
		&H_{(0)} \, \nabla \cdot \nabla \eta_{(2)}^{n+1} = \nabla \cdot \left( H_{(0)} \vv_{(0)}\right)^*, \qquad \left( H_{(0)} \vv_{(0)}\right)^* = \left( H_{(0)} \vv_{(0)} \right)^n - \dt \nabla \cdot \left( \vv^n \otimes (H_{(0)} \vv_{(0)})^n \right) \\
		&\left( H_{(0)} \vv_{(0)}\right)^{n+1} = \left( H_{(0)} \vv_{(0)}\right)^*  - \dt H_{(0)}^n \nabla \eta_{(2)}^{n+1},
	\end{align}
\end{subequations}
which can be equivalently obtained by inserting the expansions \eqref{eqn.exp2} into the rescaled scheme \eqref{eqn.FOtime_scaled2_eta}-\eqref{eqn.FOtime_scaled2_q}.
\end{proof}

\subsection{High order time discretization} 
Once the first order in time semi-discrete scheme \eqref{eqn.FOtime} is designed, its extension to high order of accuracy in time is carried out adopting the semi-implicit IMEX time integrators firstly introduced in \cite{BosFil2016}. The governing equations can be cast in the form of an autonomous system
\begin{equation}
	\label{eqn.partSyst}
\frac{\partial \U}{\partial t} = \mathcal{H}\left(\U_{E}(t), \U_{I}(t) \right), \qquad \U_0=\U(t=0),
\end{equation}
where the vector of conserved variables is $\U=(\eta,Hu,Hv,b)$ according to \eqref{eqn.SWE}. The function $\mathcal{H}$ represents any spatial approximation of the remaining terms of the shallow water system, which will be detailed in the next sections. The first argument of $\mathcal{H}$ denoted with $\U_{E}$ is discretized explicitly, and the second argument referred to as $\U_{I}$ is taken implicitly, thus obtaining a partitioned system. Looking at the semi-discrete scheme \eqref{eqn.FOtime}, the right hand side of \eqref{eqn.partSyst} results to be 
\begin{equation}
\mathcal{H}\left(\U_E, \U_I \right) = \left\{  \begin{array}{c}
-\nabla_h \cdot \q_I \\ - \nabla_h \cdot \left( \vv \otimes \q \right)_E - g H_E \nabla_h \eta_I \end{array} \right. ,
\label{eqn.partSys}
\end{equation}
with the \textit{discrete} divergence operator $\nabla_h$ that will be presented in the next section. The class of implicit-explicit (IMEX) Runge-Kutta schemes \cite{PR_IMEX} allows high order in time to be reached by performing a total number $s$ of stages which depend on the desired order of accuracy and other constraints on the asymptotic preserving property of the scheme. In this work, we use the semi-implicit IMEX schemes up to third order detailed in \cite{BosPar2021} (see \ref{app.IMEX}), which are proven to be asymptotic preserving and asymptotic accurate. Consequently, since the first order in time semi-discrete scheme is asymptotic preserving as demonstrated by Theorem \ref{th2}, its high order extension maintains the Asymptotic Preserving property, see \cite{BosFil2016} for a detailed proof. Let us also remark that the duplication of the unknowns $\U_E$ and $\U_I$ in \eqref{eqn.partSys} does not take place if judicious choices of the IMEX scheme are considered \cite{BosFil2016}.

For practical implementation, the IMEX schemes are typically represented with the double Butcher tableau:
\begin{equation}
	\label{eqn.butcher}
	\begin{array}{c|c}
		\tilde{c} & \tilde{A} \\ \hline & \tilde{b}^\top
	\end{array} \qquad
	\begin{array}{c|c}
		c & A \\ \hline & b^\top
	\end{array},
\end{equation}
with the matrices $(\tilde{A},A) \in \R^{s \times s}$ and the vectors $(\tilde{c},c,\tilde{b},b) \in \R^s$. The tilde symbol refers to the explicit scheme and matrix $\tilde{A}=(\tilde{a}_{ij})$ is a lower triangular matrix with zero elements on the diagonal, while $A=({a}_{ij})$ is a triangular matrix which accounts for the implicit scheme, thus having non-zero elements on the diagonal. A semi-implicit IMEX Runge-Kutta method is obtained as follows. Let us first set $\U_E^n=\U_I^n=\U^n$, then the stage fluxes for $i = 1, \ldots, s$ are calculated as
\begin{subequations}
	\begin{align}
		\U_E^i &= \U_E^n + \dt \sum \limits_{j=1}^{i-1} \tilde{a}_{ij} k_j, \qquad 2 \leq i \leq s, \label{eq.QE} \\[0.5pt]
		\tilde{\U}_I^i &= \U_E^n + \dt \sum \limits_{j=1}^{i-1} a_{ij} k_j, \qquad 2 \leq i \leq s, \label{eq.QI}  \\[0.5pt]
		k_i &= \mathcal{H} \left( \U_E^i, \tilde{\U}_I^i + \dt \, a_{ii} \, k_i \right), \qquad 1 \leq i \leq s. \label{eq.k} 
	\end{align}
\end{subequations}
Finally, the numerical solution is updated with
\begin{equation}
	\U^{n+1} = \U^n + \dt \sum \limits_{i=1}^s b_i k_i.
	\label{eqn.QRKfinal}
\end{equation}

\subsection{Spatial discretization of the explicit terms} \label{ssec.cweno} 
The vector of conserved variables $\U=(\eta,Hu,Hv,b)$ is stored for every time level $t^n$ within each Voronoi cell as typically done in finite volume schemes:
\begin{equation}
	\label{eqn.cellAv}
	\U_i^n:=\frac{1}{|P_i|} \int \limits_{P_i} \U(\xx,t^n) \, d\xx.
\end{equation}
The spatial discretization is composed of two main steps: (i) a high order nonlinear reconstruction, and (ii) a finite volume scheme on unstructured Voronoi meshes.

\paragraph{CWENO reconstruction on Voronoi meshes} Starting from the known cell averages \eqref{eqn.cellAv}, a reconstruction polynomial $\w(\xx,t^n)$ of arbitrary degree $M$ is computed relying on the CWENO strategy originally forwarded in \cite{LPR:99,LPR:2001} and subsequently used also in the context of unstructured meshes \cite{DGCWENO,ADER_CWENO,ArepoTN,FVBoltz,DGBoltz}. The reconstruction polynomial $\w(\xx,t^n)$ is expressed for each cell $P_i$ by means of the Taylor basis \eqref{eqn.Voronoi_modal}, that is
\begin{equation}
	\left. \w(\xx,t^n) \right|_{P_i}:=\w_i^n = \beta_l^{(i)}(\xx) \, \hat{\w}_{l,i}^n.
\end{equation}
Because of the conservative modal basis, the conservation property \eqref{eqn.modal_cons} implies $\hat{\w}_{1,i}^n=\U_i$. The reconstruction procedure is then compactly written by defining the following operator $\Rop$: 
\begin{equation}
	\label{eqn.Rop}
	\hat{\w}_{l,i}^n = \Rop(\U^n),
\end{equation}
and the explicit definition of $\Rop$, thus the details of the CWENO reconstruction algorithm, can be found in \ref{app.CWENO}. We underline that the reconstruction strategy is of arbitrary order of accuracy, thus any polynomial degree $M$ can be chosen. The higher is the accuracy, the most expensive is the computational effort and the larger is the reconstruction stencil, hence making finite volume reconstruction schemes less efficient for parallelization purposes.

\paragraph{Finite volume scheme} Once the CWENO reconstruction procedure is carried out for all the Voronoi elements, a finite volume scheme is used to discretize the nonlinear convective operators in \eqref{eqn.FOtime}. Therefore, shock capturing properties and conservation are ensured by construction for the convective sub-system \eqref{eqn.SWEc}. Integration of the discharge equation in \eqref{eqn.SWEc} over the control volume $P_i$ and application of Gauss theorem yields
\begin{equation}
	\label{eqn.SWEc_Gauss}
	\partial_t \int \limits_{P_i} \q_i \, d\xx = - \int \limits_{\partial P_i} (\vv \otimes \q)_i \cdot \nn \, dS,
\end{equation}
which is numerically approximated using a finite volume scheme:
\begin{equation}
	\label{eqn.fvscheme}
	\q_i^{n+1} = \q_i^n - \frac{\dt}{|P_i|} \sum \limits_{i=1}^{N_{v_i}} \int \limits_{\partial P_{ij}} \mathcal{F}(\w_i^n,\w_j^n,\nn_{ij}) \, dS:=\q_i^*.
\end{equation}
The numerical flux function $\mathcal{F}$ is fed by high order extrapolated values at the boundary $\partial P_{ij}$ which come from the CWENO reconstruction. We choose to use a robust Rusanov--type numerical flux, hence defining
\begin{equation}
	\mathcal{F}(\w_i^n,\w_j^n,\nn_{ij}) = \frac{1}{2} \left( \Rop(\vv_i^n) \otimes \Rop(\q_i^n) + \Rop(\vv_j^n) \otimes \Rop(\q_j^n) \right) \cdot \nn_{ij} - \frac{1}{2} |s_{\max}| \left( \Rop(\q_j^n) - \Rop(\q_i^n) \right),
\end{equation}
where the numerical dissipation $s_{\max}$ is given by the maximum eigenvalue of the convective sub-system \eqref{eqn.SWEc}, thus it is proportional to the flow velocity and not to the acoustic-gravity wave speed. In the low Froude asymptotic limit this is very important since numerical dissipation automatically tends to zero for $\varepsilon \to 0$ in \eqref{eqn.SWErr}. Furthermore, for any constant solution $\q=\q_0$, the numerical flux contribution vanishes because of Gauss theorem \eqref{eqn.SWEc_Gauss}, giving evidence that the convective term maintains the well-balance solution proved in Theorem \ref{th1}. Let us notice that the right hand side of the finite volume scheme \eqref{eqn.fvscheme} corresponds exactly to the spatial discretization of the term $\q^*$ in \eqref{eqn.Poisson}, which is referred to as $\q_i^*$. The discrete explicit fluxes $\mathcal{H}(\U_E)$ in \eqref{eqn.partSys} are then given by
\begin{equation}
	\mathcal{H}\left(\U_E\right) = \left\{  \begin{array}{c}
		0 \\ - \frac{\dt}{|P_i|} \sum \limits_{i=1}^{N_{v_i}} \int \limits_{\partial P_{ij}} \mathcal{F}(\w_i^n,\w_j^n,\nn_{ij}) \, dS \end{array} \right. .
	\label{eqn.HE}
\end{equation}

\bigskip

The result of the finite volume scheme is therefore $\q_i^*$, which formally provides a cell average of the type \eqref{eqn.cellAv}. The implicit solver for the free surface elevation, which will be described in the next section, is based on a discontinuous Galerkin representation of the numerical solution. Consequently, to make the convective numerical solution $\q_i^*$ suitable for a DG method, we need to perform a reconstruction of $\q_i^*$. The resulting CWENO polynomial is then \textit{interpreted} as a discontinuous Galerkin numerical solution within each computational cell, which has been very recently proposed in the context of IMEX solvers for the incompressible Navier-Stokes equations \cite{SIINS22}. Therefore, a CWENO reconstruction is performed after the convective terms have been updated, hence obtaining high order reconstruction polynomials on the Voronoi cells for the following quantities:
\begin{equation}
	\label{eqn.UHO}
	\hat{\eta}_{l,i}^n=\Rop(\eta^n), \qquad \hat{\q}_{l,i}^*=\Rop(\q^*), \qquad \hat{H}_{l,i}^n=\Rop(H^n), \qquad i=1,\ldots,N_P.
\end{equation}


\subsection{Spatial discretization of the implicit terms} \label{ssec.stagdg}
The implicit discretization makes use of a discontinuous Galerkin approximation on the triangular subcells and the associated staggered quadrilateral subgrid. Therefore, the free surface elevation as well as flow discharge and water depth must comply with the DG data structure. Specifically, the free surface elevation $\eta$ has a high order representation on each sub-triangle while the total water depth $H$ and flow discharge $\q$ can be represented by a high order polynomial on the staggered elements $Q_e$. This is computed starting from the high order polynomials \eqref{eqn.UHO} defined on the Voronoi cells and relying on the projections described in Section \ref{ssec.basefunc}, that is
\begin{equation}
	\label{eqn.iniDG}
	\{\hat{\eta}_{m}^n\}_m = \Topp(\{\eta_{i}^n\}_{i}), \qquad \{\hat{\q}_{e}^*\}_e=\bar{\Sop}\left(\Topp\left(\{{\q}_{i}^*\}_i\right)\right), \qquad \{\hat{H}_e^n\}_e=\bar{\Sop}\left(\Topp(\{H_i^n\}_i)\right).
\end{equation}
We remark that the input data in the above projections are given by the CWENO polynomials \eqref{eqn.UHO}. As a consequence, the discrete variables obtained using \eqref{eqn.iniDG} are explicitly approximated as follows:
\begin{eqnarray}
	\label{eqn.expDG}
	\eta(\xx)|_{\xx \in T_m}=\sum_k\phi_k^{(m)}(\xx) \, \, \hat{\eta}_{m,k} &=&\phi_k \, \hat{\eta}_m \nonumber \\
	H(\xx)|_{\xx \in Q_e}=\sum_k\psi_k^{(e)}(\xx) \, \, \hat{H}_{e,k}&=&\psi_k \, \hat{H}_e, \nonumber \\
	\q(\xx)|_{\xx \in Q_e}=\sum_k\psi_k^{(e)}(\xx) \, \, \hat{\q}_{e,k}&=&\psi_k \, \hat{\q}_e,
\end{eqnarray}
where we introduce a lighter notation, i.e. we simply use $\hat{\eta}_m$ for all $m=1,\ldots, N_m$ and $(\hat{H}_e,\hat{\q}_e)$ for all $e=1, \ldots, N_e$. We can now consider the implicit contributions in the semi-discrete scheme \eqref{eqn.FOtime}:
\begin{subequations}
	\label{eqn.FOtime2}
	\begin{align}
		\eta^{n+1} + \dt \nabla \cdot \q^{n+1} &= \eta^n, \label{eqn.FOtime_eta2} \\
		\q^{n+1} + \dt g H^n \nabla \eta^{n+1} &= \q^*. \label{eqn.FOtime_q2}
	\end{align}
\end{subequations}
A weak formulation of the momentum and continuity equation is derived following the approach presented in \cite{TavelliSWE2014}. Multiplying the mass equation \eqref{eqn.FOtime_eta2} by a test function $\phi_k$, integrating it over the control volume $T_m$ and inserting the ansatz \eqref{eqn.expDG} yields
\begin{equation}
	\bar{M}_m \hat{\eta}_m^{n+1}+\dt \sum\limits_{e \in S_m} \DD_{m,e} \hat{\q}_e^{n+1}=\bar{M}_m \hat{\eta}_m^n, \label{weakcont}
\end{equation}
where 
\begin{eqnarray}
	\bar{M}_m&=&\int\limits_{T_m} \phi_k^{(m)}(\xx) \, \phi_l^{(m)}(\xx) \, d\xx, \nonumber \\
	\DD_{m,e} &=& \int\limits_{\Gamma_e}\phi_k^{(m)}(s) \, \psi_l^{(e)}(s) \, \nn_e \, \sigma_{me} \, ds - \int\limits_{T_{m,e}} \nabla \phi_k^{(m)}(\xx) \, \psi_l^{(e)}(\xx) \,d\xx, \nonumber \\
	\sigma_{me}&=&\frac{r(e)-2m+\ell(e)}{r(e)-\ell(e)}.
\end{eqnarray}
In analogy, a weak formulation of the momentum equation may be obtained after multiplication of Equation \eqref{eqn.FOtime_q2} by a test function $\psi_k$ and integration over the staggered subcell $Q_e$:
 \begin{equation}
 	M_e \hat{\q}_e^{n+1} + \dt g \left(\QQ_{\ell(e),e} \hat{H}_e^n \hat{\eta}_{\ell(e)}^{n+1} + \QQ_{r(e),e} \hat{H}_e^n \hat{\eta}_{r(e)}^{n+1}  \right)= M_e \hat{\q}_e^*, \label{weakmom}
 \end{equation}
where
\begin{eqnarray}
	{M}_e&=&\int\limits_{Q_e} \psi_k^{(e)}(\xx) \, \psi_l^{(e)}(\xx) \, d\xx , \nonumber \\
	\QQ_{m,e} &=& -\int\limits_{\Gamma_e}\psi_k^{(e)}(s) \, \psi_l^{(e)}(s)\, \phi_r^{(m)}(s) \, \nn_e \, \sigma_{me} \,  ds + \int\limits_{T_{m,e}}\psi_k^{(e)}(\xx) \, \psi_l^{(e)}(\xx) \, \nabla \phi_r^{(m)}(\xx) \, d\xx .
\end{eqnarray}
Mimicking what done at the semi-discrete level, substitution of the weak momentum into the weak continuity equation leads to a linear system for the only unknowns $\hat{\eta}^{n+1}$, thus
\begin{eqnarray}
	\bar{M}_m \eta_m^{n+1}+g \dt^2 \sum\limits_{e \in S_m} \DD_{m,e} \left(\QQ_{\ell(e),e} \hat{H}_e^n \hat{\eta}_{\ell(e)}^{n+1} + \QQ_{r(e),e} \hat{H}_e^n \hat{\eta}_{r(e)}^{n+1}  \right) =\bar{M}_m \eta_m^{n}+ \dt \sum\limits_{e \in S_m} \DD_{m,e} \hat{\q}^*_{e}=b_m^n,
\end{eqnarray}
which is the fully discrete wave equation \eqref{eqn.Poisson}. The system can be solved using the GMRES algorithm, see \cite{TavelliSWE2014}. Once the new free surface is computed, the momentum can readily be updated with \eqref{weakmom} obtaining $\hat{\q}_e^{n+1}$. 

The conservative variables $\U$ are defined on the Voronoi tessellation and are given as cell averages in the finite volume framework according to \eqref{eqn.cellAv}. Consequently, the new DG solution is projected back to the Voronoi grid using the operator \eqref{eqn.Vop}:
\begin{equation}
	\hat{\eta}_i^{n+1} = \Vop\left(\{\hat{\eta}_m^{n+1}\}_m\right), \qquad \hat{\q}_i^{n+1}=\Vop \left( \left\{  \Sop\{ \hat{\q}_e^{n+1} \}_e \right\}_m \right).
\end{equation} 
Next, the first degree of freedom of each quantity provides the sought cell averages according to the conservation property \eqref{eqn.modal_cons}, thus
\begin{equation}
	\eta_i^{n+1}=\hat{\eta}_{1,i}, \qquad \q_i^{n+1}=\hat{\q}_{1,i}^{n+1}.
\end{equation}
Finally, to comply with the time discretization , the discrete implicit fluxes in \eqref{eqn.partSys} are simply computed by
\begin{equation}
	\mathcal{H}\left(\U_I\right) = \left\{  \begin{array}{c}
		\frac{\eta_i^{n+1}-\eta_i}{A_{ss} \, \dt} \\ \frac{\q_i^{n+1}-\q_i^*}{A_{ss} \, \dt} \end{array} \right. ,
	\label{eqn.HI}
\end{equation}
where $A_{ss}$ is the diagonally implicit coefficient of the Butcher tableau \eqref{eqn.butcher}, which is never zero (see \ref{app.IMEX}).


\section{Numerical results} \label{sec.numtest}
We present a suite of test cases that aim at assessing the robustness and the accuracy of the novel numerical method presented in this work. The label SI-FVDG (Semi-Implicit Finite Volume/Discontinuous Galerkin) is used, and the third order version of the scheme ($M=2$) is adopted by default in both space and time. The CFL number is set to $\textnormal{CFL}=0.9$ and the time step size is determined according to the stability condition \eqref{eqn.timestep}, thus it is independent of the acoustic-gravity wave speed. If the initial fluid velocity is set to zero, which would yield to a time step approaching infinity, the first time step is computed using the eigenvalues of the full system, so that after one time step an initial amount of momentum can take place in the flow that allows the condition \eqref{eqn.timestep} to be used again. Different fluid regimes are simulated with Froude numbers ranging from $\Fr=10^{-6}$ up to $\Fr=5.73$, demonstrating the ability of the SI-FVDG method to deal with multiscale flow conditions.  

\subsection{Convergence rates study}
The numerical convergence of the new SI-FVDG schemes is studied by considering the steady shallow water vortex firstly proposed in \cite{Voronoi}. The computational domain is the square $\Omega=[-5;5]\times[-5;5]$ with flat bottom ($b=0$) and periodic boundaries, while the initial condition, which also corresponds to the exact solution, is given by
\begin{equation}
	\label{eqn.SWEVortex-ini}
	\eta(\xx) = H_0 - \frac{1}{2 g} e^{-(r^2-1)}, \qquad \vv(\xx)=(u(\xx),v(\xx))=\left( -u_{\alpha} \sin(\alpha), \,  u_{\alpha} \cos(\alpha) \right),
\end{equation} 
with polar coordinates $(\alpha,r)$ defined as
\begin{equation}
	\tan(\alpha)=\frac{y}{x}, \qquad r^2=x^2+y^2.
\end{equation}
The angular velocity $u_{\alpha}=r \, e^{-\frac{1}{2}(r^2-1)}$ is prescribed, so that the momentum equation in radial direction gives rise to a balance between centrifugal and pressure forces:
\begin{equation}
	\frac{\partial \eta}{\partial r} = \frac{u_{\alpha}}{g r}.
\end{equation}
Different regimes of the Froude number can be taken into account by modifying the constant $H_0$ in the definition of the free surface elevation \eqref{eqn.SWEVortex-ini}, thus permitting to numerically verify the asymptotic preserving property of the novel schemes. To ease the computation of the Froude number, in this test case we set $g=10$. 

This test is run on a series of successfully refined computational meshes until the final time $t_f=0.1$ with four different Froude numbers $\Fr=\{0.32, 10^{-2}, 10^{-4}, 10^{-6}\}$. For the chosen values of the Froude number the associated values of the free surface constant are given by $H_0=\{10^0, 10^{3}, 10^{7}, 10^{11}\}$, hence making necessary the adoption of quadruple precision arithmetic for the computations, as already observed in \cite{Busto_SWE2022}. The errors are measured in $L_2$ norm for the free surface elevation and the horizontal velocity component, that is
\begin{equation}
	{L_2}(\eta) = \sqrt{\int_{\Omega} (\eta_h(\xx) - \eta_e(\xx))^2 \, d\xx}, \qquad  {L_2}(u) = \sqrt{\int_{\Omega} (u_h(\xx) - u_e(\xx))^2 \, d\xx},
\end{equation}  
where $\eta_h(\xx),u_h(\xx)$ is the numerical solution expressed in terms of the CWENO reconstruction polynomials, whereas the exact solution $\eta_e(\xx),u_e(\xx)$ is given by \eqref{eqn.SWEVortex-ini}. The results are reported in Table \ref{tab.conv_rate} for second and third order schemes in space and time, where the characteristic mesh size of each computational mesh is simply evaluated as $h(\Omega)=\max \limits_i \sqrt{|P_i|}$. The formal order of accuracy is obtained for both approximation degrees ($N=[1,2]$) and for all the Froude numbers, confirming that the SI-FVDG schemes are asymptotic preserving and asymptotic accurate, meaning that the achieved order of accuracy is independent of the Froude number, as expected. Figure \ref{fig.SWEVortex} shows the magnitude of the velocity field at Froude number $\Fr=0.32$ and $\Fr=10^{-6}$, where no visible differences can be noticed despite the jump of about ten orders of magnitude related to the free surface elevation.

\begin{table}[!htp]
	\begin{center}
	\caption{Numerical convergence results of the SI-FVDG scheme with second and third order of accuracy in space and time using the steady shallow water vortex problem on Voronoi meshes. The errors are measured in $L_2$ norm and refer to the free surface elevation $\eta$ and velocity component $u$ at time $t=0.1$. The asymptotic preserving (AP) property of the schemes is studied by considering different Froude numbers $\Fr=\{ 0.32, 10^{-2}, 10^{-4}, 10^{-6} \}$ with the corresponding values of the water depth $H_0$.}
	\small{
	\renewcommand{\arraystretch}{1.1}	
	\begin{tabular}{cccccccccc}
			 & \multicolumn{4}{c}{SI-FVDG $\mathcal{O}(2)$} & & \multicolumn{4}{c}{SI-FVDG $\mathcal{O}(3)$} \\
			\cline{2-5} \cline{7-10}
			 $h(\Omega)$ & ${L_2}(\eta)$ & $\mathcal{O}(\eta)$ & ${L_2}(u)$ & $\mathcal{O}(u)$ & & ${L_2}(\eta)$ & $\mathcal{O}(\eta)$ & ${L_2}(u)$ & $\mathcal{O}(u)$ \\
			\hline
			& \multicolumn{9}{c}{$\Fr=0.32$, $H_0=10^{0}$ (double precision)} \\
			\hline
			4.6405E-01 & 4.6695E-03 & -    & 4.2797E-02 & -    & & 5.2378E-03 & -    & 3.9573E-02 & -     \\
			2.4889E-01 & 8.9159E-04 & 2.66 & 8.8269E-03 & 2.53 & & 7.8042E-04 & 3.06 & 5.4838E-03 & 3.17  \\
			1.6631E-01 & 3.7358E-04 & 2.16 & 3.7142E-03 & 2.15 & & 2.1356E-04 & 3.21 & 1.6362E-03 & 3.00  \\
			1.2765E-01 & 2.0589E-04 & 2.25 & 2.0544E-03 & 2.24 & & 8.4528E-05 & 3.50 & 6.7377E-04 & 3.35  \\
			\hline
			& \multicolumn{9}{c}{$\Fr=10^{-2}$, $H_0=10^{3}$ (quadruple precision)} \\
			\hline
			4.6405E-01 & 4.4898E-03 & -    & 4.4870E-02 & -    & & 5.7524E-03 & -    & 4.1907E-02 & -     \\
			2.4889E-01 & 1.1112E-03 & 2.24 & 9.1796E-03 & 2.55 & & 7.9515E-04 & 3.18 & 5.9539E-03 & 3.13  \\
			1.6631E-01 & 4.5253E-04 & 2.23 & 3.8537E-03 & 2.15 & & 2.1886E-04 & 3.20 & 1.7490E-03 & 3.04  \\
			1.2765E-01 & 2.5926E-04 & 2.11 & 2.1324E-03 & 2.24 & & 8.5740E-05 & 3.54 & 7.2117E-04 & 3.35  \\
			\hline
			& \multicolumn{9}{c}{$\Fr=10^{-4}$, $H_0=10^{7}$ (quadruple precision)} \\
			\hline
			4.6405E-01 & 4.4782E-03 & -    & 4.4604E-02 & -    & & 5.7562E-03 & -    & 4.1909E-02 & -     \\
			2.4889E-01 & 1.1180E-03 & 2.23 & 9.3643E-03 & 2.51 & & 7.9785E-04 & 3.17 & 5.9544E-03 & 3.13  \\
			1.6631E-01 & 4.8607E-04 & 2.07 & 3.9086E-03 & 2.17 & & 2.1972E-04 & 3.20 & 1.7491E-03 & 3.04  \\
			1.2765E-01 & 2.6381E-04 & 2.31 & 2.1140E-03 & 2.32 & & 8.5924E-05 & 3.55 & 7.2122E-04 & 3.35  \\
			\hline
			& \multicolumn{9}{c}{$\Fr=10^{-6}$, $H_0=10^{11}$ (quadruple precision)} \\
			\hline
			4.6405E-01 & 4.4782E-03 & -    & 4.4604E-02 & -    & & 5.7562E-03 & -    & 4.1909E-02 & -     \\
			2.4889E-01 & 1.1180E-03 & 2.23 & 9.3643E-03 & 2.51 & & 7.9785E-04 & 3.17 & 5.9544E-03 & 3.13  \\
			1.6631E-01 & 4.8607E-04 & 2.07 & 3.9086E-03 & 2.17 & & 2.1972E-04 & 3.20 & 1.7491E-03 & 3.04  \\
			1.2765E-01 & 2.6381E-04 & 2.31 & 2.1140E-03 & 2.32 & & 8.5924E-05 & 3.55 & 7.2122E-04 & 3.35  \\
	\end{tabular}
    }
	\label{tab.conv_rate}
	\end{center}
\end{table}

\begin{figure}[!htbp]
	\begin{center}
		\begin{tabular}{cc}
			\includegraphics[width=0.47\textwidth]{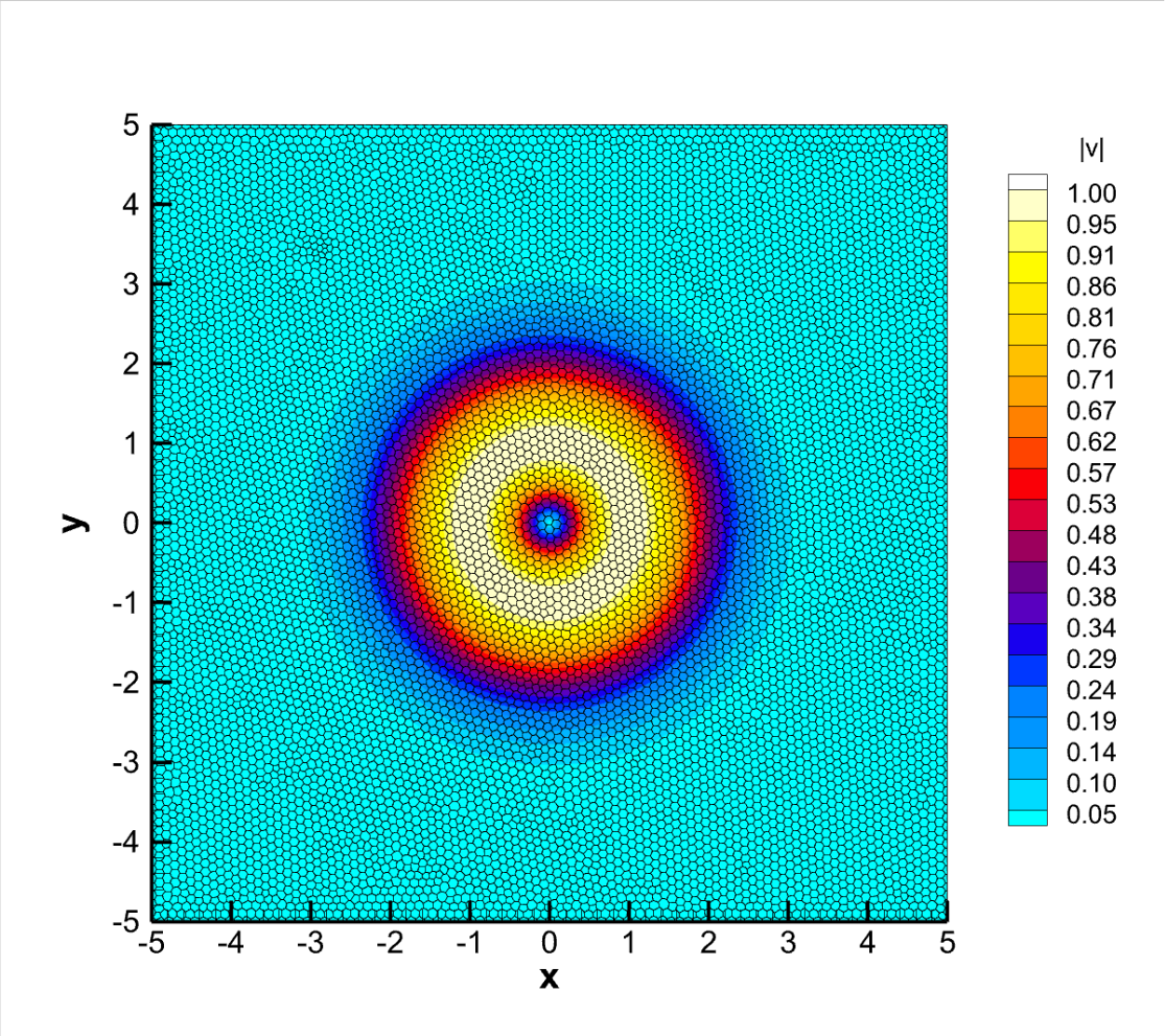}  &          
			\includegraphics[width=0.47\textwidth]{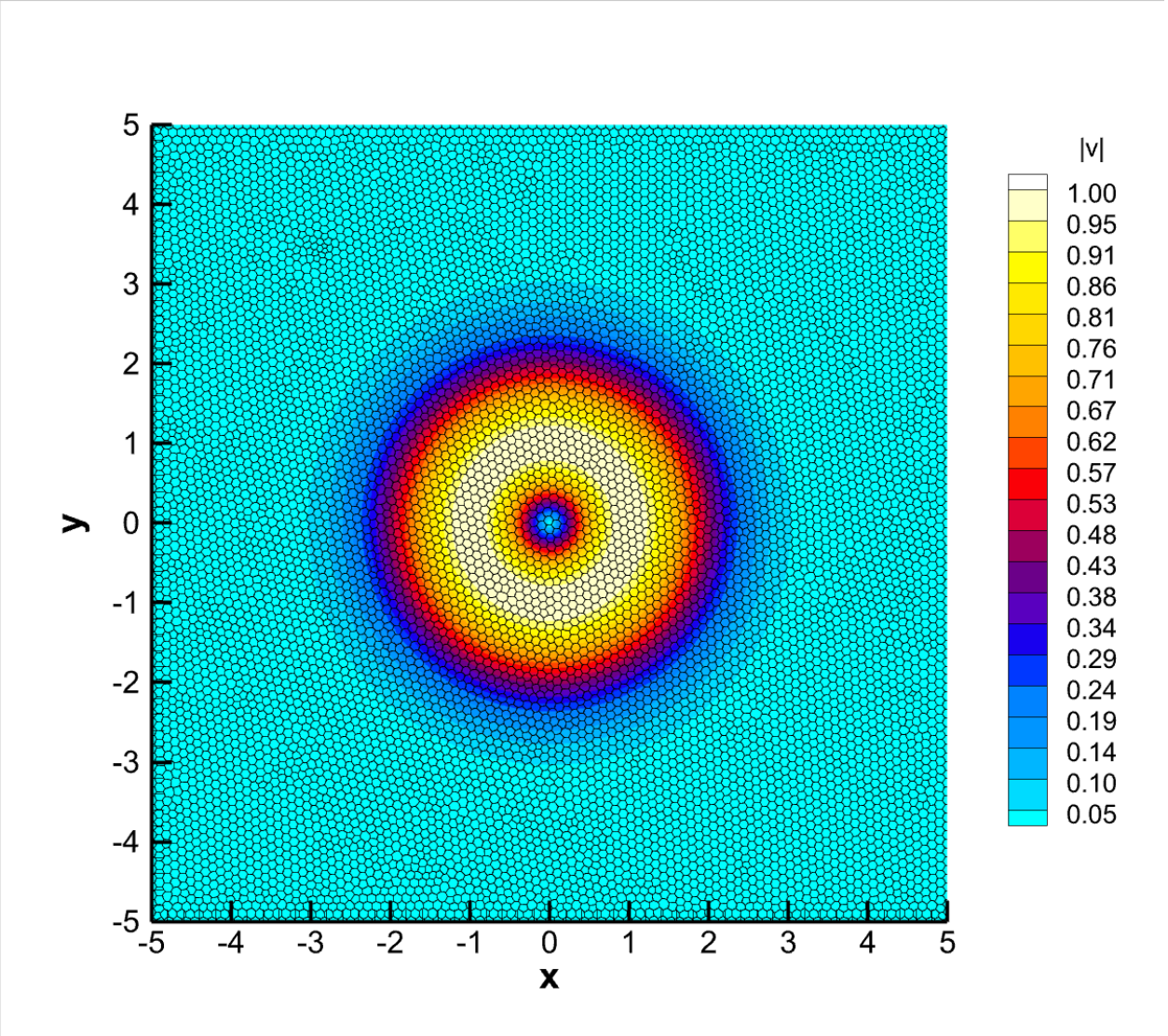}    \\
		\end{tabular}
		\caption{Steady shallow water vortex problem at $t_f=0.1$ with a Voronoi mesh of characteristic size $h\approx 1/2$. Magnitude of the velocity field obtained with $Fr=0.32$ (left) and $Fr=10^{-6}$ (right) using quadruple precision arithmetic for the computations.}
		\label{fig.SWEVortex}
	\end{center}
\end{figure}

\subsection{Well-balance test}
To numerically verify the well-balance property of the SI-FVDG schemes proved in Theorem \ref{th1}, which is also referred to as C-property, we consider the benchmark devised in \cite{WBLeVeque}. The setting of this test allows to assess whether a numerical scheme is able to preserve stationary equilibrium solutions of the governing equations up to machine precision. Specifically, equilibrium solutions of the shallow water equations are characterized by a constant free surface elevation $\eta(\xx,t)=0$ and zero fluid velocity, i.e. $\vv(\xx,t)=\mathbf{0}$, while prescribing an arbitrary bottom topography different from the trivial profile $b(\xx)=0$. Following \cite{WBLeVeque}, we consider a computational domain $\Omega=[-2;1]\times[-0.5;0.5]$ with Dirichlet boundary conditions in $x-$direction and periodic boundaries in $y-$direction, which is discretized with a mesh size of $h=1/50$, hence resulting in a total number of $N_P=8633$ Voronoi cells. The bathymetry and the initial free surface elevation are then given by
\begin{equation}
	\label{eqn.WBini}
	b(\xx) = 0.5 \cdot e^{-5 \, (x+0.1)^2 - 50 y^2}, \qquad \eta(\xx,0) = \left\{ \begin{array}{lc}
		1 + \delta & \textnormal{ if } -0.95 \leq x \leq -0.85 \\
		1 & \textnormal{elsewhere} 
	\end{array}\right. .
\end{equation} 
The fluid is initially at rest and we set the perturbation amplitude $\delta=0$. The simulation is run until the final time $t_f=0.1$ using double and quadruple finite arithmetic, and the errors with respect to the initial condition are reported in Table \ref{tab.Cproperty}. One can notice that the novel SI-FVDG scheme are well-balanced up to machine accuracy.

\begin{table}[!htp]
	\begin{center}
		\caption{Well-balance test with double and quadruple finite arithmetic precision. Errors measured in $L_2$ and $L_{\infty}$ norms for the free surface elevation $\eta$ and velocity component $u$ at the final time $t_f=0.1$.}
			\begin{tabular}{c|cccc}
				Precision & $L_2(\eta)$ & $L_{\infty}(\eta)$ & $L_2(u)$ & $L_{\infty}(u)$ \\
				\hline
				Double    & 1.3933E-15 & 3.8857E-15  & 3.4734E-14 & 8.1454E-13  \\
				Quadruple & 1.0690E-33 & 3.9481E-33  & 3.3339E-32 & 7.3247E-31
			\end{tabular}
		\label{tab.Cproperty}
	\end{center}
\end{table}

Next, as proposed in \cite{WBLeVeque}, a small perturbation is put in the free surface elevation, namely we set $\delta=10^{-2}$ in \eqref{eqn.WBini}. Here, a fixed time step of $\dt=0.01$ is adopted in order to properly follow the wave propagation. The results are depicted in Figure \ref{fig.Ctest2} at different output times, showing that no spurious oscillations are generated by the presence of the bottom bump. The flow structure is qualitatively in excellent agreement with the results available in the literature \cite{TavelliSWE2014,Canestrelli2010,Busto_SWE2022}.

\begin{figure}[!htbp]
	\begin{center}
		\begin{tabular}{cc}
			\includegraphics[width=0.49\textwidth]{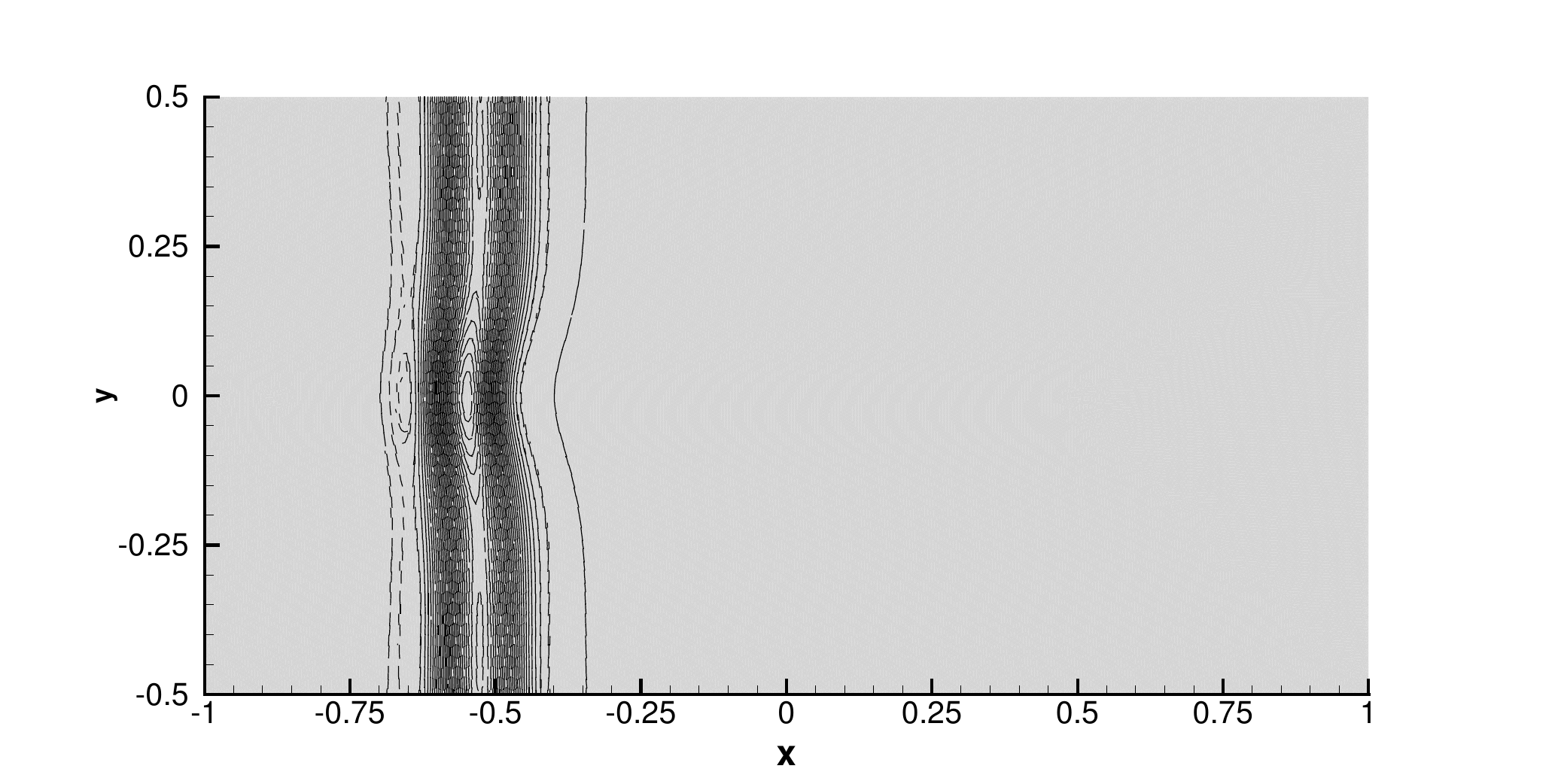}  &          
			\includegraphics[width=0.49\textwidth]{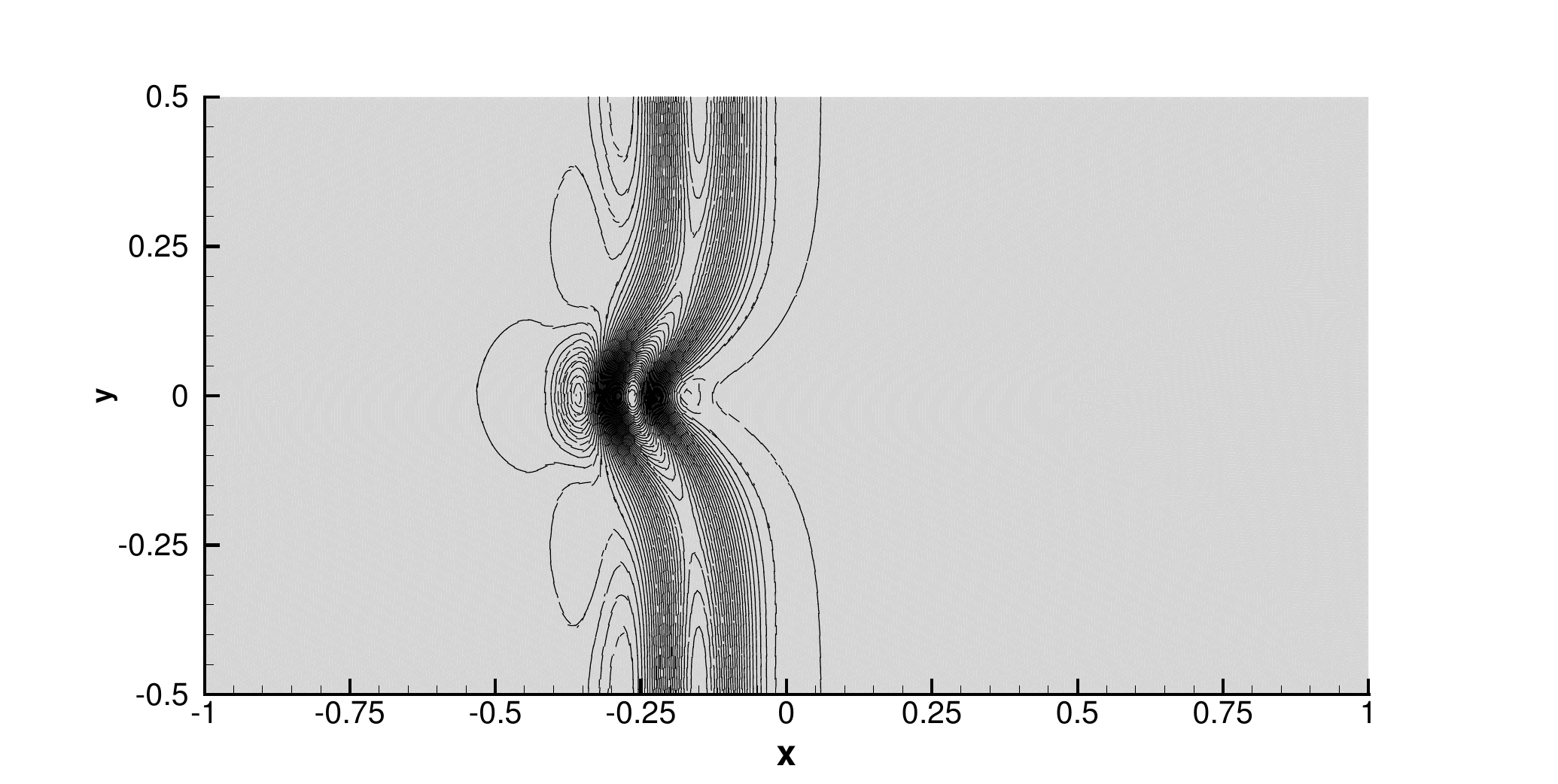}    \\
			\includegraphics[width=0.49\textwidth]{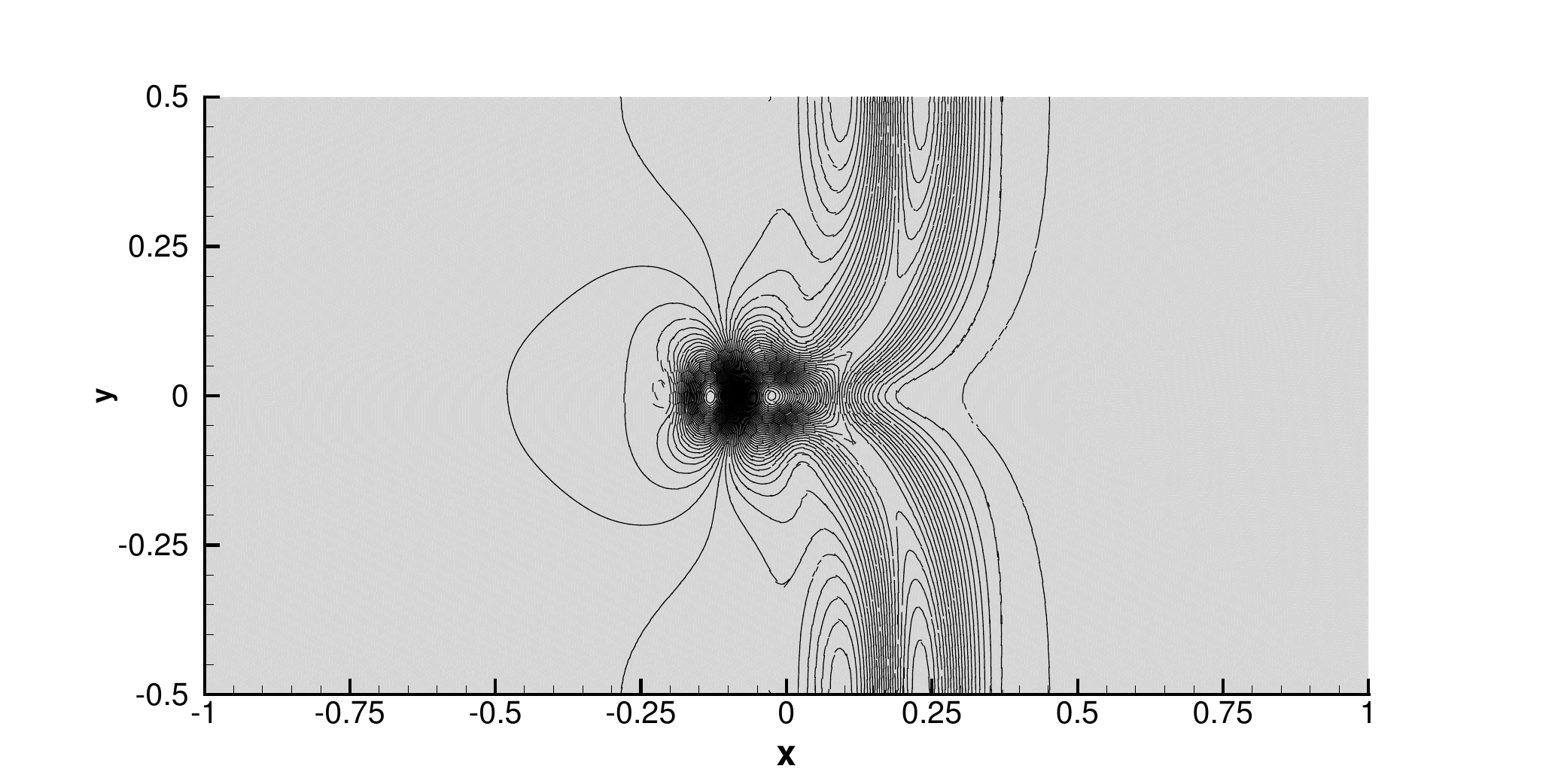}  &          
			\includegraphics[width=0.49\textwidth]{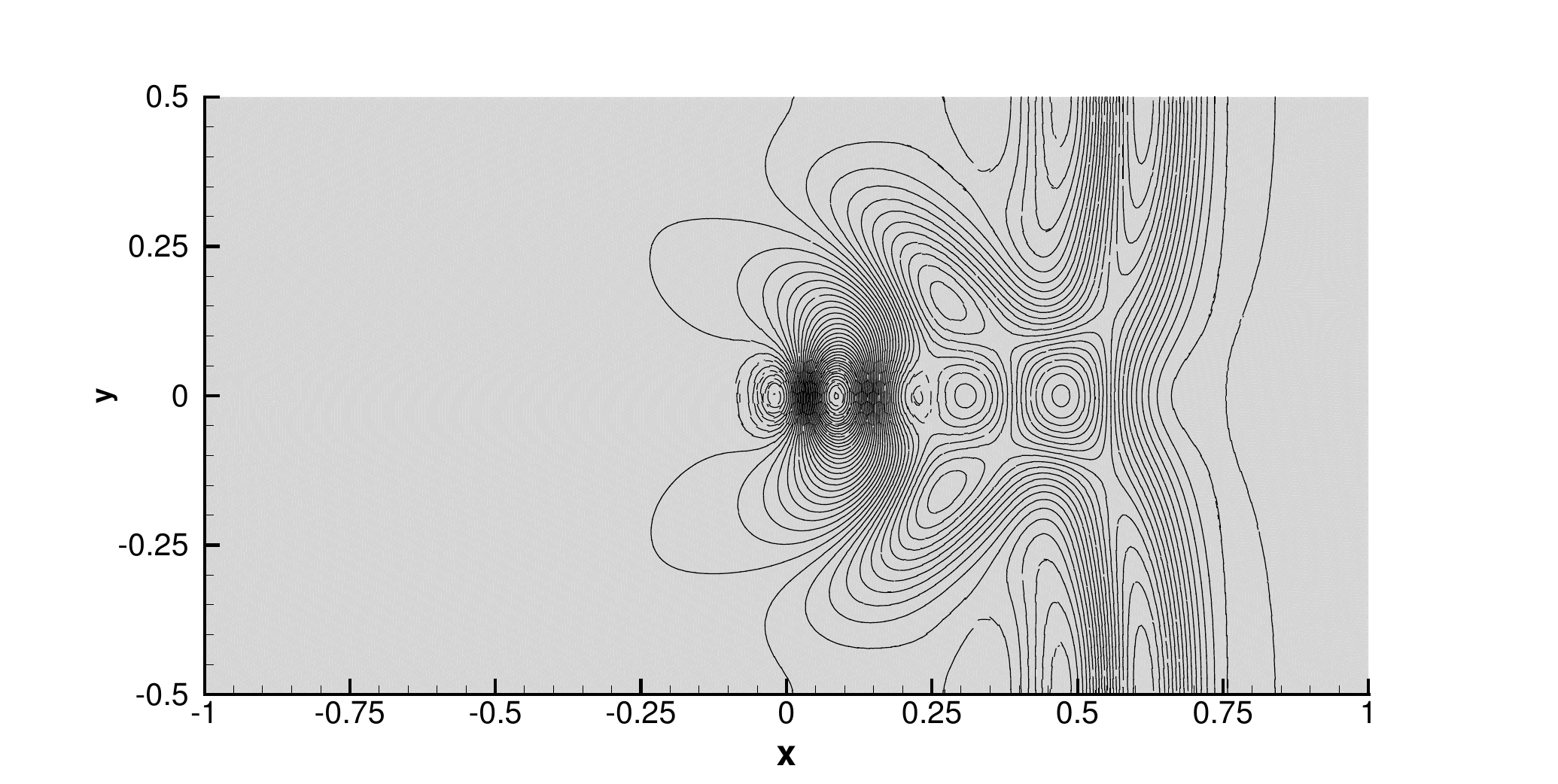}    \\
		\end{tabular}
		\caption{Well-balance test with small perturbation of the free surface ($\delta=10^{-2}$). 80 equidistant contour lines in the interval $\eta=[0.993;1.008]$ are shown at output times $t=0.12$, $t=0.24$, $t=0.36$ and $t=0.48$ (from top left to bottom right panel).}
		\label{fig.Ctest2}
	\end{center}
\end{figure}

\subsection{Circular dambreak}
In order to simulate 2D problems with shock waves, let us consider the circular dambreak problem over a bottom step forwarded in \cite{TavelliSWE2014,StagDG_Dumbser2013}. The computational domain is given by the circle $\Omega=\{\xx \in \mathds{R}^2 \, \, | \, \, r=|\xx|\leq 2\}$ with Dirichlet boundary conditions everywhere. The computational grid counts a total number of $N_P=34477$ with characteristic mesh size $h=1/50$, and the following initial condition is considered:
\begin{equation}
	\eta(\xx,0) = \left\{ \begin{array}{cc}
		1.0 & \textnormal{ if } r \leq 1 \\
		0.5 & \textnormal{ if } r > 1
	\end{array}\right., \qquad b(\xx) = \left\{ \begin{array}{cc}
	0.2 & \textnormal{ if } r \leq 1 \\
	0.0 & \textnormal{ if } r > 1
\end{array}\right., \qquad \vv(\xx,0) = \mathbf{0}.
\end{equation}
The final time of the simulation is $t_f=0.2$, at which the solution exhibits a contact wave traveling towards the center of the domain, as well as a shock wave that is heading the outer boundary. Furthermore, due to the presence of the bottom step, an additional discontinuity is present in the flow at $r=1$. The results are depicted in Figure \ref{fig.dambreak} together with a comparison against the reference solution, which has been computed by solving the one-dimensional shallow water equations in radial direction with geometric reaction source terms, using a classical shock capturing MUSCL-TVD finite volume scheme with 10000 cells \cite{ToroBook}. Overall one can appreciate a very good matching between numerical and reference solution, and no spurious oscillations occur in the plateau between the two shocks. We underline that numerical dissipation is only present in the CWENO finite volume solver for the convective terms and not in the pressure Poisson solver as needed in \cite{Busto_SWE2022}. 

\begin{figure}[!htbp]
	\begin{center}
		\begin{tabular}{cc}
			\multicolumn{2}{c}{\includegraphics[width=0.65\textwidth]{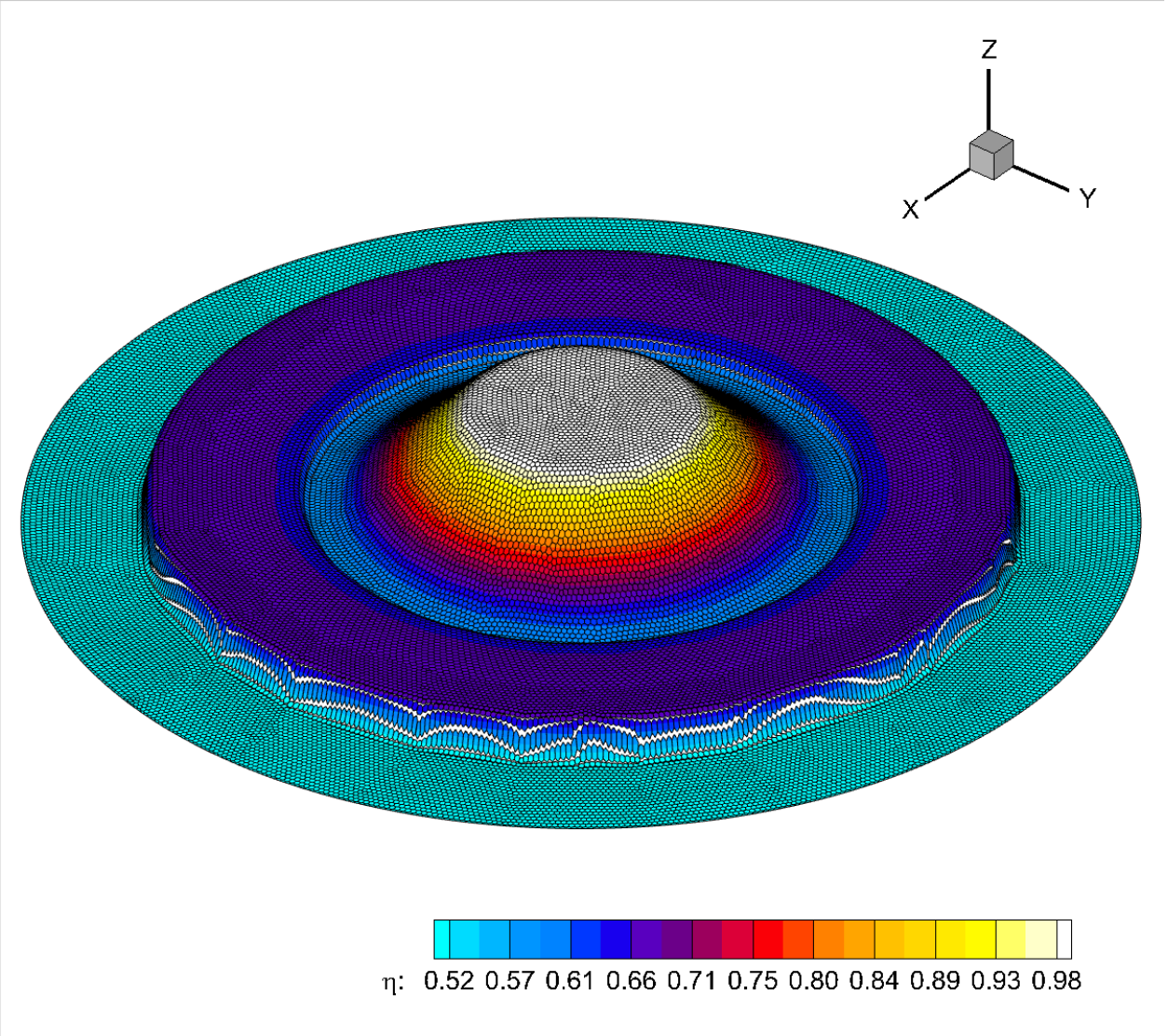}}  \\   	   
			\includegraphics[width=0.47\textwidth]{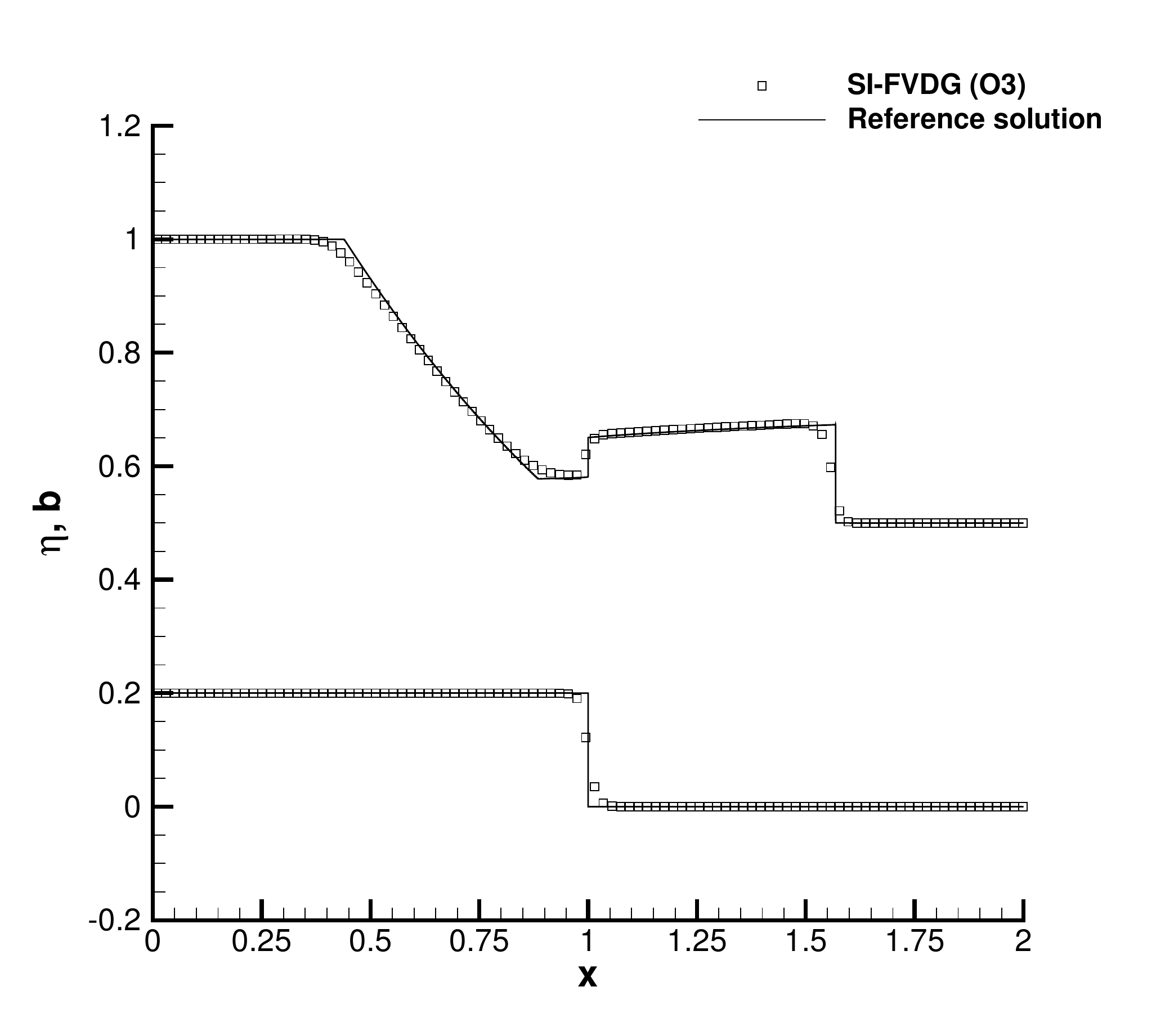} &  			\includegraphics[width=0.47\textwidth]{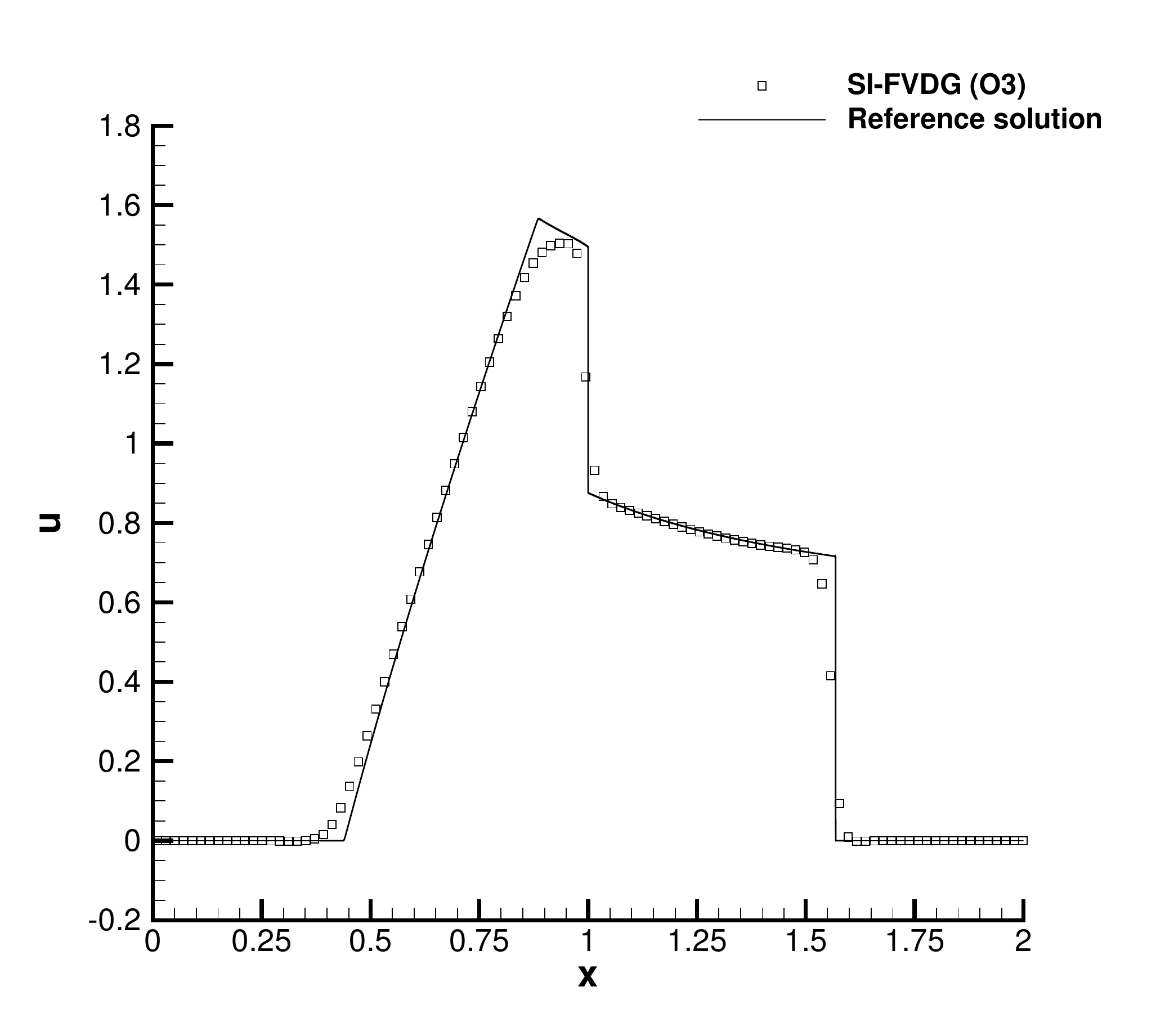} \\
		\end{tabular}
		\caption{Circular dambreak problem at time $t_f=0.2$. Top: three-dimensional view of the free surface elevation. Bottom: one-dimensional cut along the line $y=0$ of the numerical solution compared against the reference solution for the free surface and bottom profile (left) as well as for the horizontal velocity component (right).}
		\label{fig.dambreak}
	\end{center}
\end{figure}

\subsection{Riemann problems}
The SI-FVDG schemes are validated against a set of Riemann problems which take into account flat and variable bottom topography. The exact solution is computed relying on the Riemann solver presented in \cite{ToroBookSWE} and \cite{ToroSWERP} for flat and variable bottom, respectively. The initial condition is given in terms of two states $\Q_L=(\eta_L,u_L,b_L)$ and $\Q_R=(\eta_R,u_R,b_R)$ separated by a discontinuity located at position $x=x_d$:
\begin{equation}
	\Q(\xx,0) = \left\{ \begin{array}{ccc} \Q_L & \textnormal{if} & x \leq x_d \\
		\Q_R & \textnormal{if} & x > x_d
	\end{array}\right. .
\end{equation}
Table \ref{tab:initRP} summarizes the extension of the computational domain as well as the initial condition for free surface elevation, horizontal velocity and bottom elevation for four Riemann problems. 

\begin{table}[!htbp]  
	\caption{Initialization of Riemann problems. Initial states left (L) and right (R) are reported as well as the final time of the simulation $t_f$, the computational domain $[x_L;x_R]$, the position of the initial discontinuity $x_d$ and the characteristic mesh size $h$.}  
	\begin{center} 
		\begin{small}
			\renewcommand{\arraystretch}{1.0}
			\begin{tabular}{llllllllllll} 
				\hline
				Test & $\eta_L$ & $u_L$ & $b_L$ & $\eta_R$ & $u_R$ & $b_R$ & $x_L$ & $x_R$ & $x_d$ & $h$ & $t_{f}$\\
				\hline
				RP1 \cite{ToroBookSWE}  & 1 & 0 & 0 & 2 & 0 & 0 & -0.5 & 0.5 & 0 & 1/200 & 0.075\\
				RP2                 & $10^3$ & 0 & 0 & 1 & 0 & 0 & -15 & 15 & 0 & 1/200 & 0.09\\ 
				RP3 \cite{ToroSWERP} & 1 & 0 & 0.2 & 0.5 & 0 & 0 & -5 & 5 & 0 & 1/200 & 1\\
				RP4 \cite{ToroSWERP}  & 1.46184 & 0 & 0 & 0.30873 & 0 & 0.2 & -0.5 & 0.5 & 0 & 1/200 & 1\\ 
				\hline
			\end{tabular}
		\end{small}
	\end{center}
	\label{tab:initRP}
\end{table}

Despite the one-dimensional setup of these test cases, the computational domain is given by $\Omega=[x_L;x_R] \times [x_L/10;x_R/10]$ and it is discretized with an unstructured Voronoi mesh of size $h$, hence making the computation intrinsically multidimensional. The results are depicted in Figure \ref{fig.RP} where the numerical solution is compared against the reference solution through a one-dimensional cut of 200 equidistant points along the $x-$axis of the computational domain at $y=0$. The first two Riemann problems (RP1 and RP2) assume a constant flat bathymetry, and deal with shock and rarefaction waves. The remaining Riemann problems (RP3 and RP4) deal with a jump in the bottom elevation of height $\Delta b=0.2$, which is responsible of the generation of shock waves. Overall an excellent agreement can be noticed, demonstrating that the novel SI-FVDG schemes can also handle supercritical flows with Froude numbers greater than one, namely for RP2 the maximum Froude number is $\Fr=5.73$. This is achieved thanks to the very robust CWENO finite volume scheme for the discretization of the nonlinear convective terms. The implicit treatment of the free surface elevation is enough to guarantee a stable scheme for all the four Riemann problems. The moving shocks are correctly captured as well as the values of the plateau between two discontinuities, confirming that the SI-FVDG schemes are conservative by construction. Finally, the one-dimensional symmetry of the problem is perfectly retrieved even in the context of arbitrary shaped polygonal cells, as confirmed by the three-dimensional views of the free surface elevation in Figure \ref{fig.RP}.

\begin{figure}[!htbp]
	\begin{center}
		\begin{tabular}{ccc} 
			\includegraphics[width=0.33\textwidth]{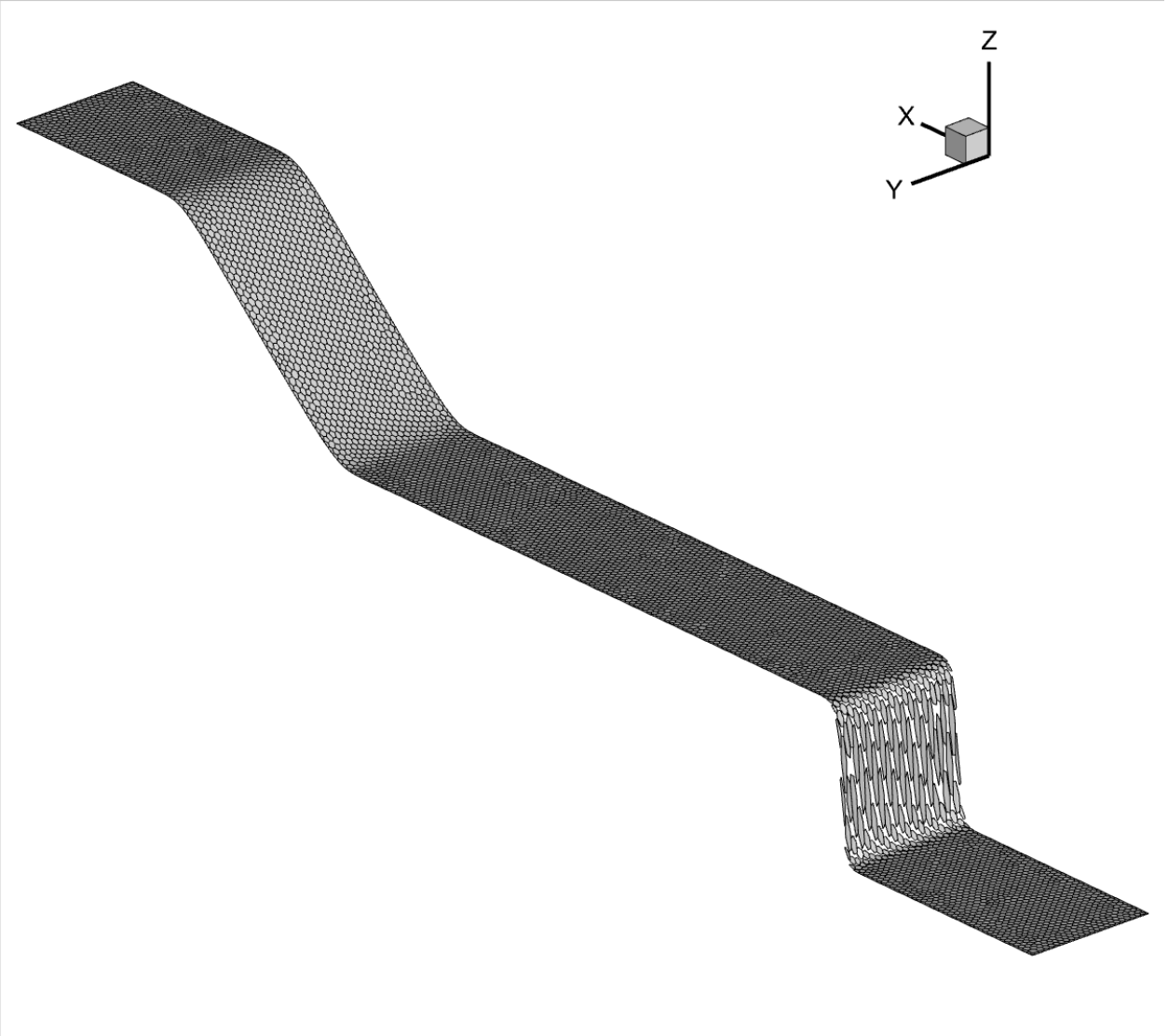} &  
			\includegraphics[width=0.33\textwidth]{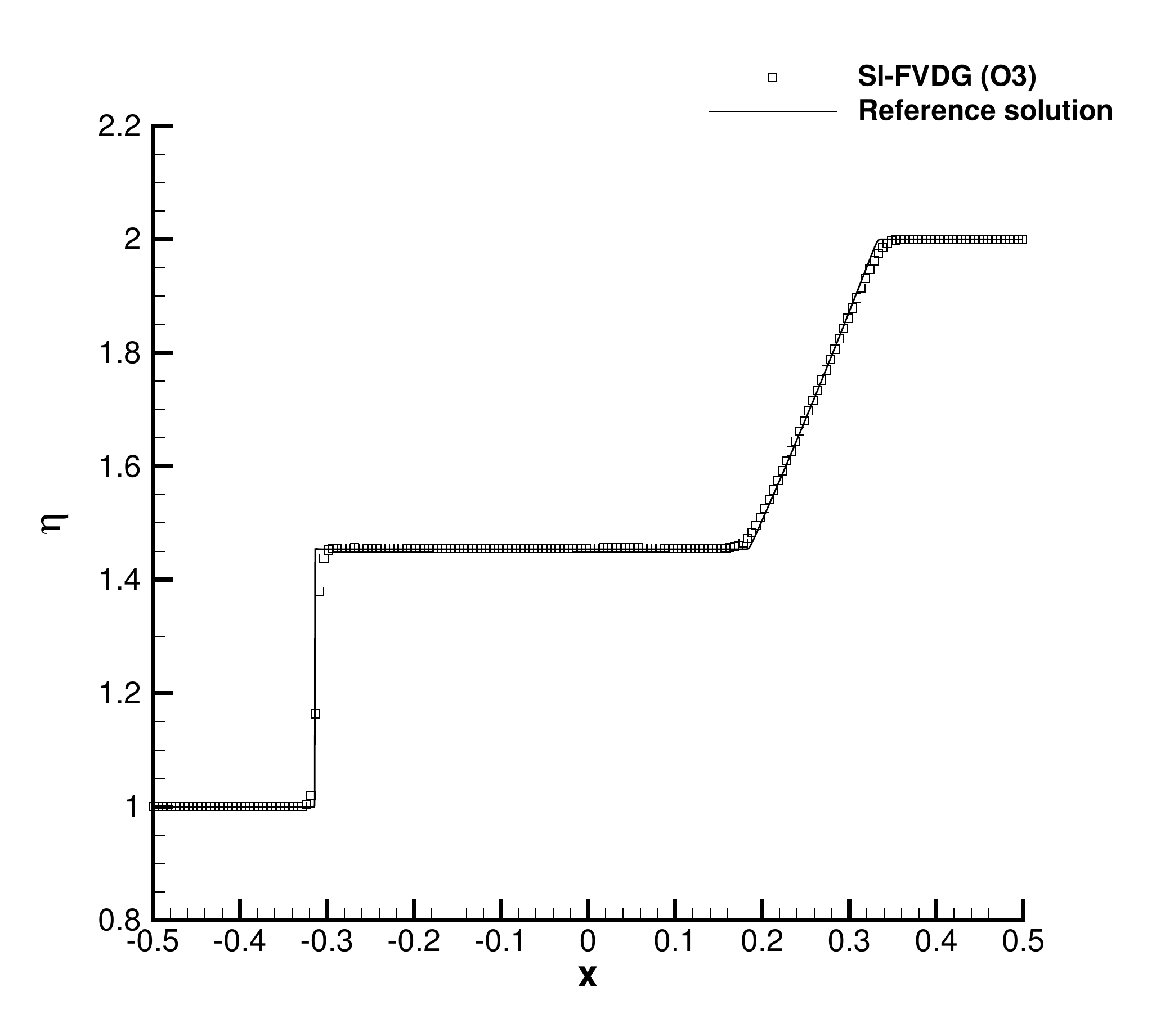} &  			\includegraphics[width=0.33\textwidth]{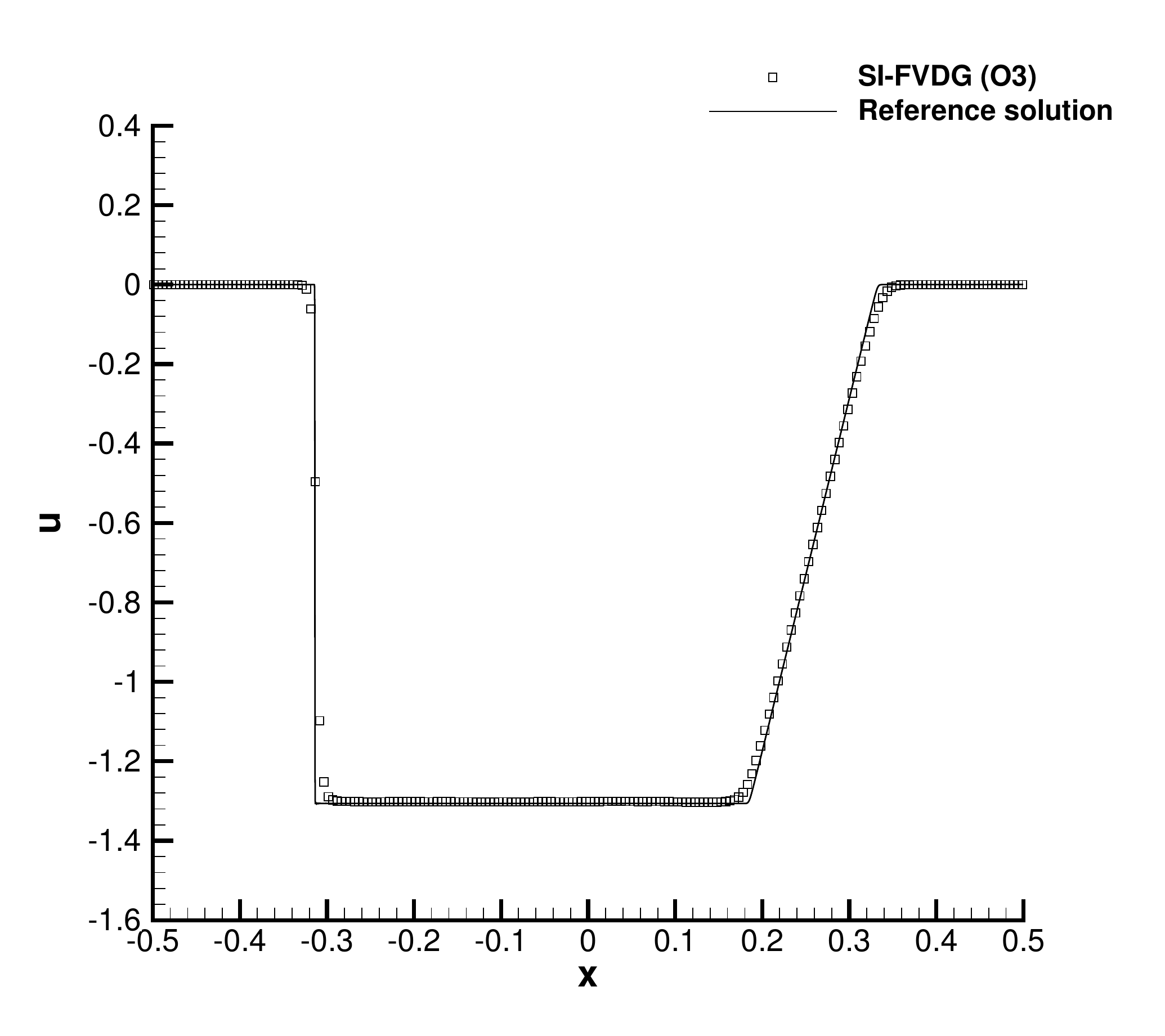} \\
			\includegraphics[width=0.33\textwidth]{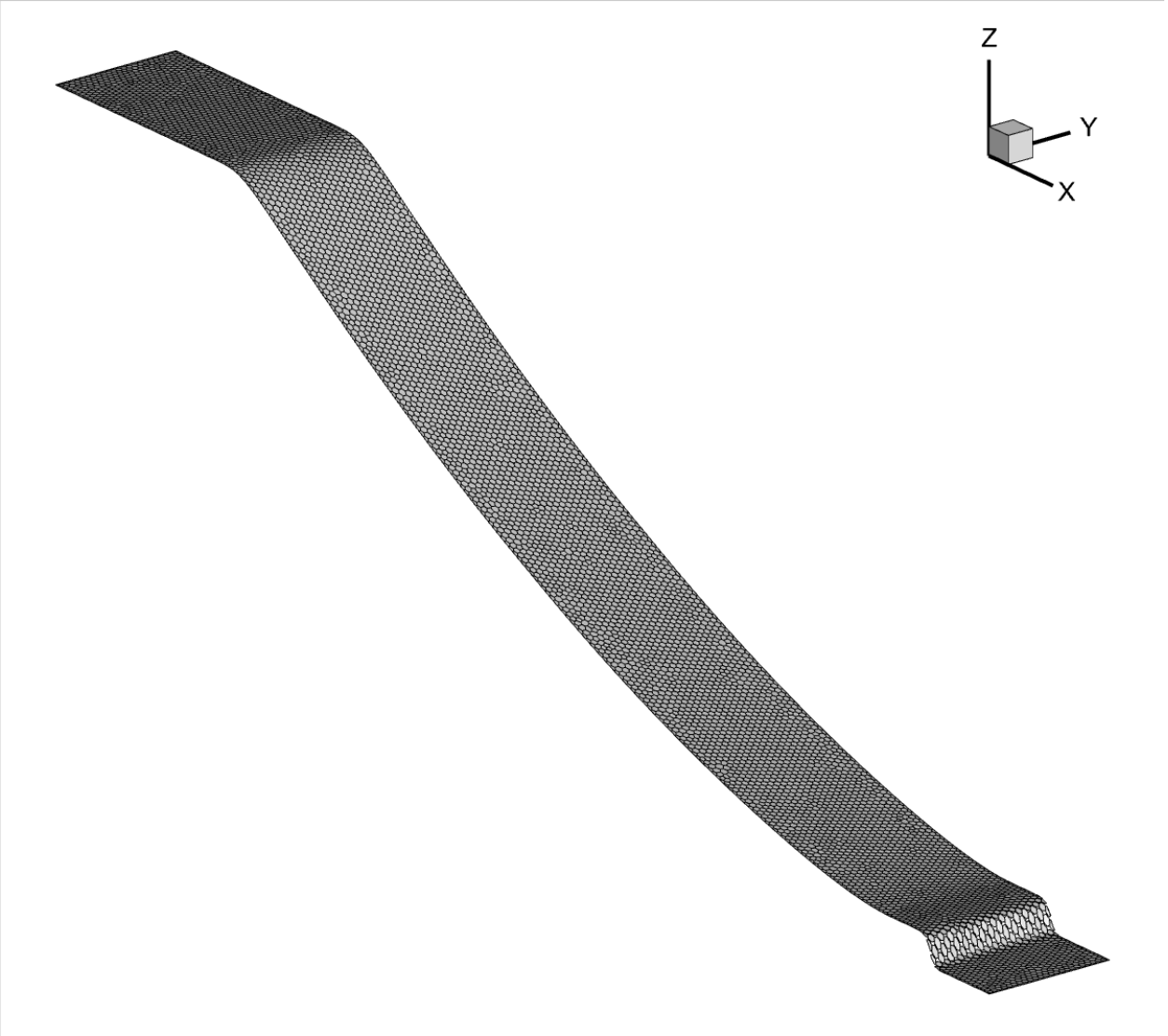} &  
			\includegraphics[width=0.33\textwidth]{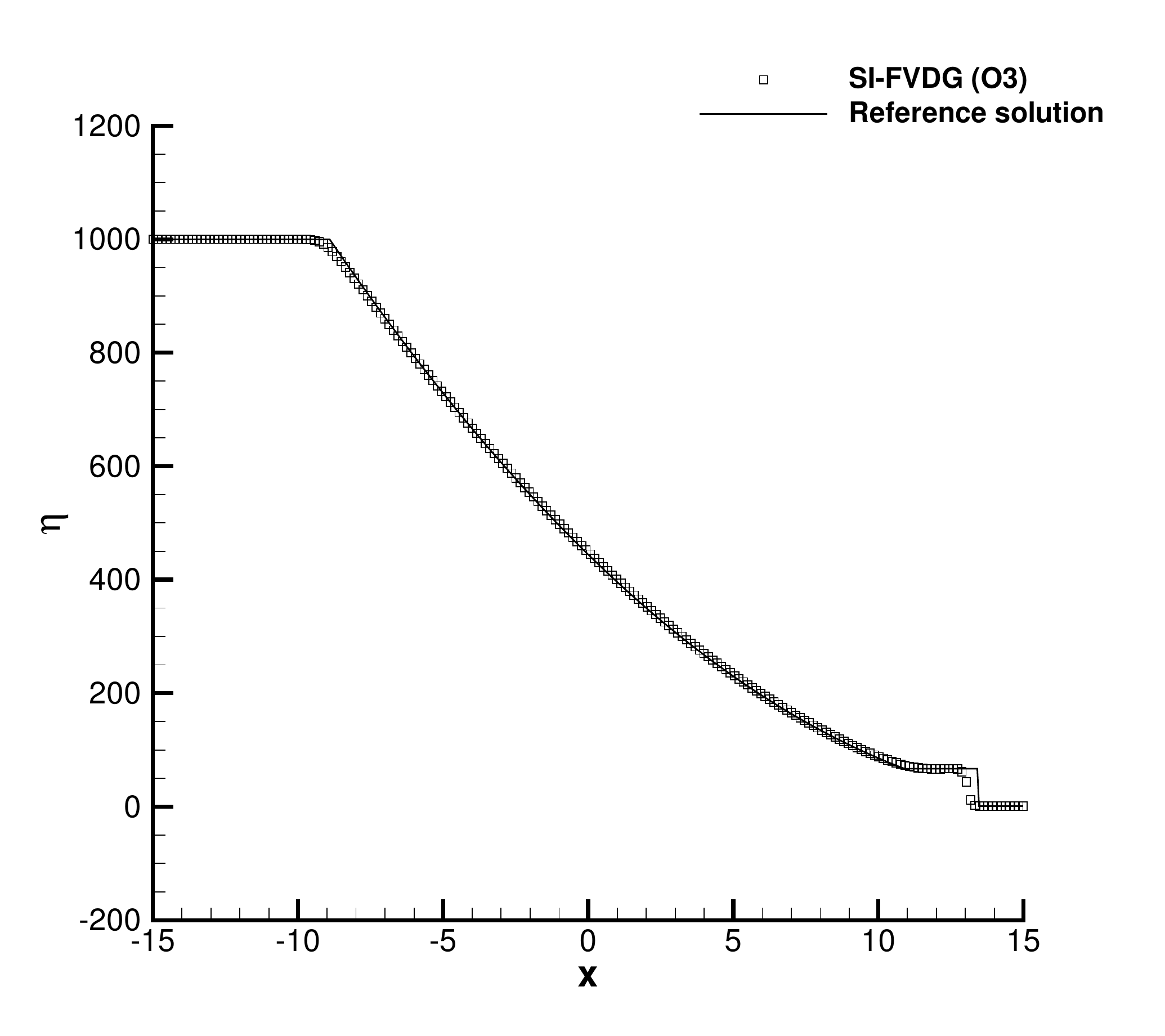} &  			\includegraphics[width=0.33\textwidth]{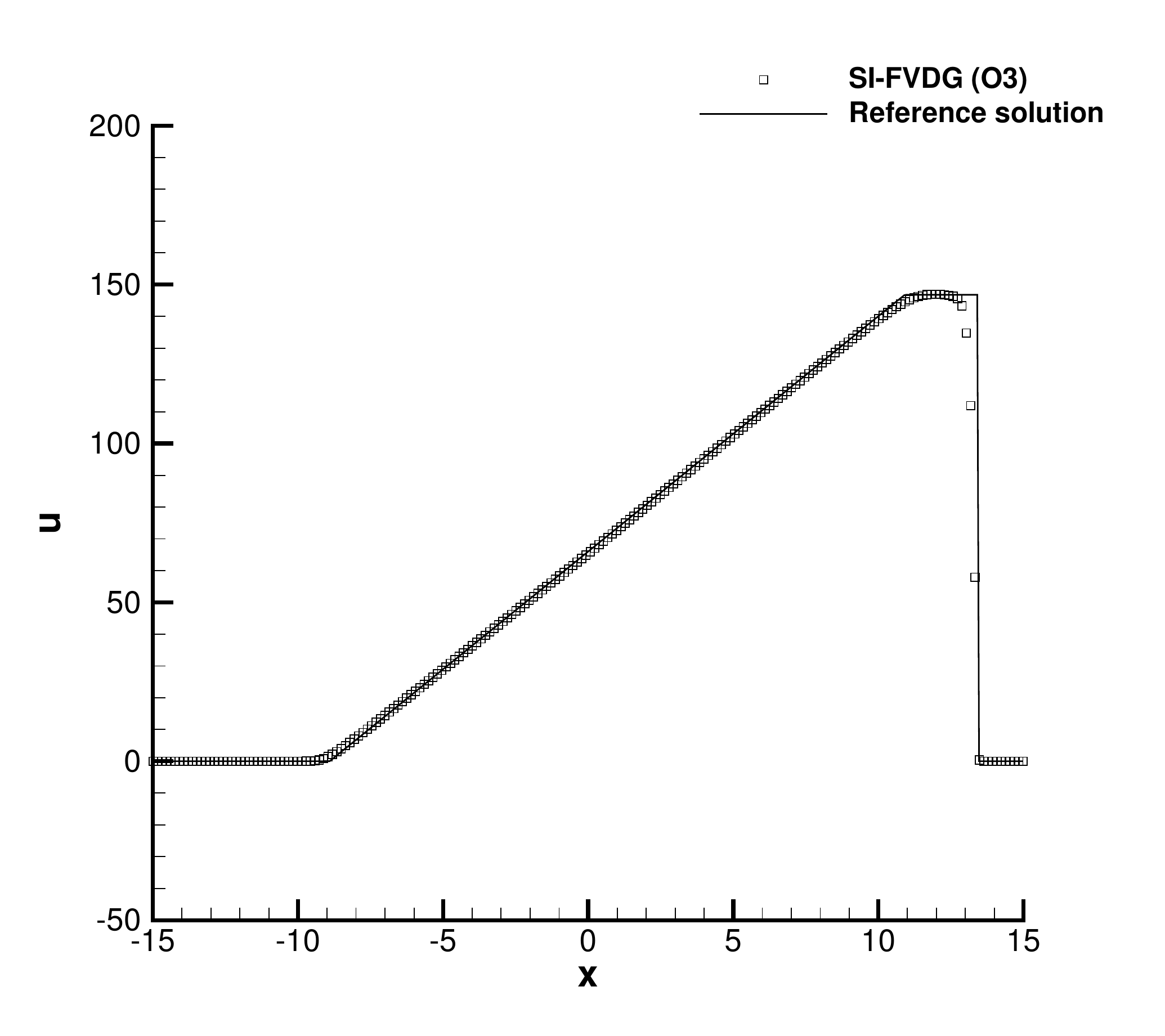} \\
			\includegraphics[width=0.33\textwidth]{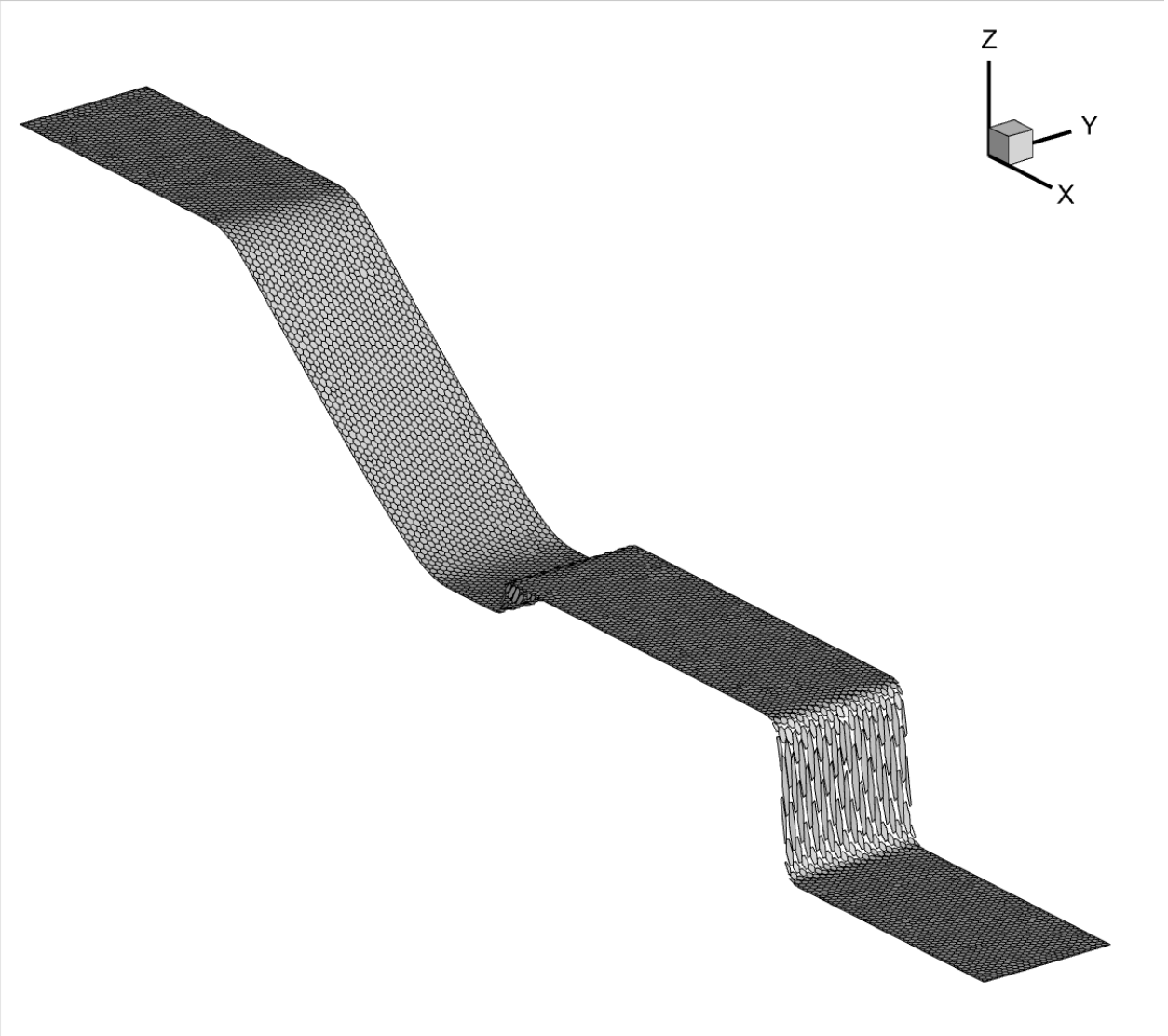} &  
			\includegraphics[width=0.33\textwidth]{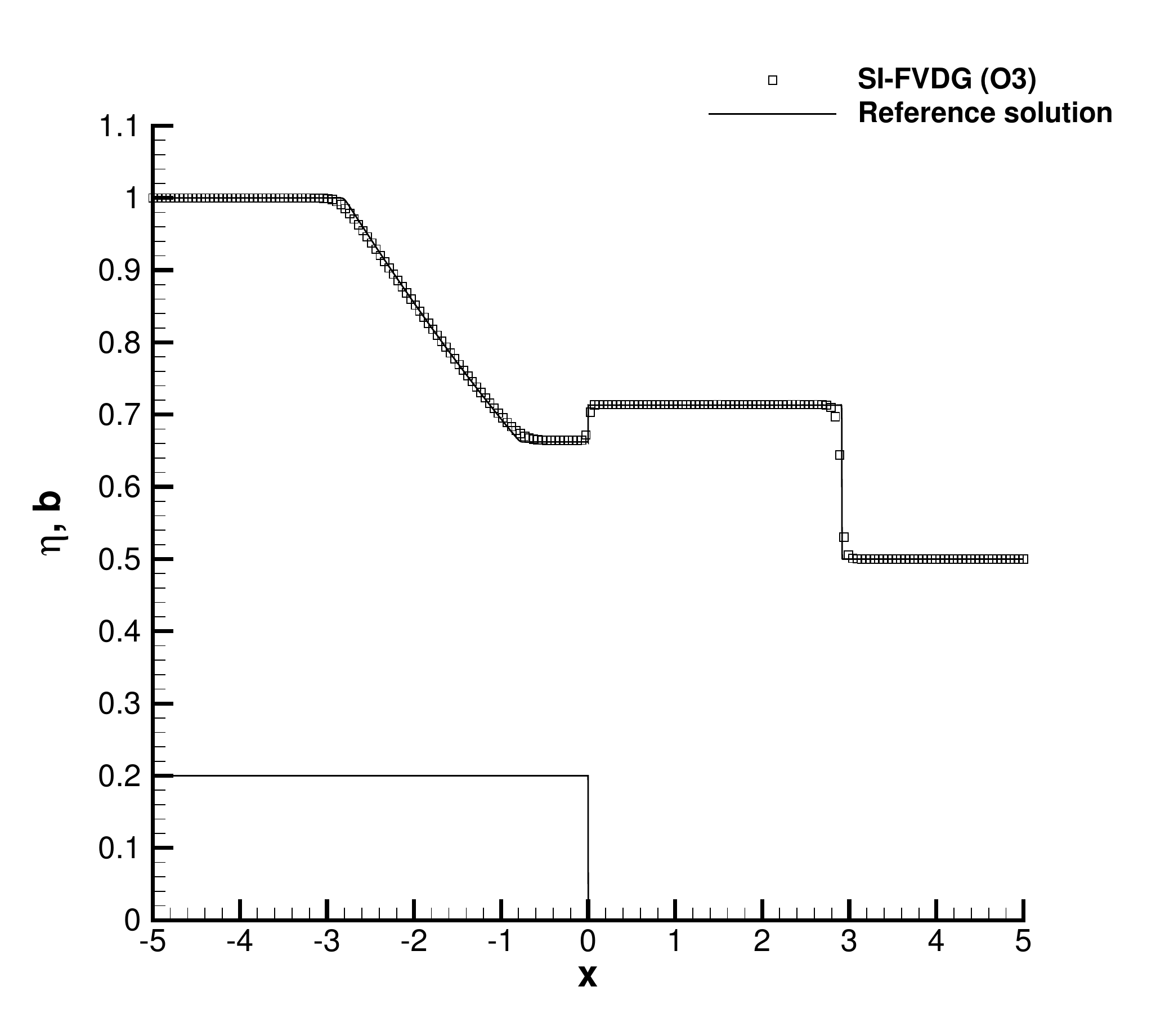} &  			\includegraphics[width=0.33\textwidth]{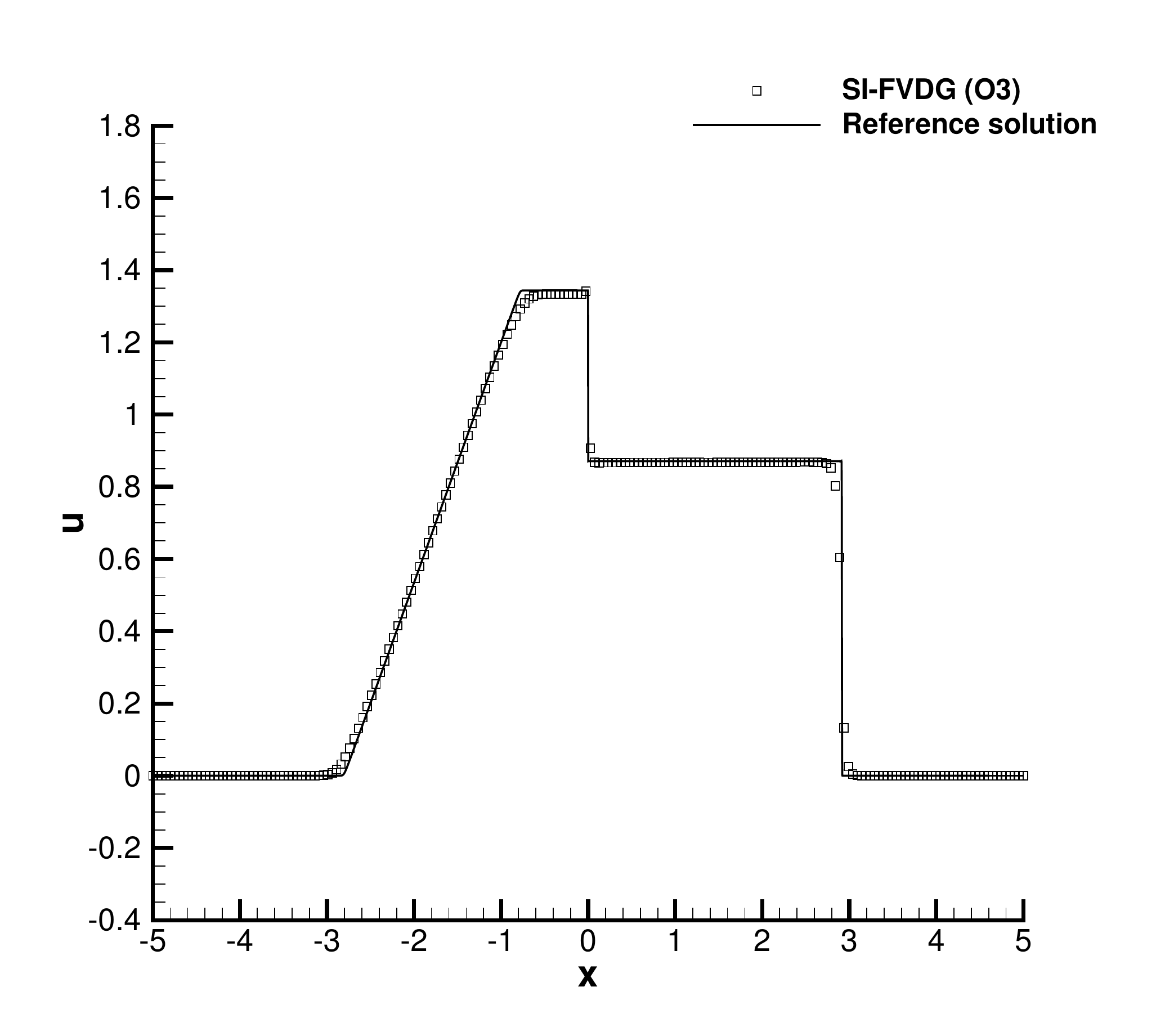} \\
			\includegraphics[width=0.33\textwidth]{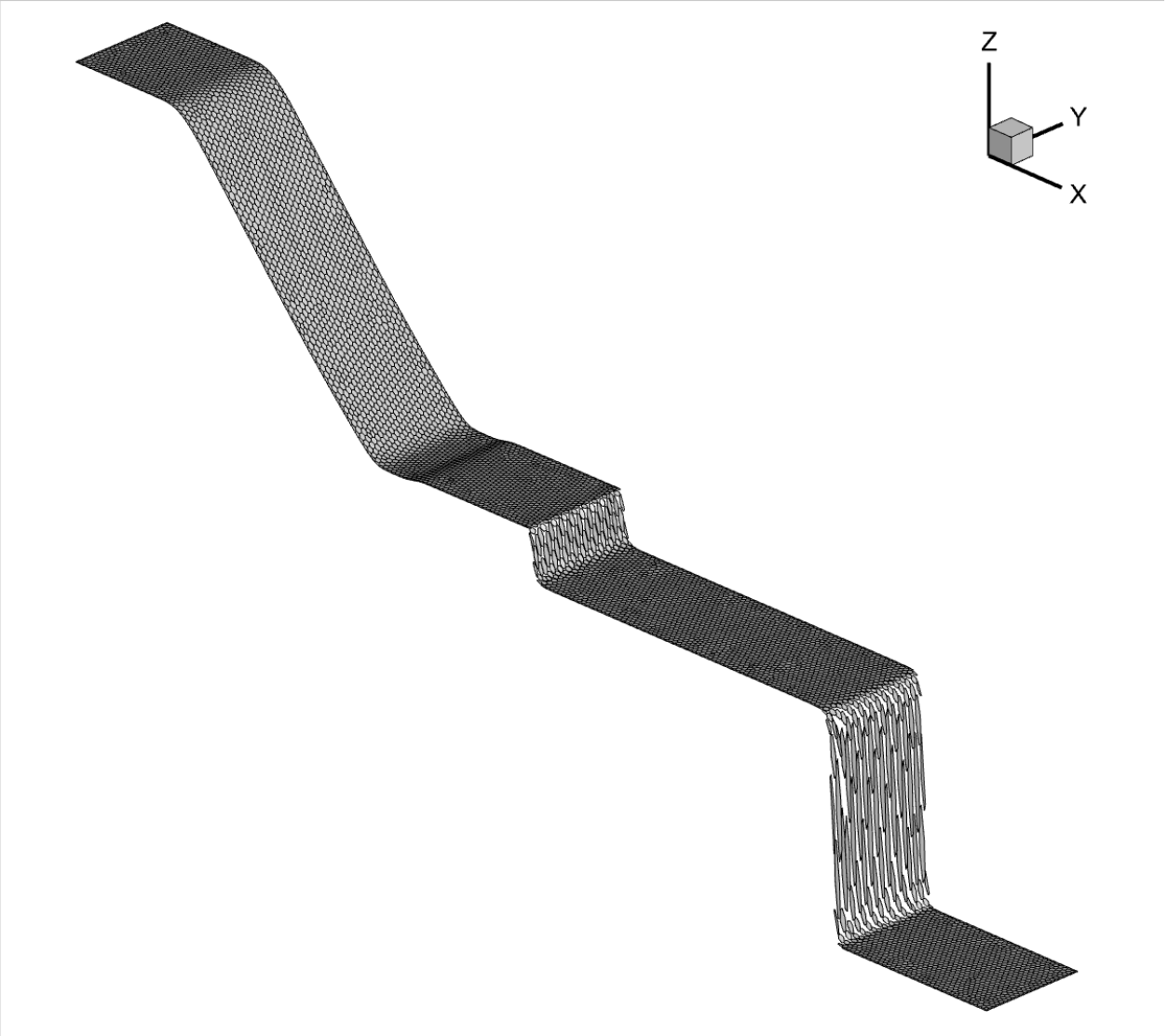} &  
			\includegraphics[width=0.33\textwidth]{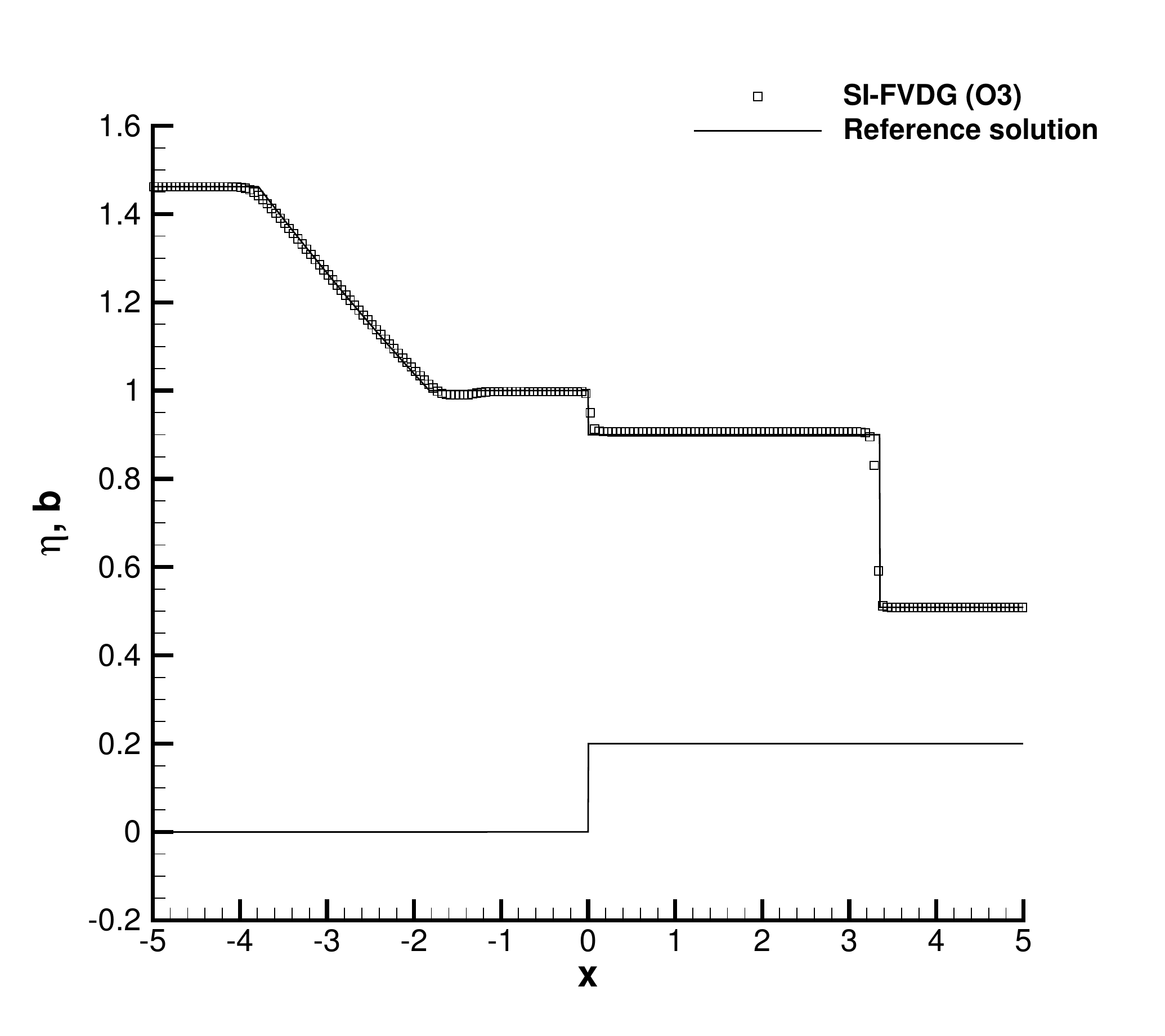} &  			\includegraphics[width=0.33\textwidth]{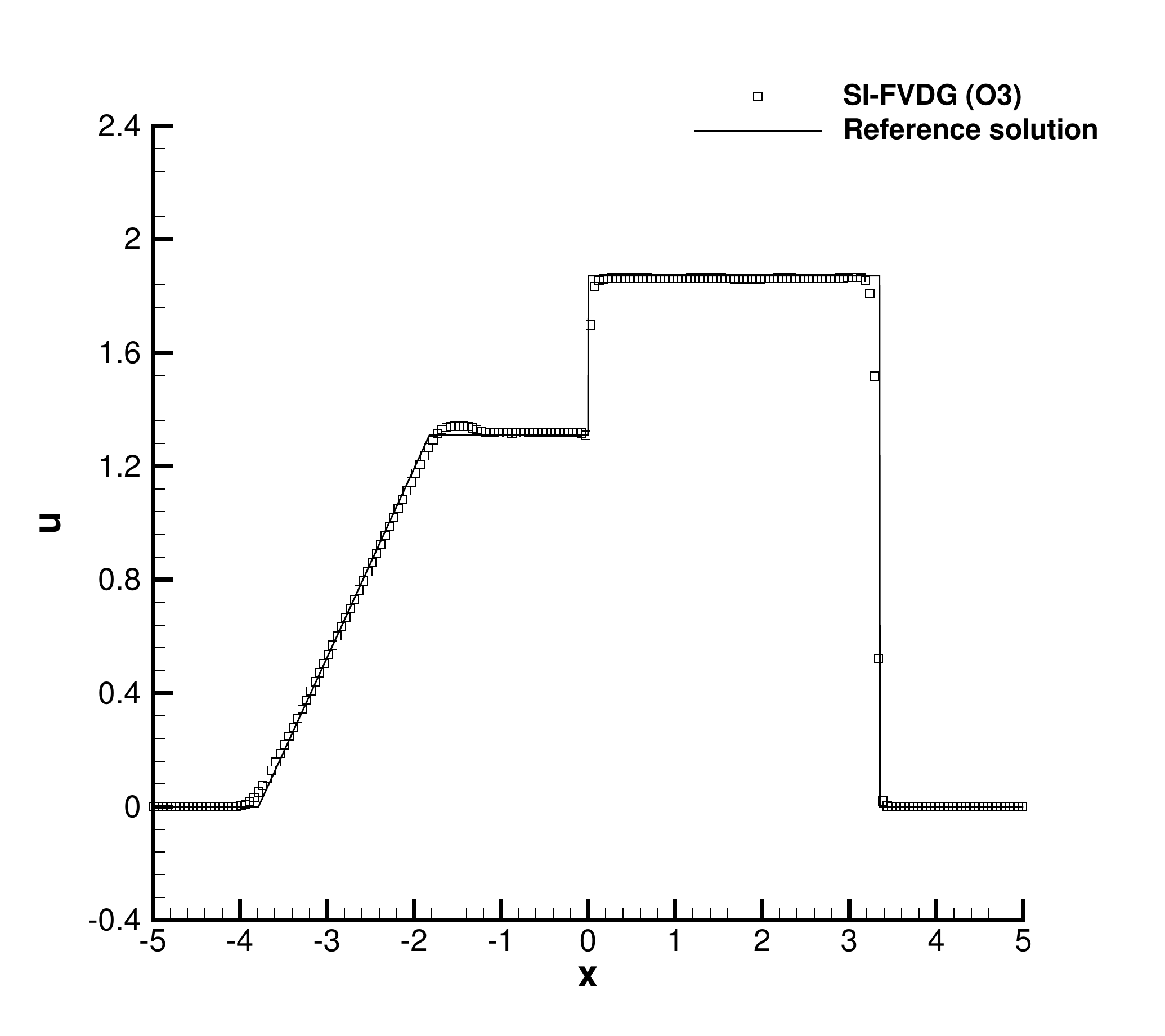} \\
		\end{tabular}
		\caption{Riemann problems RP1, RP2, RP3 and RP4 (from top to bottom row). Three-dimensional view of the free surface elevation with the Voronoi computational mesh (left) and comparison against the reference solution for the variables $\eta$ (middle) and $u$ (right).}
		\label{fig.RP}
	\end{center}
\end{figure}

\subsection{Smooth surface wave propagation}
Here, a wave propagation of the free surface elevation is considered following the setup presented in \cite{StagDG_Dumbser2013}. The computational domain is the square $\Omega=[-1;1]^2$ with Dirichlet boundary conditions imposed everywhere, which is discretized with a total number of $N_P=15717$ Voronoi cells with characteristic mesh size $h=1/50$. The initial condition is given by
\begin{equation}
	\label{eqn.smoothwave_ini}
	\eta(\xx,0) = 1 + e^{-\frac{1}{2}(r^2/\sigma^2)}, \qquad \vv(\xx,0) = \mathbf{0}, \qquad b(\xx)=0,
\end{equation}
with $\sigma=0.1$. The time step is fixed to $\dt=0.001$ and the final time of the simulation is $t_f=0.15$, so that the wave profile becomes stiff and a shock wave starts. Figure \ref{fig.smoothwave3D} depicts a three-dimensional view of the free surface elevation at different output times, highlighting the capability of the SI-FVDG schemes of maintaining the symmetry of the solution despite the unstructured computational grid. As already done for the circular dambreak problem, the reference solution is computed relying on a one-dimensional MUSCL-TVD scheme on a very fine mesh, and a comparison against the numerical solution for the free surface elevation and the horizontal velocity component is plot in Figure \ref{fig.smoothwaveXY}. An overall good agreement can be observed, especially until time $t=0.1$, when the flow is still smooth. At time $t=0.15$, the shock is smeared by the CWENO reconstruction technique and the SI-FVDG scheme is stable and does not present spurious oscillations.

\begin{figure}[!htbp]
	\begin{center}
		\begin{tabular}{cc} 
			\includegraphics[width=0.40\textwidth]{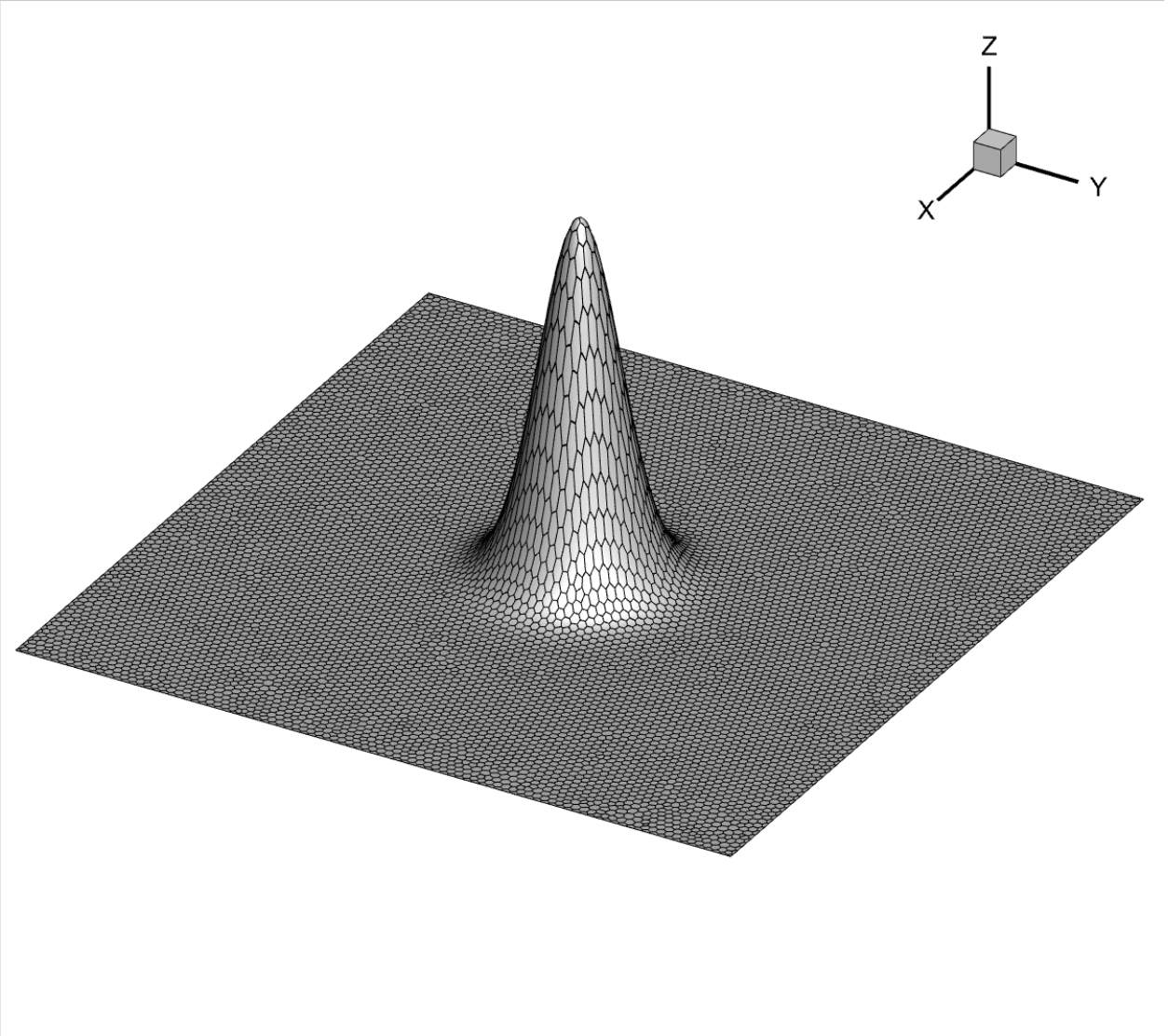} &  			\includegraphics[width=0.40\textwidth]{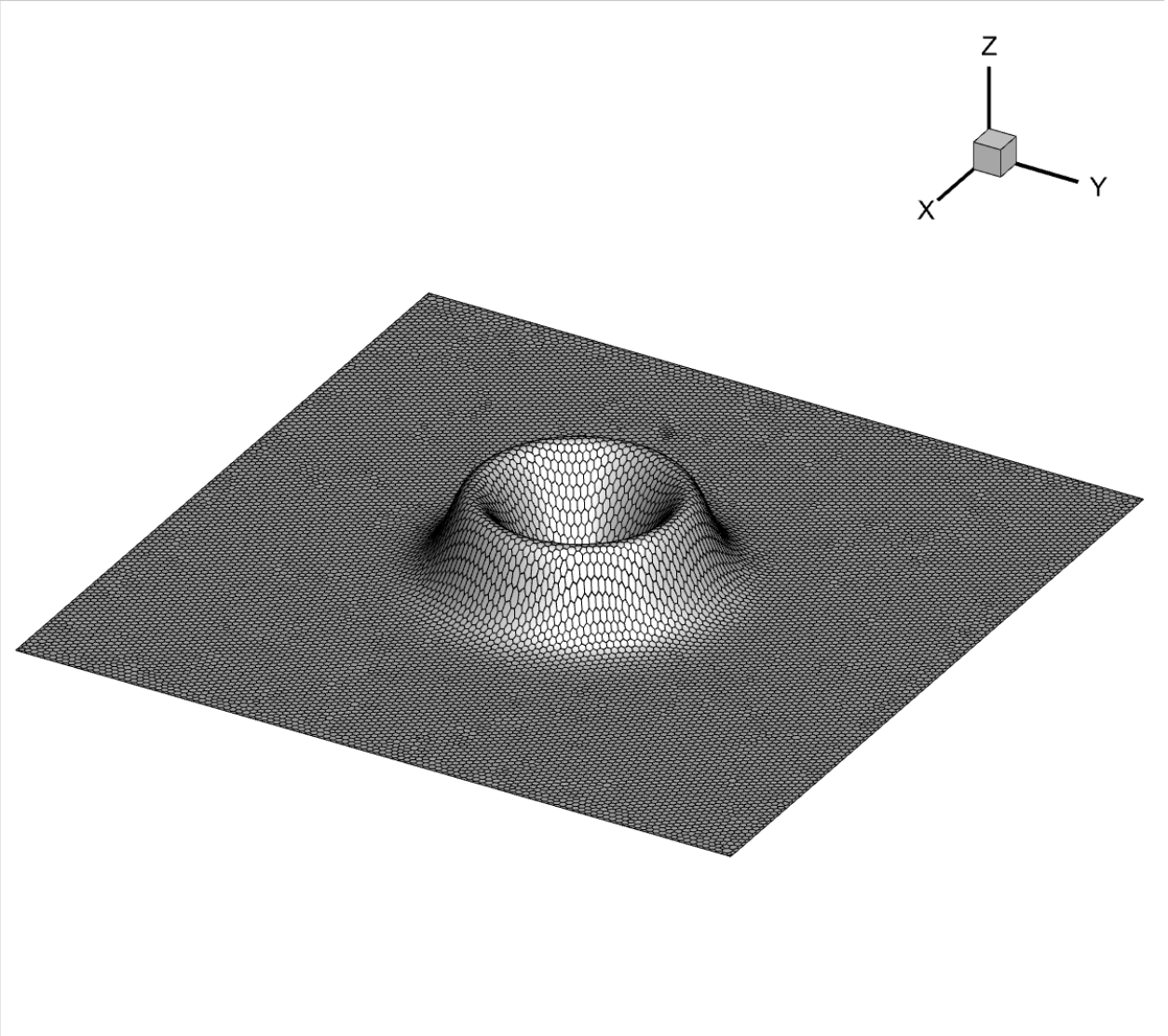} \\
			\includegraphics[width=0.40\textwidth]{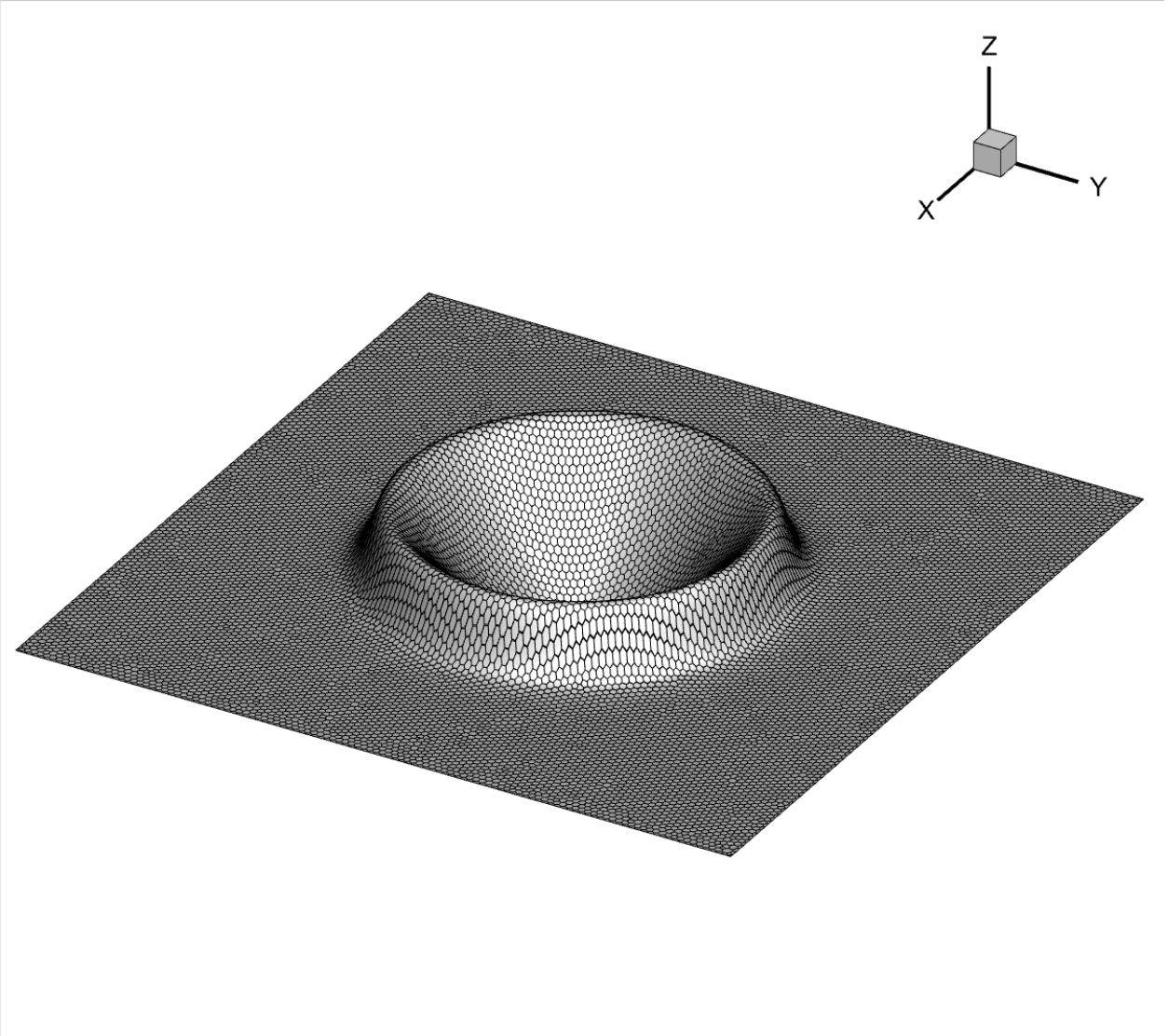} &  			\includegraphics[width=0.40\textwidth]{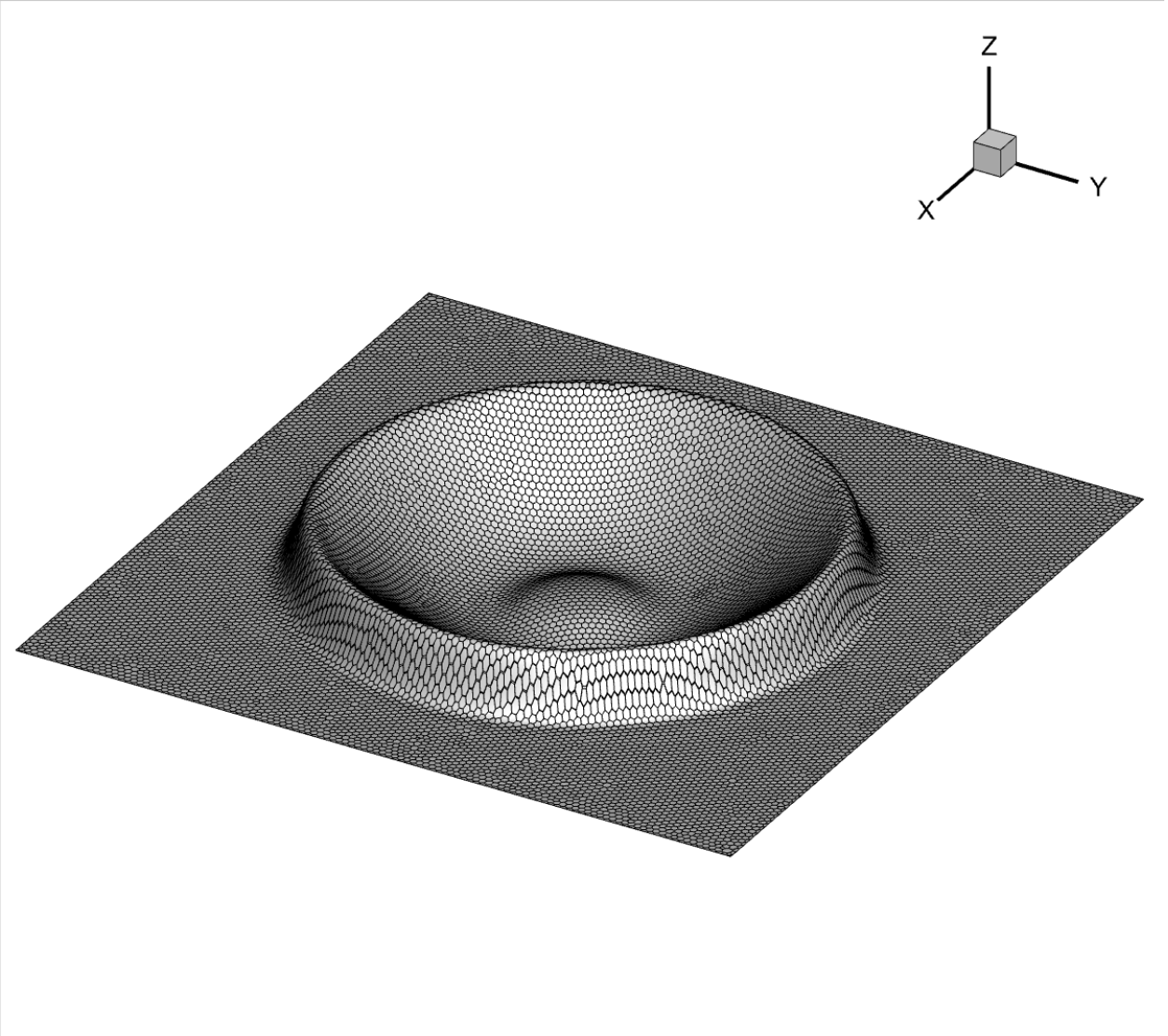} \\
		\end{tabular}
		\caption{Smooth surface wave propagation problem. Three-dimensional view of the free surface elevation and Voronoi computational mesh at output times $t=0$, $t=0.05$, $t=0.1$ and $t=0.15$ (from top left to bottom right panel).}
		\label{fig.smoothwave3D}
	\end{center}
\end{figure}

\begin{figure}[!htbp]
	\begin{center}
		\begin{tabular}{cc} 
			\includegraphics[width=0.47\textwidth]{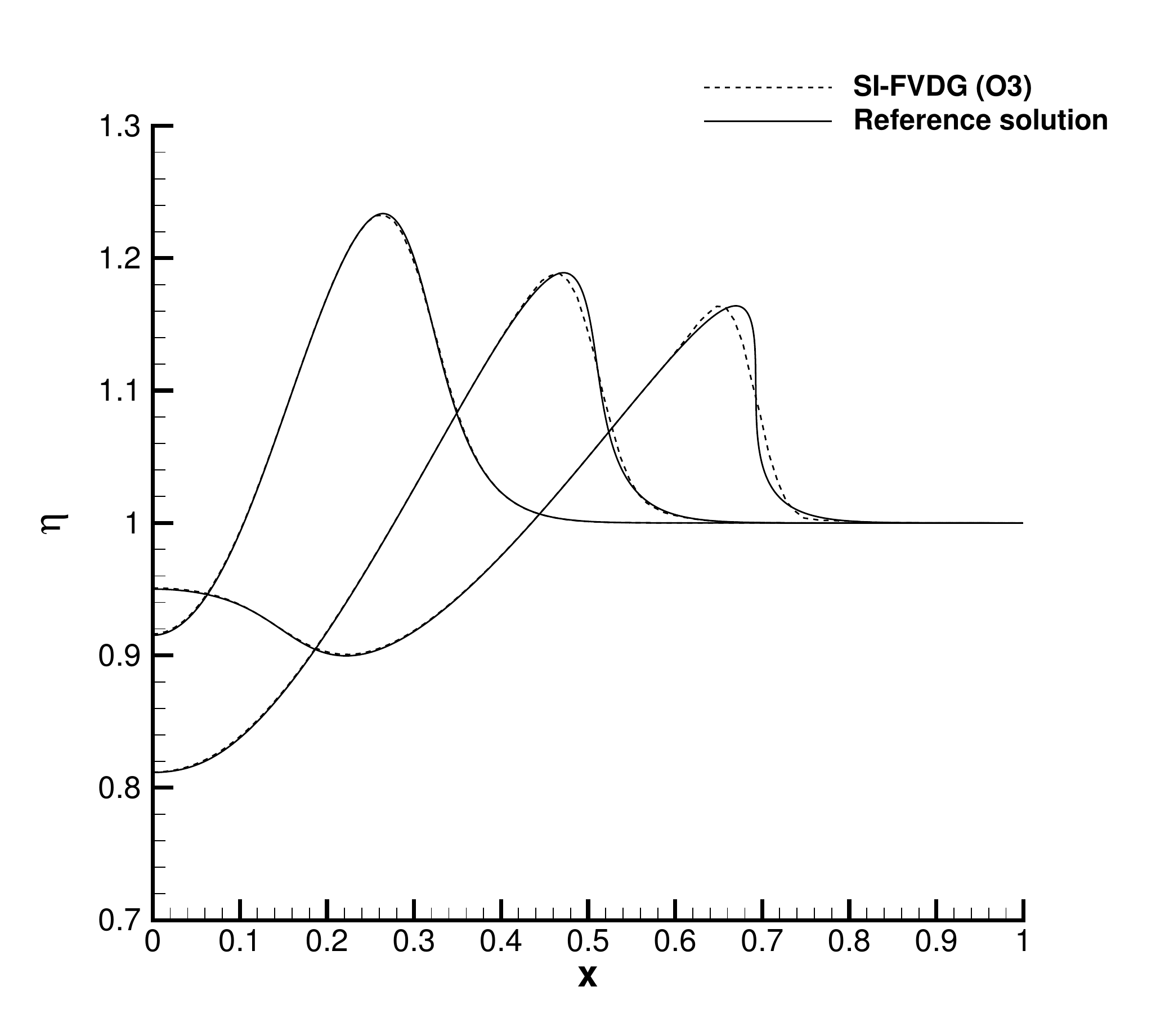} &  			\includegraphics[width=0.47\textwidth]{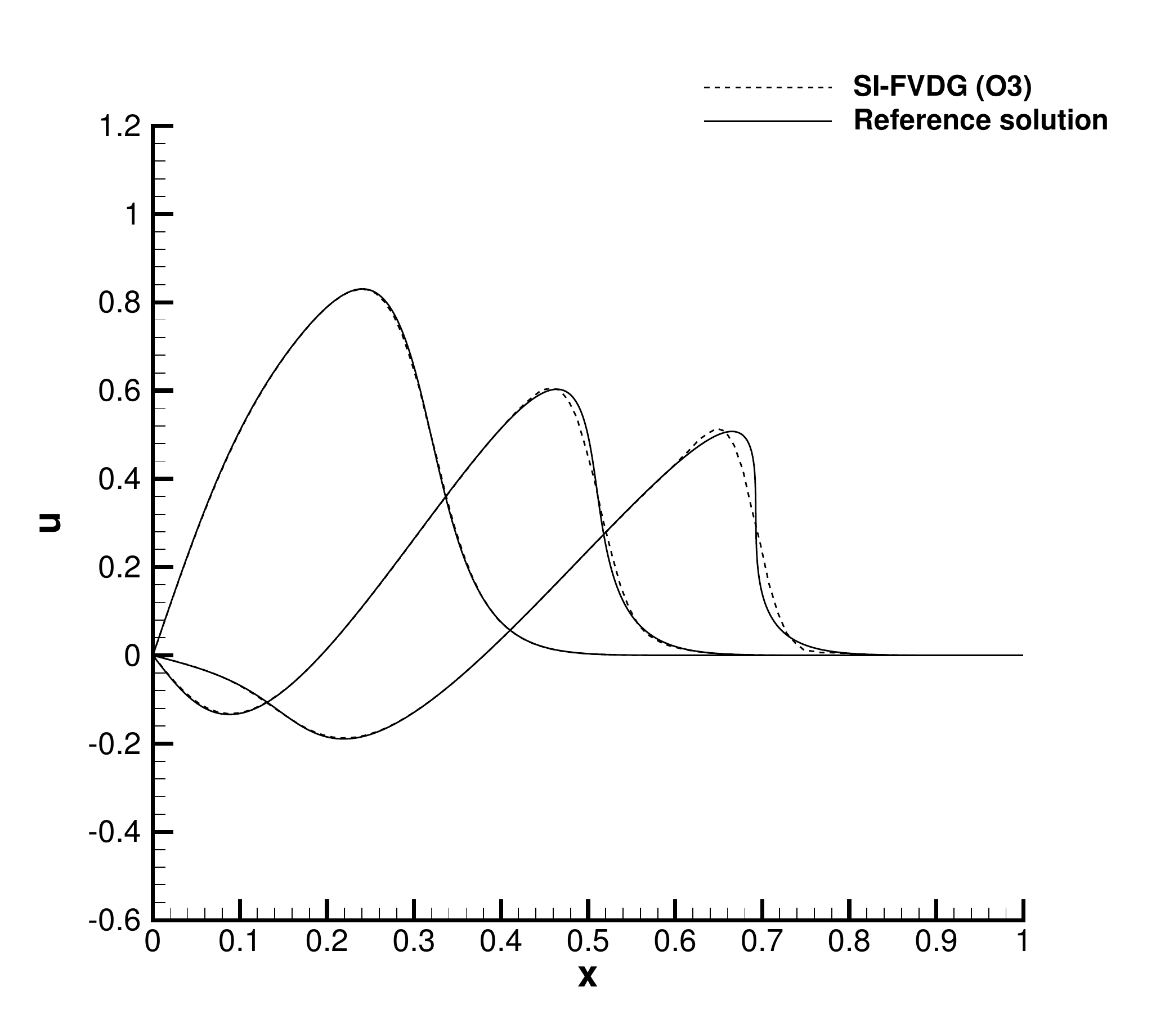} \\
		\end{tabular}
		\caption{Smooth surface wave propagation problem. Comparison between numerical (dashed lines) and reference (solid lines) solution at output times $t=0.05$, $t=0.1$ and $t=0.15$ for the free surface elevation $\eta$ (left) and the horizontal velocity component $u$ (right).}
		\label{fig.smoothwaveXY}
	\end{center}
\end{figure}

\subsection{Low Froude number flow around a circular cylinder}
As a last test case we propose to simulate a low Froude flow with $\Fr=3.19 \cdot 10^{-3}$ that passes around a circular cylinder of radius $r_c=1$ \cite{Busto_SWE2022,BassiRebay97}. The computational domain is $\Omega=[-16;16]^2 \symbol{92} {\xx \in \R^2 | r \leq r_c}$, with the generic radial coordinate given by $r=\sqrt{x^2+y^2}$, and the bottom is assumed to be flat ($b=0$). The mesh counts a total number of $N_P=15516$ and it is made of Voronoi cells with characteristic mesh size of $h=1/20$ close to the border of the cylinder which regularly increase their diameter until $h=1/2$ on the domain boundaries, see Figure \ref{fig.Cylinder1}. This is needed in order to properly approximate the geometry of the cylinder without resorting to an isoparametric description of the physical boundaries as forwarded in \cite{TavelliSWE2014}.

\begin{figure}[!htbp]
	\begin{center}
		\begin{tabular}{cc} 
			\includegraphics[width=0.47\textwidth]{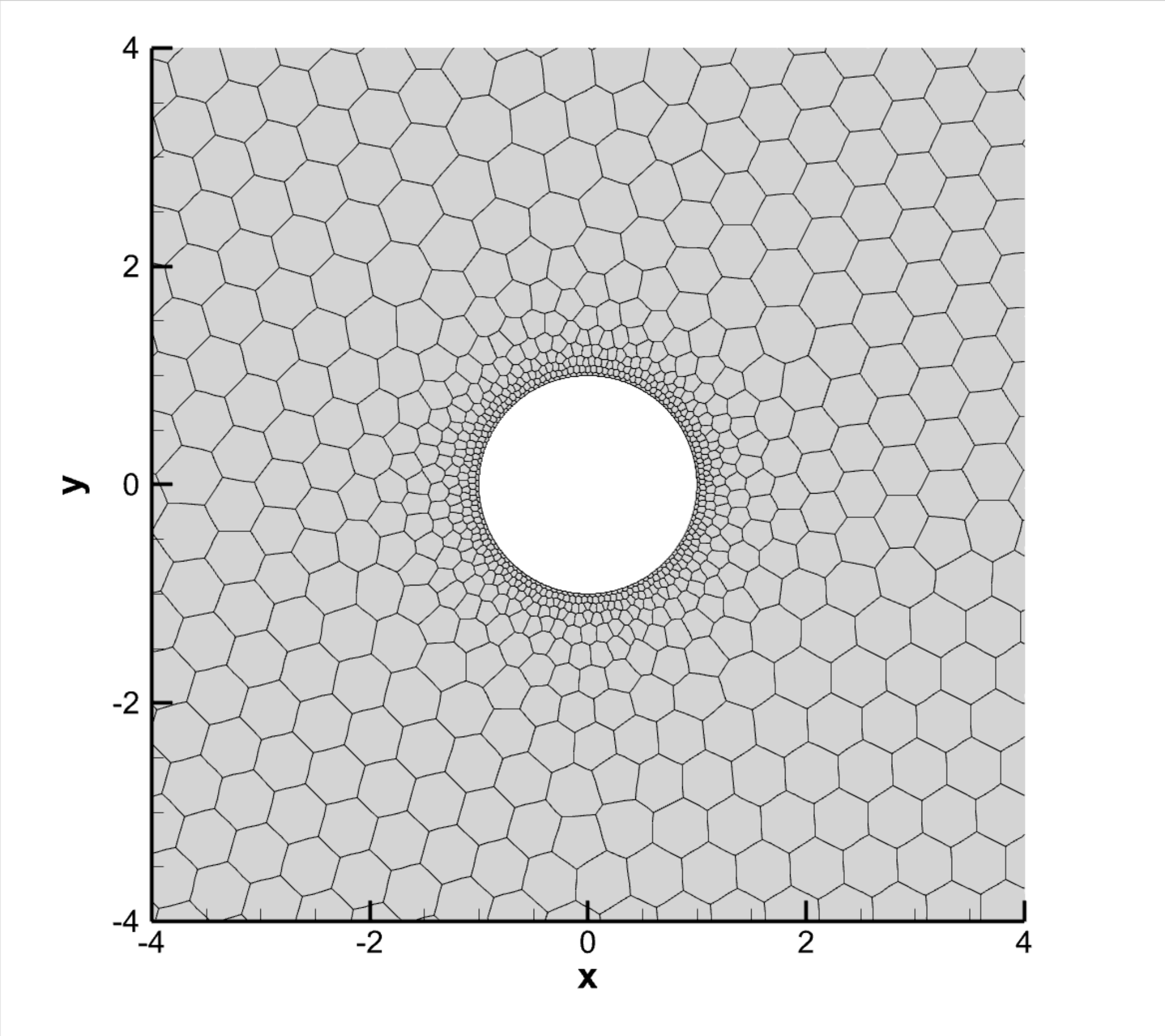} &  			
			\includegraphics[width=0.47\textwidth]{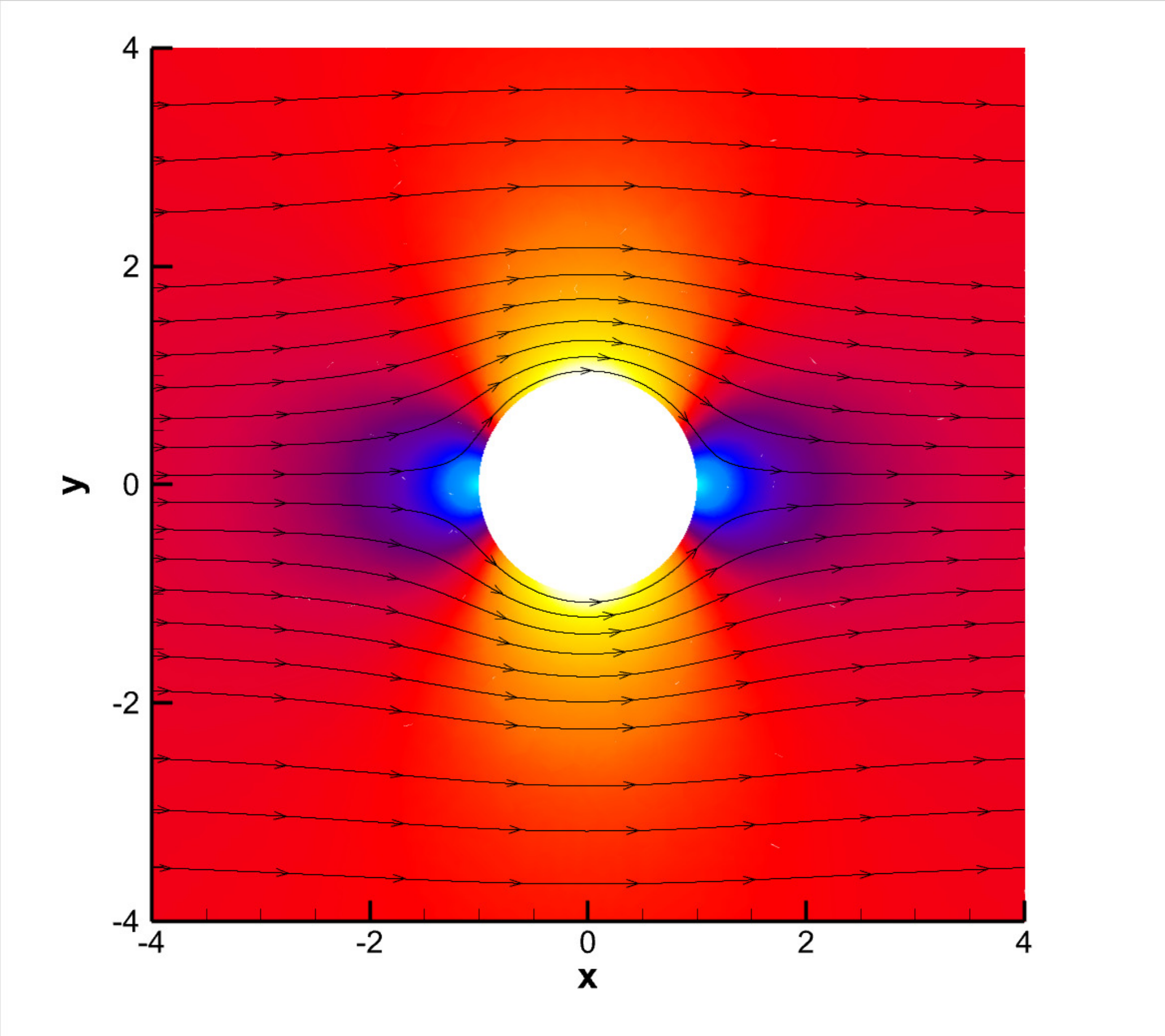} \\
		\end{tabular}
		\caption{Flow around a circular cylinder. Left: Voronoi computational mesh around the cylinder. Right: contours of the magnitude of the velocity field with associated streamtraces computed with the fourth order accurate SI-FVDG scheme.}
		\label{fig.Cylinder1}
	\end{center}
\end{figure}

The analytical solution for this test problem can be derived both for the velocity field in terms of polar coordinates $(r,\theta)$ as well as for the free surface elevation:
\begin{equation}
	\label{eqn.Cyl_ini}
	v_r = v_m \left( 1-\frac{r_c^2}{r^2} \right) \cos(\theta), \qquad v_{\theta} = -v_m \left( 1+\frac{r_c^2}{r^2} \right) \sin(\theta), \qquad \eta = \eta_0 + \frac{1}{2} v_m^2 g \left( 2\frac{r_c^2}{r^2} \cos(2\theta) -\frac{r_c^4}{r^2} \right),
\end{equation}
where we set $\eta_0=1$ and $v_m=10^{-2}$. In order to avoid the generation of strong initial transient waves we impose as initial condition the exact velocity field, but a flat free surface, namely $\eta(\xx,0)=\eta_0$. The exact solution is imposed on all boundaries, apart from the rightmost side of the domain ($x=16$) where an outflow condition is set. The simulation is run until the final time $t_f=10$, so that the stationary state has been reached. To enhance the advantages of the high order discretization proposed in this work, we run the SI-FVDG scheme using the second and the fourth order version in space, while keeping a first order time discretization for the sake of comparison. Figure \ref{fig.Cylinder1} shows the computational mesh around the cylinder as well as the magnitude of the velocity field with the associated streamlines at the final time. A comparison against the reference solution is plot in Figure \ref{fig.Cylinder2} along the circumference of radius $r = 1.01$ centered at the origin, where the fourth order accurate scheme retrieves much better the exact profile of the free surface elevation. The velocity field is resolved rather well by both schemes because the initial condition already provides the exact solution according to \eqref{eqn.Cyl_ini}.

\begin{figure}[!htbp]
	\begin{center}
		\begin{tabular}{cc} 
			\includegraphics[width=0.47\textwidth]{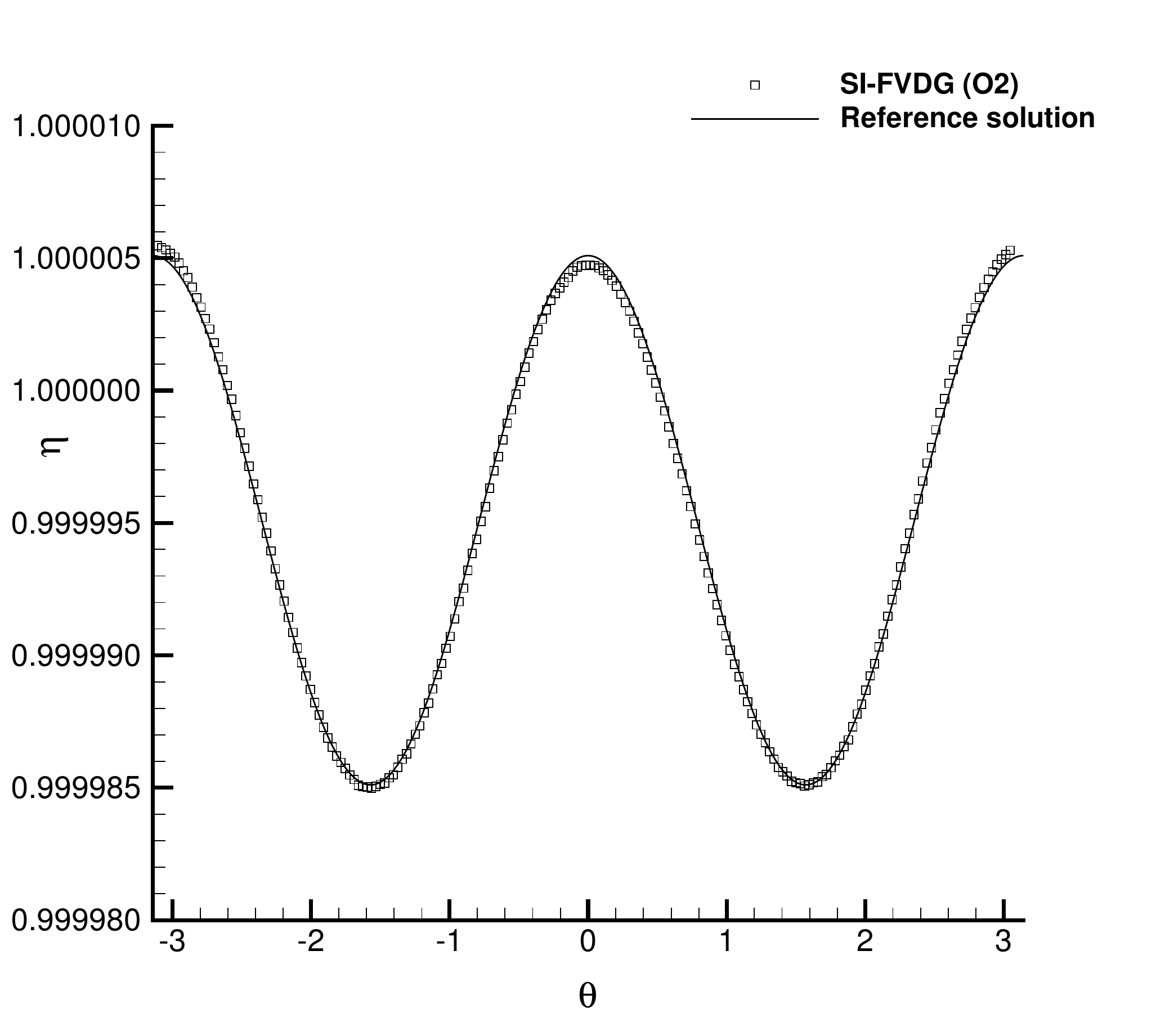} &  			
			\includegraphics[width=0.47\textwidth]{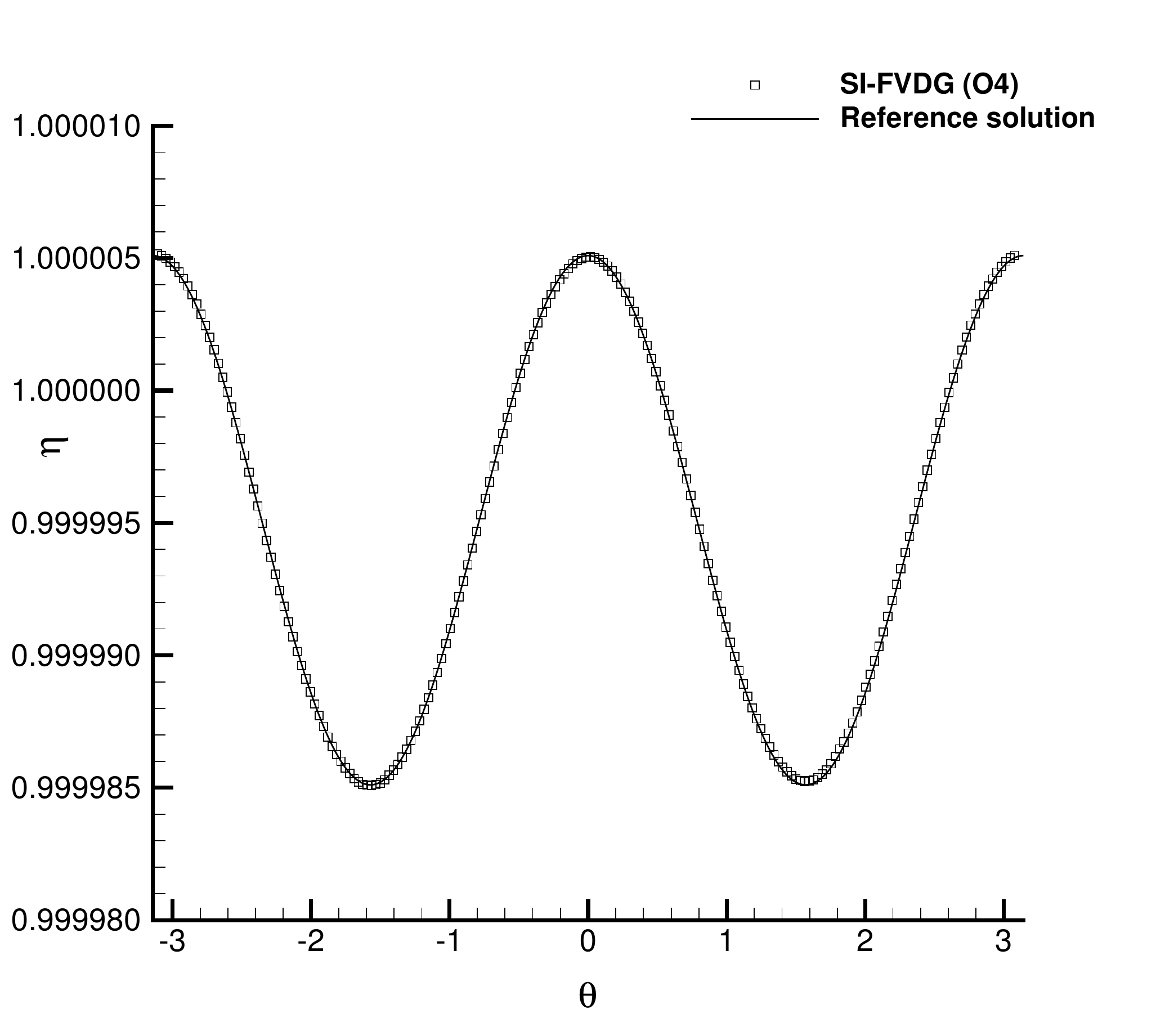} \\
			\includegraphics[width=0.47\textwidth]{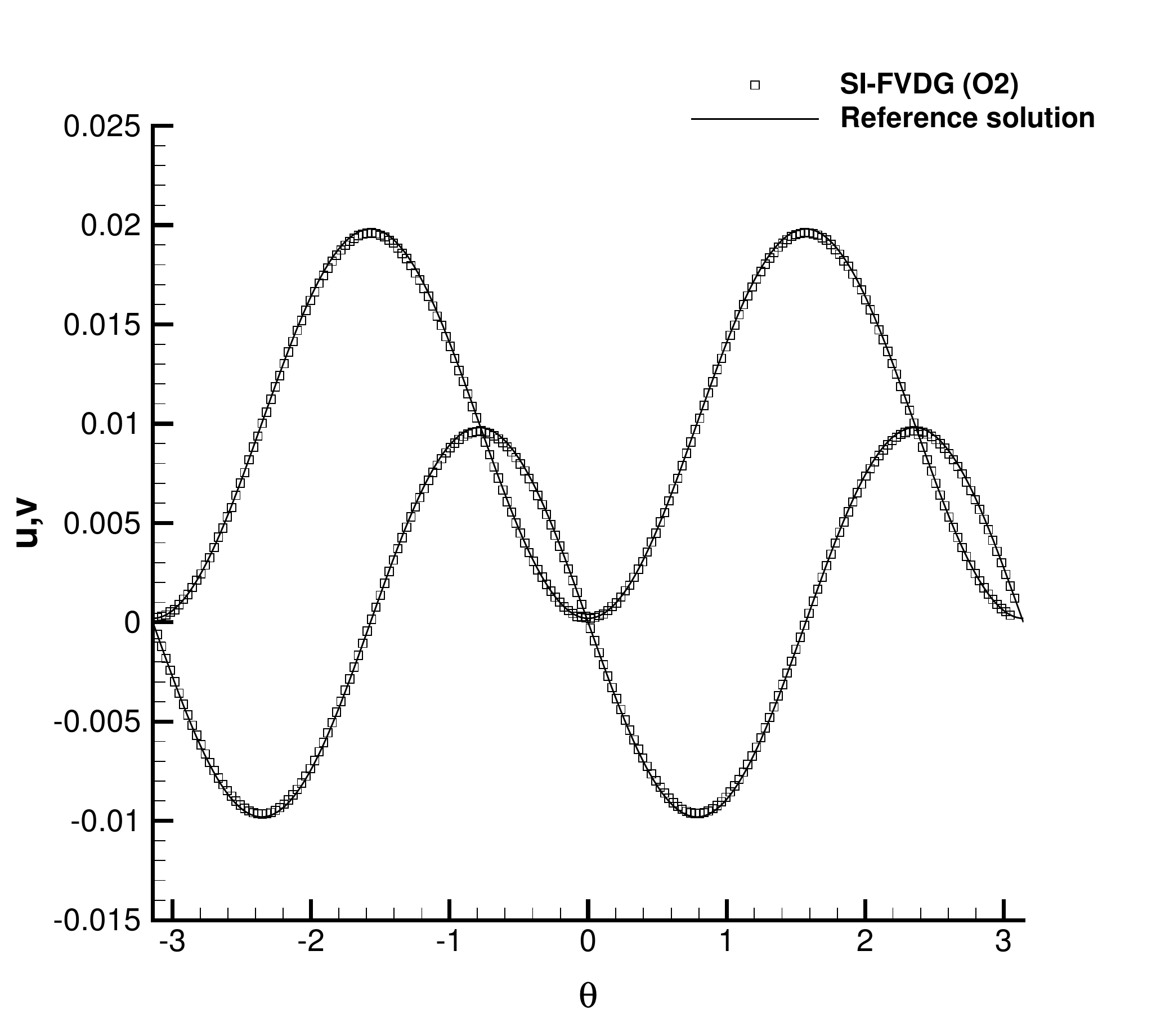} &  			
			\includegraphics[width=0.47\textwidth]{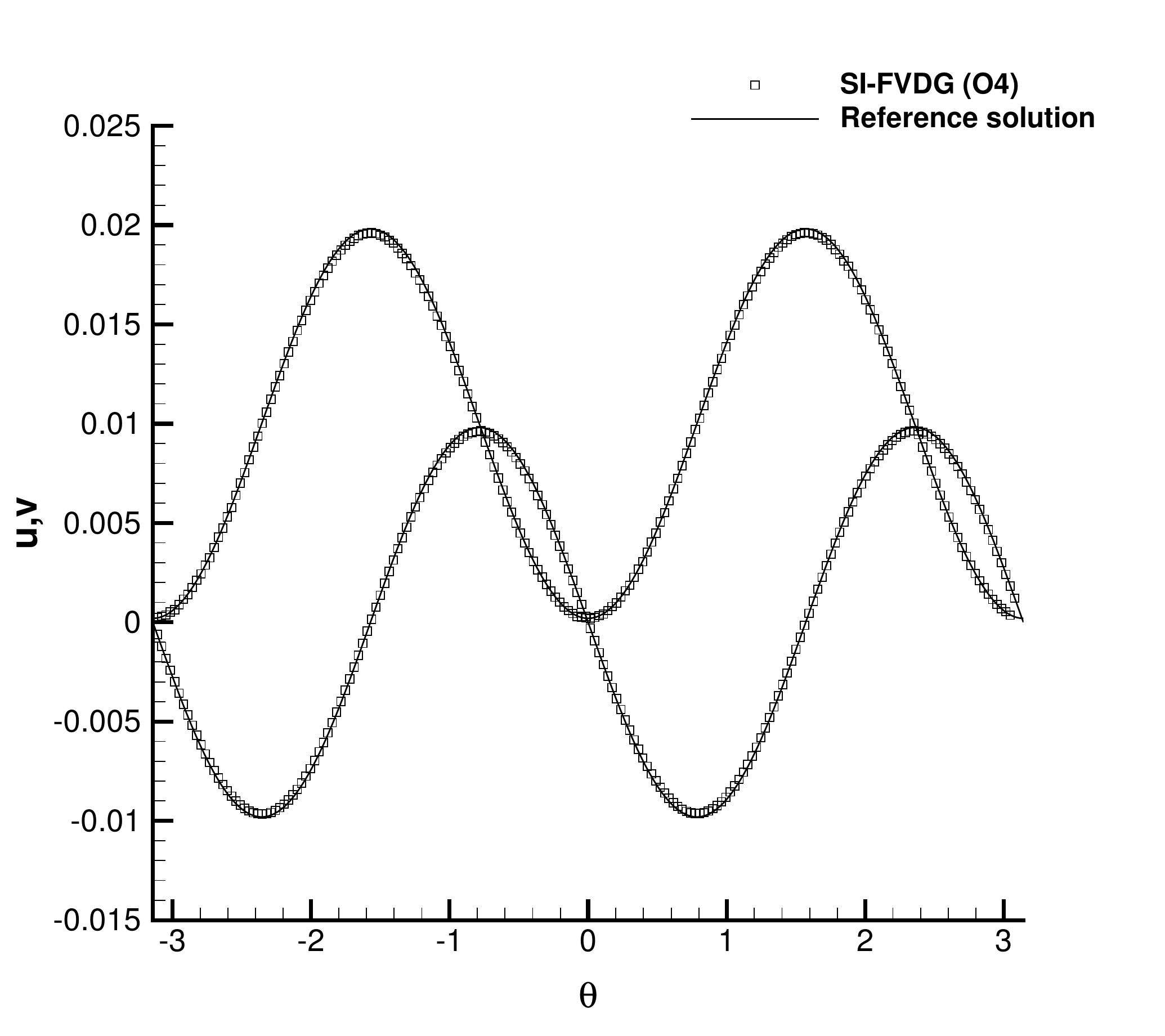} \\
		\end{tabular}
		\caption{Flow around a circular cylinder. Comparison of the numerical and exact solutions obtained for the free surface (top) and velocity field (bottom) at radius $r = 1.01$ using the second (left) and fourth (right) order SI-FVDG scheme.}
		\label{fig.Cylinder2}
	\end{center}
\end{figure}

\section{Conclusions} \label{sec.concl}
In this work we presented a high order all Froude regime IMEX well-balanced scheme for the two-dimensional shallow water model on unstructured polygonal meshes. In order to implement an accurate and efficient numerical scheme we combined a flux splitting formulation with an implicit-explicit discretization for the acoustic and advection waves, respectively. A high order numerical solution in space and time is obtained by a CWENO spatial reconstruction and a IMEX Runge-Kutta time integrator. Well-balanced and asymptotic preserving properties of the first order semi-discrete scheme have been demonstrated. The accuracy and robustness of the new proposed numerical scheme have been validated by solving six test problems. The first test problem deals with a convergence study where second and third order accuracy was reached considering four Froude regimes ($\Fr=0.32$,  $\Fr=10^{-2}$, $\Fr=10^{-4}$ and  $\Fr=10^{-6}$). The second test problem considered two scenarios, the first one testing the well-balanced property with a variable bathymetry and a lake at rest condition. This test shows that the scheme is capable to preserve the initial condition up to machine accuracy. The second scenario introduced a perturbation of the free surface producing a wave traveling over a variable bathymetry without generating any spurious artifact. The third numerical test problem simulates a circular dambreak case with a step on the bathymetry. The numerical solution agrees well with the reference solution. The fourth test problem deals with four Riemann problems and exact numerical solution, all of them solved accurately. The fifth one is a two-dimensional wave propagation initiated by a smooth perturbation of the free surface. In this test we observe that the third order numerical method is in very good agreement with the reference solution. Finally, in the sixth test problem we see a low Froude flow across a cylinder and the second and fourth order numerical solution matching the exact solution in an unstructured grid.

In the future we plan to apply and extend the novel schemes to the incompressible Navier-Stokes (INS) equations, since the wave equation for the pressure would look very similar to the one solved for the shallow water model, with the only difference lying in the metric term $H$ that will simply become a unity constant for the INS equations \cite{TavelliIncNS}. Further investigations will be devoted to treat also compressible viscous flows along the lines of \cite{SICNS22}, including an implicit discretization of the viscous terms. Finally, the inclusion of a mobile bottom bathymetry would require the coupling of the shallow water equations with the Exner equation, that also represents an interesting research direction.

\section*{Acknowledgments}

WB acknowledges financial support through from PRIN Project 2017 No. 2017KKJP4X granted by the Italian Ministry of Instruction, University and Research (MIUR). WB and MT are members of the GNCS-INdAM (\textit{Istituto Nazionale di Alta Matematica}) group. 

\newpage
\appendix

\section{IMEX schemes} \label{app.IMEX}
The Butcher tableau for the IMEX schemes used in this work are reported hereafter. They have been derived in \cite{PR_IMEX,PR_IMEXHO} and each IMEX scheme is described with a triplet $(s,\tilde{s},p)$ which characterizes the number $s$ of stages of the implicit method, the number $\tilde{s}$ of stages of the explicit method and the order $p$ of the resulting scheme. The acronym SA stands for Stiffly Accurate, while DIRK refers to Diagonally Implicit Runge-Kutta schemes.

\begin{itemize}
	\item SP(1,1,1)
	
	\begin{equation}
		\begin{array}{c|c}
			0 & 0 \\ \hline & 1
		\end{array} \qquad
		\begin{array}{c|c}
			1 & 1 \\ \hline & 1
		\end{array}
		\label{eqn.IMEX1}
	\end{equation}
	
	\item LSDIRK2(2,2,2) \hspace{0.2cm} $\gamma=1-1/\sqrt{2}$, \hspace{0.05cm} $\beta=1/(2\gamma)$
	
	\begin{equation}
		\begin{array}{c|cc}
			0 & 0 & 0 \\ \beta & \beta & 0 \\ \hline & 1-\gamma & \gamma
		\end{array} \qquad
		\begin{array}{c|cc}
			\gamma & \gamma & 0 \\ 1 & 1-\gamma & \gamma \\ \hline & 1-\gamma & \gamma
		\end{array}
		\label{eqn.IMEX2}
	\end{equation}
	
	\item SA DIRK (3,4,3) \hspace{0.2cm} $\gamma=0.435866$
	
	\begin{equation}
		\begin{array}{c|cccc}
			0 & 0 & 0 & 0 & 0 \\ \gamma & \gamma & 0 & 0 & 0 \\ 0.717933 & 1.437745 & -0.719812 & 0 & 0 \\ 1 & 0.916993 & 1/2 & -0.416993 & 0 \\ \hline  & 0 & 1.208496 & -0.644363 & \gamma
		\end{array} \qquad
		\begin{array}{c|cccc}
			\gamma & \gamma & 0 & 0 & 0 \\ \gamma & 0 & \gamma & 0 & 0 \\ 0.717933 & 0 & 0.282066 & \gamma & 0  \\ 1 & 0 & 1.208496 & -0.644363 & \gamma \\ \hline  & 0 & 1.208496 & -0.644363 & \gamma
		\end{array}
		\label{eqn.IMEX3}
	\end{equation}

\end{itemize} 
	
\section{CWENO reconstruction} \label{app.CWENO}
The piecewise reconstruction polynomials of degree $M$ have a total number of {unknown} degrees of freedom $\mathcal{M}=(M+1)(M+2)/2$ which are determined for each variable of the state vector $\U$ starting from the known cell averages $\U_i^n$. Let us consider a central reconstruction stencil $\mathcal{S}_i^c$ that is composed by the cell under consideration and by all the associated Neumann neighbors, hence
\begin{equation}
	\mathcal{S}_i^c = \bigcup \limits_{l=1}^{n_e} P_{j(l)}, 
	\label{stencil}
\end{equation}
where $j=j(l)$ denotes a mapping from the set of integers $l\in[1,n_e]$ to the global indexes $j$ used to sort the cells in the mesh. We assume that $j(1)=i$ so that the first cell in the  stencil is always the element for which we are computing the reconstruction. To avoid ill-conditioning of the resulting reconstruction matrices, the stencil contains a total number of elements $n_e$ that is greater than the smallest number $\mathcal{M}$ needed to reach the formal second order of accuracy (see \cite{barthlsq}).

The reconstruction polynomial $\w_i^{opt}(\xx)$ for the central stencil $\mathcal{S}_i^c$ is called \textit{optimal} polynomial and is expressed through the following conservative expansion
\begin{equation}
	\label{eqn.recpolydef} 
	\w_i^{opt}(\xx) = \sum \limits_{l=1}^\mathcal{M} 
	\beta_l^{(i)}(\xx) \, {\hat \w}^{opt}_{l,i} ,   
\end{equation}
with ${\hat \w}^{opt}_{l,i}$ representing the unknown expansion coefficients and the basis functions given by \eqref{eqn.Voronoi_modal}. The reconstruction procedure is built upon conservation on each element $P_j \in \mathcal{S}_i^c$, hence yielding an overdetermined linear system that is solved with a least-squares approach \cite{Dumbser2007693}, which reads
\begin{equation}
	\label{CWENO:Popt}
	\mathbf{w}_i^{opt} = \underset{{\mathbf{w}_i\in\mathcal{W}_i}}{\argmin}  
	\sum_{P_j \in \mathcal{S}_i} 
	\left( \U_j^n-\frac{1}{|P_{j}|} \int_{P_{j}} \mathbf{w_i}(\xx) \, d\xx \right)^2,
\end{equation}
where $\mathcal{W}_i $ is the set of all polynomials $\mathbb{P}_M$ satisfying
\begin{equation}\label{CWENO:Popt1}
	\mathcal{W}_i = \left\{\mathbf{w}_i \in \mathbb{P}_M: \bar \U_j^n=\frac{1}{|P_{i}|} \int_{P_{i}} \mathbf{w_i}(\xx) \, d\xx \right\}
	\subset \mathbb{P}_M.
\end{equation}
The optimal polynomial $\mathbf{w}_i^{opt}$ is chosen among all the possible polynomials of degree $M$ so that it exhibits the property of sharing the same cell average of the finite volume data $\U_i^n$ in the cell $P_i$ while being close in the least-square sense to the other cell averages in the stencil $\mathcal{S}_i$.

The polynomial $\mathbf{w}_i^{opt}$ is generated from a linear arbitrary high order reconstruction procedure, thus it needs to be stabilized by a nonlinear operator, which will be done following the CWENO approach. A set of $N_{S_i}$ interpolating polynomials of degree one are also computed in order to make a nonlinear hybridization among the resulting polynomials. These are called \textit{lateral} reconstruction polynomials $\mathbf{w}_i^L$ that are obtained by considering one-sided stencils $\mathcal{S}_i^L$ always composed by three elements, namely the element itself $P_i$, one direct neighbor $P_j$ and the other Neumann neighbor that is a direct neighbor of both $P_i$ and $P_j$. For each stencil $\mathcal{S}_i^L$ with $L=1,\ldots,N_{S_i}$, the linear polynomial $\mathbf{w}_i^L$ is obtained through the unique solution of the system 
\begin{equation}
	\label{CWENO:Ps}
	\mathbf{w}_i^L\in\mathbb{P}_1 
	\, \, \text{  s.t. } \, \, \forall P_{j} \in \mathcal{S}_i^L: \,\, 
	\bar \U_j^n=\frac{1}{|P_{j}|} \int_{P_{j}} \mathbf{w}_i^L(\xx) \, d\xx,
\end{equation}
where $j$ indicates the mesh element belonging to the stencil $\mathcal{S}_i^L$ and the polynomial $\mathbf{w}_i^L$ is defined again relying on the same conservative Taylor expansion \eqref{eqn.Voronoi_modal}. The central polynomial corresponding to $L=0$ is then derived on the basis of conservation principles as
\begin{equation}
	\label{CWENO:P0}
	\mathbf{w}_i^0= \frac{1}{\delta_0} \, \mathbf{w}_{i}^{opt} - \sum_{L=1}^{N_{S_i}} \frac{\delta_0}{\delta_{L}} \mathbf{w}^L_{i} \in\mathbb{P}_1,
\end{equation}
where $\delta_{0},\ldots,\delta_{N_{S_i}}$ are positive coefficients such that
\begin{equation}
	\sum_{s=0}^{N_{S_i}}\delta_{s}=1.
	\label{eqn.sumCWENO}
\end{equation}
A linear combination of the polynomials $\mathbf{w}_i^0,\ldots,\mathbf{w}_i^{N_{S_i}}$ with the linear weights $\delta_{0},\ldots,\delta_{N_{S_i}}$ yields the optimal polynomial $\mathbf{w}_i^{opt}$. In this way the accuracy of the CWENO reconstruction does not depend on the choice of the coefficients, which must only represent a normalization that sums up to unity. In order to achieve essentially non-oscillatory properties, the final CWENO reconstruction polynomial is computed from the reconstruction polynomials obtained on each single stencil. Therefore, the following oscillation indicators $\sigma_s$ are introduced
\begin{equation}
	\sigma_s = \sum \limits_{l=2}^{\mathcal{M}} \left(\hat{\w}^{s}_{l,i}\right)^2,
	\label{eqn.OI}
\end{equation}
where $\hat{\w}^{s}_{l,i}$ denote the expansion coefficients \eqref{eqn.recpolydef} of the polynomial defined on stencil $s$. The nonlinear weights $\omega_s$ are then given by
\begin{equation}
	\label{eqn.weights}
	\omega_s = \frac{\tilde{\omega}_s}{\sum \limits_{s=0}^{N_{S_i}} \tilde{\omega}_s}, 
	\qquad \textnormal{ with } \qquad 
	\tilde{\omega}_s = \frac{\delta_s}{\left(\sigma_s + \epsilon \right)^r}, 
\end{equation} 
where $\epsilon=10^{-14}$ and $r=4$ are chosen according to \cite{Dumbser2007693}. Furthermore, we set $\delta_0=200/\delta_{sum}$ and $\delta_{L}=1/\delta_{sum}$ with $\delta_{sum}=200+N_{S_i}$ for the definition of the positive coefficients. The final nonlinear CWENO reconstruction polynomial and its coefficients are then given by 
\begin{eqnarray}
	\label{CWENO:Prec}
	\mathbf{w}(\xx,t^n)&=& \sum_{s=0}^{N_{S_i}} \omega_s \mathbf{w}_{i}^s(\xx) \nonumber \\
	& =& \beta_l^{(i)} \, \hat{\w}_{l,i}^n.
\end{eqnarray}	
The reconstruction operator $\Rop$ given by \eqref{eqn.Rop} carries out the CWENO procedure detailed above and provides the sought high order expansion coefficients $\hat{\w}_{l,i}^n$. Further details can be found in \cite{ArepoTN,FVBoltz}.

	%
	\bibliographystyle{plain}
	\bibliography{biblio}

\end{document}